\pdfoutput=1
\RequirePackage{silence}
\documentclass[a4paper,11pt]{amsart}
\usepackage[hmarginratio={1:1},vmarginratio={1:1},lmargin=60.0pt,tmargin=60.0pt]{geometry}

\allowdisplaybreaks

\usepackage[utf8]{inputenc}

\usepackage{latexsym,exscale,mathtools,textcomp}
\usepackage{amssymb,amsmath,amsthm,amsfonts,mathrsfs,bbm,enumitem,stmaryrd,dsfont}
\usepackage[table]{xcolor}
\usepackage{graphicx}
\usepackage{mathtools}
\usepackage{ytableau}
\usepackage{relsize}
\usepackage{xfrac}
\usepackage{caption}
\usepackage{float}
\usepackage{xparse}
\usepackage{graphbox}
\usepackage{verbatim}

\usepackage{dynkin-diagrams}

\usepackage{array}
\newcolumntype{C}{>{$}c<{$}}

\definecolor{mygray}{gray}{0.6}
\definecolor{mygraydark}{gray}{0.4}
\definecolor{mygraylight}{gray}{0.85}
\definecolor{spinach}{RGB}{46,139,87}
\definecolor{tomato}{RGB}{255,99,71}
\definecolor{orchid}{RGB}{143,40,194}
\definecolor{neon}{RGB}{77,77,255}
\definecolor{pumpkin}{RGB}{224,180,80}
\definecolor{citron}{RGB}{190,180,90}

\definecolor{lava}{RGB}{207,16,32}
\definecolor{cream}{RGB}{255,253,208}
\definecolor{verdigris}{RGB}{67,179,174}
\definecolor{Black}{RGB}{0,0,0}
\definecolor{mydarkblue}{RGB}{10,10,170}
\definecolor{darkspinach}{RGB}{20,70,20}
\definecolor{darktomato}{RGB}{155,40,30}
\definecolor{darkorchid}{RGB}{50,10,100}
\definecolor{darklava}{RGB}{150,8,16}

\usepackage{todonotes}

\usepackage{enumitem}
\setlist[enumerate]{itemsep=0.15cm,label=\emph{\upshape(\alph*)}}
\setlist[enumerate,2]{itemsep=0.15cm,label=\emph{\upshape(\roman*)}}
\setlist[enumerate,3]{itemsep=0.15cm,label=\emph{\upshape(\Alph*)}}

\let\emph\relax
\DeclareTextFontCommand{\emph}{\bfseries\em}

\newcommand{\placeholder}{{}_{-}}

\renewcommand{\dots}{\text{...}}

\let\<=\langle
\let\>=\rangle

\newcommand{\acts}{\centerdot}

\renewcommand{\dots}{\text{...}}

\DeclarePairedDelimiterX{\set}[1]{\{}{\}}{\setargs{#1}}
\NewDocumentCommand{\setargs}{>{\SplitArgument{1}{|}}m}{\setargsaux#1}
\NewDocumentCommand{\setargsaux}{mm}
{\IfNoValueTF{#2}{#1} {#1\,\delimsize|\,\mathopen{}#2}}

\newcommand{\ie}{\text{i.e.}}
\newcommand{\eg}{\text{e.g.}}
\newcommand{\cf}{\text{cf.}}
\newcommand{\etc}{\text{etc.}}
\newcommand{\aka}{\text{a.k.a.}}

\newcommand{\N}{\mathbb{Z}_{\geq 0}}

\newcommand{\K}{\mathbb{K}}

\newcommand{\Hom}{\mathrm{\Hom}}

\newcommand{\algebra}[1][A]{\mathscr{#1}}

\newcommand{\cellbasis}[1][\mathscr{A}]{B_{\mathscr{A}}}

\newcommand{\rad}[1][\jcell]{\mathrm{Rad}}

\newcommand{\lcell}{\mathcal{L}}
\newcommand{\rcell}{\mathcal{R}}
\newcommand{\jcell}{\mathcal{J}}
\newcommand{\hcell}{\mathcal{H}}

\newcommand{\jideal}[1][\lambda]{\algebra^{>_{lr}\lambda}}

\newcommand{\module}[1][M]{{#1}}
\newcommand{\one}{\mathbbm{1}}
\newcommand{\oneb}{\mathbbm{1}_{b}}
\newcommand{\onet}{\mathbbm{1}_{t}}

\newcommand{\jb}{\mathcal{J}_{b}}
\newcommand{\jt}{\mathcal{J}_{t}}

\newcommand{\monoid}[1][M]{\mathrm{#1}}

\newcommand{\onemon}{\mathrm{1}}

\newcommand{\group}[1][G]{\mathrm{#1}}

\usepackage[all]{xy}
\usepackage{tikz}
\usetikzlibrary{cd}
\usetikzlibrary{decorations}
\usetikzlibrary{decorations.markings}
\usetikzlibrary{decorations.pathreplacing}
\usetikzlibrary{decorations.pathmorphing}
\usetikzlibrary{arrows.meta,shapes,positioning,matrix,calc}
\usetikzlibrary{shapes.callouts}

\tikzset{
anchorbase/.style={baseline={([yshift=#1]current bounding box.center)}},
anchorbase/.default={-0.5ex},
tinynodes/.style={font=\tiny,text height=0.25ex,text depth=0.05ex},
smallnodes/.style={font=\scriptsize,text height=0.75ex,text depth=0.15ex},
mor/.style={line width=0.75,color=black,fill=cream},
mor2/.style={line width=0.75,color=black,fill=tomato},
mor3/.style={line width=0.75,color=black,fill=spinach},
usual/.style={line width=1.2,color=black},
crossline/.style={preaction={draw=white,line width=5.0pt,-},preaction={draw=black,line width=0.9pt,-}},
dot/.style = {
decoration={markings,
post length=0.25mm,
pre length=0.25mm,
mark=at position #1 with {\node[circle,radius=0.15cm,inner sep=-1.2pt,color=black,fill=black]{};}
},
postaction={decorate}
},
dot/.default=1,
}
\tikzstyle directed=[postaction={decorate,decoration={markings,
mark=at position #1 with {\arrow[line width=0.25mm, black]{>}}}}]

\usepackage{aliascnt,etoolbox}
\def\NewTheorem#1{%
\newaliascnt{#1}{equation}%
\newtheorem{#1}[#1]{#1}%
\aliascntresetthe{#1}%
\expandafter\def\csname #1autorefname\endcsname{#1}%
}
\def\equationautorefname~#1\null{(#1)\null}

\numberwithin{equation}{subsection}

\NewTheorem{Proposition}
\NewTheorem{Theorem}
\NewTheorem{Corollary}
\AtEndEnvironment{Corollary}{\null\hfill$\square$}%
\NewTheorem{Lemma}
\NewTheorem{Conjecture}
\NewTheorem{Speculation}
\NewTheorem{Observation}
\NewTheorem{Assumption}
\theoremstyle{definition}
\NewTheorem{Definition}
\AtEndEnvironment{Definition}{\null\hfill$\Diamond$}%
\NewTheorem{Classification Problem}
\AtEndEnvironment{Classification Problem}{\null\hfill$\Diamond$}%
\NewTheorem{Notation}
\AtEndEnvironment{Notation}{\null\hfill$\Diamond$}%
\NewTheorem{Example}
\AtEndEnvironment{Example}{\null\hfill$\Diamond$}%
\NewTheorem{Examples}
\AtEndEnvironment{Examples}{\vskip-10mm\null\hfill$\Diamond$}%

\theoremstyle{remark}
\NewTheorem{Remark}
\AtEndEnvironment{Remark}{\null\hfill$\Diamond$}%
\NewTheorem{Question}
\AtEndEnvironment{Question}{\null\hfill$\Diamond$}%

\usepackage[hypertexnames=false]{hyperref}
\usepackage{bookmark}
\hypersetup{
pdftoolbar=true,
pdfmenubar=true,
pdffitwindow=false,
pdfstartview={FitH},
pdftitle={Representation gaps of rigid planar diagram monoids},
pdfauthor={Willow Stewart and Daniel Tubbenhauer},
pdfsubject={},
pdfcreator={Willow Stewart and Daniel Tubbenhauer},
pdfproducer={Willow Stewart and Daniel Tubbenhauer},
pdfkeywords={},
pdfnewwindow=true,
colorlinks=true,
linkcolor=mydarkblue,
citecolor=teal,
filecolor=magenta,
urlcolor=orchid,
linkbordercolor=lava,
citebordercolor=teal,
urlbordercolor=orchid,
linktocpage=true
}

\def\makeautorefname#1#2{\csdef{#1autorefname}{#2}}

\makeautorefname{section}{Section}%
\makeautorefname{subsection}{Section}%
\makeautorefname{subsubsection}{Section}%

\title[Representation gaps of rigid planar diagram monoids]{Representation gaps of rigid planar diagram monoids}
\author[W. Stewart and D. Tubbenhauer]{Willow Stewart and Daniel Tubbenhauer}

\address{W.S.: The University of Sydney, School of Mathematics and Statistics F07, Office Carslaw 807, NSW 2006, Australia, \href{https://www.maths.usyd.edu.au/ut/people?who=W_Stewart}{www.maths.usyd.edu.au/ut/people?who=W\_Stewart}, \href{https://orcid.org/0009-0000-2854-2256}{ORCID: 0009-0000-2854-2256}}
\email{W.Stewart@maths.usyd.edu.au}
\address{D.T.: The University of Sydney, School of Mathematics and Statistics F07, Office Carslaw 827, NSW 2006, Australia, \href{http://www.dtubbenhauer.com}{www.dtubbenhauer.com}, \href{https://orcid.org/0000-0001-7265-5047}{ORCID 0000-0001-7265-5047}}
\email{daniel.tubbenhauer@sydney.edu.au}

\begin{document}

\begin{abstract} 
We define non-pivotal analogs of the Temperley--Lieb, Motzkin, and planar rook monoids, and compute bounds for the sizes of their nontrivial simple representations. From this, we assess the two types of monoids in their relative suitability for use in cryptography by comparing their representation gaps and gap ratios. We conclude that the non-pivotal monoids are generally worse for cryptographic purposes.
\end{abstract}


\maketitle

\tableofcontents

\section{Introduction}

The main question of this paper is how representations behave when a pivotal structure is relaxed to a rigid, yet non-pivotal, one. Our primary objects of study are (partially) novel diagram categories that serve as non-pivotal analogs of classical diagram categories, such as the Temperley–-Lieb category.

\subsection{Setting the stage: diagram monoids and representation gaps}

Motivated by the problem of developing cryptographic protocols based on finite noncommutative monoids, \cite{khovanov-monoidal-2024} developed a framework to connect monoidal categories and cryptography. Noncommutative monoids are often vulnerable to \textit{linear attacks}: an attacker can convert a protocol operating in such a monoid into a linear algebra problem, enabling the use of powerful linear algebra techniques to break the encryption; see e.g. \cite{MyRo-linear-attack} for examples. Converting such a protocol to linear algebra requires a nontrivial, say simple or faithful, representation of the monoid. In short, a cryptographic protocol based on a noncommutative monoid can be bypassed by targeting the nontrivial simple representations, so if these are too low in dimension, the protocol will not be secure. The security of such a protocol can be roughly measured by the \emph{representation gap (RepGap)}. As such, it is important to find monoids with large RepGaps. There are also many other reasons to study these numerical invariants of monoids, see, for example, \cite{BoGa-bounds-cayley-graphs-slfp,Go-quasirandom-groups}.

Monoidal categories provide a wealth of such monoids, since each object $X$ provides an endomorphism monoid $\textnormal{End}(X^{\otimes n})$ for each $n \in \mathbb{Z}_{\geq 0}$, and by the Schur--Weyl dual of \cite{CoOsTu-growth} it is expected that their representations grow exponentially. In particular, \cite[Section 4 \& 5]{khovanov-monoidal-2024} discusses \textit{diagram monoids}, which arise from endomorphisms of subcategories of the partition category. These are split into the \textit{symmetric} monoids, namely the partition monoid, rook Brauer monoid, Brauer monoid, rook monoid, and symmetric monoid (group), and their \textit{planar} counterparts, the planar partition monoid, Motzkin monoid, Temperley--Lieb monoid, planar rook monoid, and planar symmetric monoid (the trivial group). The main tools for computing the representations of these monoids and their dimensions are based on \textit{cell theory (a.k.a. Green's relations)} as in \cite{Gr-structure-semigroups}. 

The paper \cite{khovanov-monoidal-2024} primarily explores planar monoids, noting that symmetric monoids may be less suitable for cryptographic applications than their planar counterparts. While this hypothesis remains unproven, it raises the question of how categorical structures, such as symmetry, in a monoidal category influence the RepGap (see e.g. \cite{egno-tensor-2015} or, more diagrammatic, \cite{Tu-qt} for background on such structures).
This is where our story starts:
The monoidal categories above (excluding the rook categories and those for the groups) are \textit{pivotal}, meaning every object $X$ has a left and right dual, which essentially means that the category is \textit{rigid}, and further the double dual of $X$ is isomorphic to $X$. 
A natural generalization is to consider \textit{non-pivotal}, rigid monoidal categories as a potential source of monoids for cryptographic purposes. These are, in a sense, more ``complicated'' than pivotal monoidal categories: For example, the Temperley--Lieb category is equivalent to the category $Rep(SL_2)$ of finite dimensional representations of the special linear group (a classic result; see \cite{RuTeWe-sl2}) or quantum versions of the group $SL_2$, and Schur--Weyl duality arguments \cite{erd,soe5,AnStTu-cellular-tilting} then imply that the RepGap problem can be solved using representation theory of $SL_2$.
No such equivalence exists for any non-pivotal, rigid monoidal category and finite dimensional representations of any group $G$. 

\begin{Remark}
To avoid confusion, note that in the classical diagram categories under discussion, the object $X$ is self-dual. Relaxing this assumption, while retaining pivotality, would, for example, require introducing an orientation. However, this is \textit{not} what we do here: we break the symmetry of left and right duals, and, as far as we can tell, the categories and monoids we study are (partially) new. Importantly, we expect that passing to oriented versions of the classical diagram categories does not affect the RepGap. For a related result, see \cite{GT} for the Schur--Weyl dual statement in the context of the Brauer category.
\end{Remark}

This paper begins the study of this with non-pivotal, rigid analogs of some of the planar diagram monoids above, namely the Temperley--Lieb monoid, Motzkin monoid, and the planar rook monoid. The rigid Temperley--Lieb category, as we define it, has been studied as the free rigid monoidal category, for example in \cite{adjoin-dual-2021,St-honours-2021,Tu-qt}. As far as we are aware, there is no existing definition for a rigid analog for the Motzkin monoid. Since the planar rook category has no duals, we treat it as a two-color submonoid of the rigid Motzkin monoid. We directly ignore the planar partition monoid and the planar symmetric monoid, as the former is isomorphic to the even-strand Temperley--Lieb monoids, see \cite{HaRa-partition-algebras}, and the latter is isomorphic to the trivial group.  

\subsection{Results}

\begin{figure}[ht]
\begin{gather*}
\includegraphics[scale=0.6]{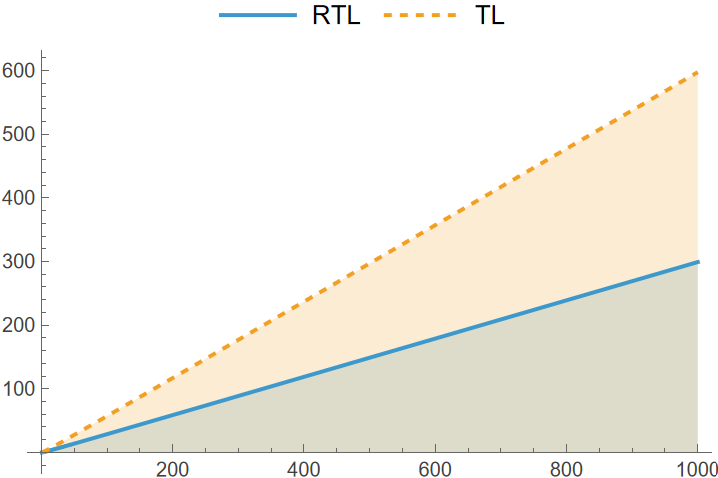}
\includegraphics[scale=0.6]{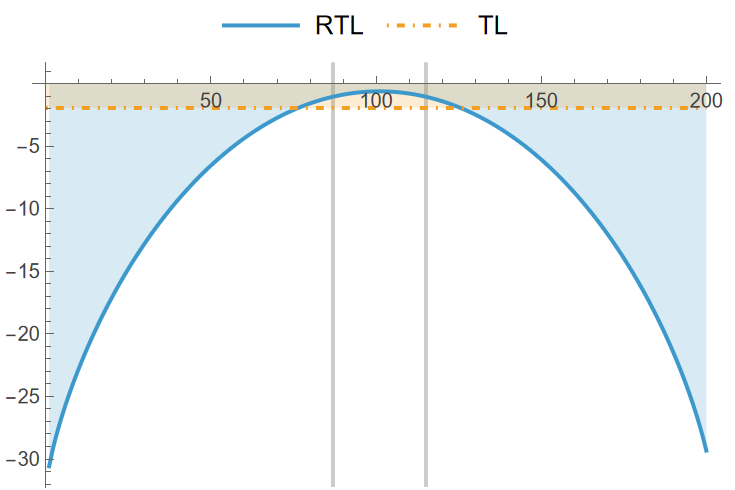}
\end{gather*}
\caption{Comparison of the Temperley--Lieb and rigid Temperley--Lieb monoids
on base-ten logarithmic scales. Left, representation gaps up to $n=1000$,
using the exact gap for the Temperley--Lieb monoid. Right, normalized simple
dimensions by rank at $n=100$. The vertical lines mark the truncation
window.}
\label{Fig:Intro}
\end{figure}

One might assume that a more complicated monoidal category would translate to more secure cryptographic protocols; however, \emph{surprisingly}, we find that the RepGap and gap ratio (the RepGap divided by the square root of the order of the truncated monoid) are largely worse for the non-pivotal diagram monoids. The asymptotic behavior of these is summarized in the table below (more precise results are given in the main body of the paper), with the pivotal monoids on the left and the non-pivotal monoids on the right. We also colored the quantities that are larger when comparing between non-pivotal and pivotal. The table uses the truncations defined below. The rigid Motzkin entries are semisimple bounds.
\begin{gather*}
\begin{tabular}{c|c|c||c|c|c}
\arrayrulecolor{tomato}
Monoid & RepGap & Gap Ratio 
& Monoid & RepGap & Gap Ratio
\\
\hline
\hline
$TL_n^T$ & \cellcolor{spinach!50}$\Theta(n^{-3/2}4^n)$ & $\Theta(n^{-3/4})$ & $rTL_n^T$ & $\Theta(n^{-1/2}2^n)$ & \cellcolor{spinach!50}$\Theta(n^{-1/4})$
\\
\hline 
$Mo_n^T$ &\cellcolor{spinach!50}$\textnormal{Gap}^{1/n}\to 9$ & \cellcolor{spinach!50}$\cellcolor{orchid!50}\textnormal{Ratio}^{1/n}\to 1$ & $rMo_n^T$ & $\textnormal{ssGap}^{1/n}\to4$ & $\cellcolor{orchid!50}$$\textnormal{ssRatio}^{1/n}\to1$
\\
\hline 
$pRo_n^T$ & \cellcolor{spinach!50}$\Theta(n^{-1/2}4^n)$ & \cellcolor{spinach!50}$\Theta(n^{-1/4})$ & $rpRo_n^T$ & $\Theta(n^{-1/2}2^n)$ & $O(0.87^n)$
\\
\end{tabular}
\end{gather*}
For Temperley--Lieb, passing to the rigid truncation increases the gap ratio by a factor of order $n^{1/2}$ and decreases the RepGap by a factor of order $2^n/n$, as illustrated in \autoref{Fig:Intro}. For the Motzkin monoids our results are not fine enough to come to a conclusion for the gap ratio, but the RepGap is again much worse.

We ideally aim for a more secure protocol over improved relative security for the computational complexity, so in this case we still say that the pivotal Temperley--Lieb monoid is better for cryptographic purposes than the non-pivotal analog; ditto for the Motzkin monoids. For the remaining case, it is clear that the pivotal monoids are better. With this, we build on the results from \cite{khovanov-monoidal-2024} to determine that pivotal planar diagram monoids are generally better for cryptographic purposes than non-pivotal planar diagram monoids.

The results related to the symmetric monoids in \cite[Section 5]{khovanov-monoidal-2024} are much less precise than those of the planar monoids. This is because there is greater difficulty in computing the dimensions of simple representations for these monoids, although some work toward that has been done in, for example, \cite{dwh,cdm,Tu-sandwich}. Future work based on this paper should thus begin from these symmetric monoids, both in the pivotal and non-pivotal cases. 

\begin{Remark}
We assume some familiarity with diagrammatic algebras, monoids, and categories, and (sandwich) cellularity a.k.a. Green's relations. There are many references for these, see e.g. \cite[Sections 3 and 4]{khovanov-monoidal-2024}. (For readers acquainted with diagram algebras such as the Temperley--Lieb algebra, the monoids we study can be obtained by specializing all parameters, i.e., floating components to $1$.)
\end{Remark}

\begin{Remark}
Some parts of this paper are based on computer calculations. For the reader who wants to run these calculations themselves, we have collected all relevant material on GitHub following the link \cite{St-github-2025}. This applies to all Mathematica and GAP computations in this paper. That page also contains an erratum, which may be empty.
\end{Remark}

\noindent\textbf{Acknowledgments.}
This paper is intended to be part of the first author's PhD thesis.

We would like to thank Kevin Coulembier for bringing the rigid Temperley–Lieb category to our attention, and Matthias Fresacher for providing assistance with definitions and Gram matrices, and a referee for some useful feedback.
WS was supported by the Postgraduate Research Scholarship in Mathematics and Statistics (SC4238).
DT was supported by the ARC Future Fellowship FT230100489 and contemplates the various ways their life went awry, taking full responsibility.

\noindent\textbf{AI declaration.}
When this paper was written, we did not use any AI models. For the revision we used Astra 6.2 to polish the existing arguments.

\section{Cells and monoid RepGaps}

The main representation theory of monoids that we will be using is discussed in detail in \cite[Section 2 \& 3]{khovanov-monoidal-2024}, which we follow closely; see also \cite{St-rep-monoid}. In this section, we will simply state the main results and tools that we require without proof. 

\subsection{Monoid representations}

Let $\monoid$ be a finite monoid, and let $\mathbb{K}$ be a field. Except where specified, $\mathbb{K}$ has arbitrary characteristic and representations are left representations over $\mathbb{K}$. We define preorders
on $\monoid$, called \textit{left, right and two-sided cell order}, by
\begin{align*}
(a\leq_{l}b)&\Leftrightarrow
\exists c:b=ca,
\\
(a\leq_{r}b)&\Leftrightarrow
\exists c:b=ac,
\\
(a\leq_{lr}b)&\Leftrightarrow
\exists c,d:b=cad.
\end{align*}

\begin{Remark}\label{R:CellsOrder}
There are similar notions of cells in the study of cellular algebras as, e.g., in \cite{GrLe-cellular,KoXi-affine-cellular,TuVa-handlebody,Tu-sandwich}. 
As in \cite[Remark 3.1]{khovanov-monoidal-2024}, these orders are
in-line with the most common convention used in the theory of cellular algebras but the opposite of the one usually used in monoid theory.
\end{Remark}

We define equivalence relations, 
the \textit{left, right and two-sided equivalence}, by
\begin{align*}
(a\sim_{l}b)&\Leftrightarrow
(a\leq_{l}b\text{ and }b\leq_{l}a),
\\
(a\sim_{r}b)&\Leftrightarrow
(a\leq_{r}b\text{ and }b\leq_{r}a),
\\
(a\sim_{lr}b)&\Leftrightarrow
(a\leq_{lr}b\text{ and }b\leq_{lr}a).
\end{align*}
The respective equivalence classes are called 
left, right respectively two-sided \textit{cells}. 
We denote all these by $\lcell$, $\rcell$ and $\jcell$
and call two-sided cells \textit{$J$-cells}. 
Finally, an \textit{$H$-cell} $\hcell=\hcell(\lcell,\rcell)=\lcell\cap\rcell$ is an intersection of a 
left $\lcell$ and a right cell $\rcell$. 

The picture to keep in mind is
\begin{gather}\label{Eq:CellsIllustration}
\begin{tikzpicture}[baseline=(A.center),every node/.style=
{anchor=base,minimum width=1.4cm,minimum height=1cm}]
\matrix (A) [matrix of math nodes,ampersand replacement=\&] 
{
\hcell_{11} \& \hcell_{12} 
\& \hcell_{13} \& \hcell_{14}
\\
\hcell_{21} \& \hcell_{22} 
\& \hcell_{23} \& \hcell_{24}
\\
\hcell_{31} \& \hcell_{32} 
\& \hcell_{33} \& \hcell_{34}
\\
};
\draw[fill=blue,opacity=0.25] (A-3-1.north west) node[blue,left,xshift=0.15cm,yshift=-0.5cm,opacity=1] 
{$\rcell$} to (A-3-4.north east) 
to (A-3-4.south east) to (A-3-1.south west) to (A-3-1.north west);
\draw[fill=red,opacity=0.25] (A-1-3.north west) node[red,above,xshift=0.7cm,opacity=1] 
{$\lcell$} to (A-3-3.south west) 
to (A-3-3.south east) to (A-1-3.north east) to (A-1-3.north west);
\draw[very thick,black,dotted] (A-1-2.north west) to (A-3-2.south west);
\draw[very thick,black,dotted] (A-1-3.north west) to (A-3-3.south west);
\draw[very thick,black,dotted] (A-1-4.north west) to (A-3-4.south west);
\draw[very thick,black,dotted] (A-2-1.north west) to (A-2-4.north east);
\draw[very thick,black,dotted] (A-3-1.north west) to (A-3-4.north east);
\draw[very thick,black] (A-1-1.north west) node[black,above,xshift=-0.5cm] {$\jcell$} to 
(A-1-4.north east) to (A-3-4.south east) 
to (A-3-1.south west) to (A-1-1.north west);
\draw[very thick,black,->] ($(A-1-1.north west)+(-0.4,0.4)$) to (A-1-1.north west);
\draw[very thick,black,->] ($(A-3-4.south east)+(0.5,0.2)$) 
node[right]{$\hcell(\lcell,\rcell)=\hcell_{33}$} 
to[out=180,in=0] ($(A-3-3.south east)+(0,0.2)$);
\end{tikzpicture}
,
\end{gather}
where we use matrix notation for the twelve $H$-cells in $\jcell$. 
In this notation 
left cells are columns, right cells are rows, the $J$-cell is the whole block and 
$H$-cells are the small blocks.

We will also write $<_{l}$ or $\geq_{r}$ 
{\etc}, having the evident meanings.
Note that the three preorders also give rise to 
preorders on the set of cells, as well as between elements of $\monoid$ and cells.
For example, the notations $\lcell\geq_{l}a$ or $\lcell\leq_{l}\lcell^{\prime}$ 
make sense. In particular, for a fixed left cell $\lcell$ we can define
\begin{gather*}
\monoid_{\geq_{l}\lcell}=\{a\in\monoid|a\geq_{l}\lcell\},
\end{gather*}
as well as others which we will distinguish by the subscript.
We write $\hcell(e)$ if $\hcell$ contains an idempotent $e\in \monoid$. 
The $H$-cells of the form $\hcell(e)$ are called \textit{idempotent $H$-cells}, and the $J$-cells $\jcell(e)$ containing these 
$\hcell(e)\subset\jcell(e)$ are called \textit{idempotent $J$-cells}.

\begin{Notation}\label{N:CellsHCells}
When $\hcell(e)\cong\onemon$ is the trivial group, we say $\hcell(e)$ is \textit{trivial}. This is the most relevant case for planar diagram monoids (in fact, one might define planar diagram monoids to be those with trivial $H$-cells). In the following, we will be specifically addressing the case $\hcell(e)\cong\onemon$.
\end{Notation}

We have minimal and maximal $J$-cells in the $\leq_{lr}$-order.
In our illustrations the minimal cell will be at the bottom, so we call it the \textit{bottom cell} $\jb$, while the maximal cell will be at the top, so 
we call it the \textit{top cell} $\jt$.

\begin{Lemma}\label{L:CellsMaximal}
Every monoid has a unique bottom and top $J$-cell 
which are minimal respectively maximal in the $\leq_{lr}$-order. 
Both are idempotent $J$-cells.\null\hfill$\square$
\end{Lemma}

Cells can be considered $\monoid$-representations, called \textit{cell representations} or \textit{Sch{\"u}tzenberger representations}, up to higher order terms:

\begin{Lemma}\label{L:CellsMod}
Each left cell $\lcell$ of $\monoid$ gives rise to a left $\monoid$-representation $\Delta_{\mathcal{L}}=\K\lcell$ by
\begin{gather*}
a\acts l\in\Delta_{\mathcal{L}}=
\begin{cases}
al&\text{if }al\in\lcell,
\\
0&\text{else.}
\end{cases}
\end{gather*}
Similarly, right cells give right $\monoid$-representations $\Delta_{\mathcal{R}}$ and $J$-cells give $\monoid$-birepresentations. We have $\textnormal{dim}_{\mathbb{K}}(\Delta_{\mathcal{L}})=|\lcell|$ and $\textnormal{dim}_{\mathbb{K}}(\Delta_{\mathcal{R}})=|\rcell|$.\null\hfill$\square$
\end{Lemma}

The top and bottom cells of \autoref{L:CellsMaximal} correspond to the two trivial representations of the monoid $\monoid$ as in the following definition.

\begin{Definition}\label{D:RepGapTrivial}
Let $\group\subset\monoid$ be the subgroup of all invertible
elements of $\monoid$, {\ie} $\group$ is the group of units. Then we define \textit{trivial representations}
\begin{gather*}
\oneb\colon\monoid\to\K,\quad s\mapsto
\begin{cases}
1&\text{if }s\in\group,
\\
0&\text{else},
\end{cases}
\quad
\onet\colon\monoid\to\K,\quad s\mapsto 1.
\end{gather*}
An $\monoid$-representation $\module$ is called \textit{trivial} 
if $\module\cong\oneb$ or $\module\cong\onet$.
\end{Definition}

The subscripts \textit{b} and \textit{t} are short for \textit{bottom} and \textit{top}, respectively. The top trivial representation $\onet$ is also what is called the trivial representation $\one$ of $S$, the unit object of the monoidal category of representations of $S$ with $\one\otimes\module\cong\module$ for any $\monoid$-representation $\module$.

\begin{Remark}\label{R:RepGapOrder}
With respect to the above, we again warn the reader familiar with monoid theory that the order we use for $J$-cells ({\aka} Green's $J$-classes) is opposite of the one often used in monoid theory. Thus, what we call bottom/top is usually the top/bottom in monoid theory. In contrast, our convention matches most of the cellular algebra literature.
\end{Remark}

The annihilator $\mathrm{Ann}_{\monoid}(\module)=\{s\in\monoid|s\acts\module=0\}$ 
of an $\monoid$-representation $\module$ is a two-sided ideal of $\monoid$. 
An \textit{apex} of $\module$ is a $J$-cell $\jcell$ 
such that, firstly, $\jcell\cap\mathrm{Ann}_{\monoid}(\module)=\emptyset$, and 
secondly, all $J$-cells $\jcell^{\prime}$ 
with $\jcell^{\prime}\cap\mathrm{Ann}_{\monoid}(\module)=\emptyset$ satisfy 
$\jcell^{\prime}\leq_{lr}\jcell$.  
In other words, an apex is the $\leq_{lr}$-maximal $J$-cell not annihilating $\module$.
The following justifies the terminology of the \textit{apex of a simple $\monoid$-representation}:

\begin{Lemma}\label{L:CellsApex}
Every simple $\monoid$-representation has a unique apex, and apexes correspond 1:1 to idempotent $J$-cells.\null\hfill$\square$
\end{Lemma}

Then, we can classify the simple representations with the \textit{Clifford--Munn--Ponizovski\u{\i} theorem} 
or \textit{$H$-reduction}:

\begin{Theorem}\label{P:CellsSimples}
For a monoid $\monoid$ with only trivial $H$-cells:
\begin{gather*}
\{\text{simple $\monoid$-representations}\}/\cong\;
\xleftrightarrow{1{:}1}
\{\text{idempotent $J$-cells}\}\;
,
\end{gather*}
and the association can be chosen (and we will do this) such that simple $\monoid$-representations are mapped to their apex.\null\hfill$\square$
\end{Theorem}

We can define a partial order, also denoted by $\leq_{lr}$, on the set of simple $\monoid$-representations by saying that one simple is strictly smaller than another if its apex is strictly smaller.

We can estimate the size of simple representations with the following upper bound:

\begin{Lemma}\label{L:CellsDimensions}
The dimension of the simple $\monoid$-representation $L_\jcell$ 
associated to the apex $\jcell$ via \autoref{P:CellsSimples} can be bounded by
\begin{gather*}
\textnormal{dim}_{\mathbb{K}}(L_\jcell)\leq
\text{min}\{|\lcell|,|\rcell|\},
\end{gather*}
where $\lcell,\rcell \subset\jcell$ are arbitrary.
\null\hfill$\square$
\end{Lemma}

We call $\textnormal{ssdim}_{\mathbb{K}}(L_\jcell):=\min\{|\lcell|,|\rcell|\}$ the \textit{semisimple} dimension of $\jcell$. This name is justified by the following:

\begin{Proposition}\label{P:CellsSemisimple}
The following are equivalent for a monoid $\monoid$ with only trivial $H$-cells.
\begin{enumerate}

\item The monoid $\monoid$ is semisimple over $\K$.

\item All $J$-cells are idempotent and square, and we have
$\textnormal{dim}_{\mathbb{K}}(L_\jcell)=\textnormal{ssdim}_{\mathbb{K}}(L_\jcell)$ for all 
$J$-cells $\jcell$, and $L_\jcell\cong\Delta_{\mathcal{L}}$.\null\hfill$\square$

\end{enumerate}
\end{Proposition}

A method for finding more precise dimensions of simple representations is via the construction of Gram matrices, also called sandwich matrices in monoid theory, see {\eg} \cite[Section 5.4]{St-rep-monoid}.
Fix an idempotent $H$-cell $\hcell(e)=\lcell\cap\rcell$ 
in some idempotent $J$-cell $\jcell$. Then 
$\lcell$ is a free right $\hcell(e)$-set and $\rcell$ is a
free left $\hcell(e)$-set, so we can let $\{l_{1},\dots,l_{R}\}$ 
and $\{r_{1},\dots,r_{L}\}$ complete sets of representatives for
$\lcell/\hcell(e)$ respectively for $\hcell(e)\text{\textbackslash}\rcell$.
Here $R$ is the number of right cells and $L$ is the number of left cells in $\jcell$.

The \emph{Gram matrix} $P^{\jcell}=(P^{\jcell}_{i,j})_{i,j}$ is the matrix with values in $\K\hcell(e)$ defined by
\begin{gather*}
P^{\jcell}_{i,j}=
\begin{cases}
r_{i}l_{j}&\text{if }r_{i}l_{j}\in\hcell(e),
\\
0&\text{else}.
\end{cases}
\end{gather*}
Note that $P^{\jcell}$ depends on choices, but one can show that 
its important properties do not depend on these choices, see the references above.

Gram matrices are particularly useful for $\hcell(e)\cong\onemon$ 
as justified by the following (which the reader familiar with \cite{GrLe-cellular} might recognize):

\begin{Proposition}\label{P:CellsGram}
Fix an idempotent $J$-cell $\jcell$. All cells in the statement 
are within $\jcell$.
\\
Assume
$\hcell(e)\subset\jcell$ satisfies $\hcell(e)\cong\onemon$. 
The matrix $P^{\jcell}$ may be rectangular.
Let $L_{\jcell}$ 
denote the associated simple $\monoid$-representation, see \autoref{P:CellsSimples}. Then:
\begin{gather*}
\textnormal{dim}_{\mathbb{K}}(L_{\jcell})=\textnormal{Rank}(P^{\jcell}).
\end{gather*}

\end{Proposition}
\begin{proof}The kernel of the sandwich pairing is the radical of the corresponding cell module: modulo the kernel, multiplication in the principal factor gives rank-one operators with coefficients from $P^{\jcell}$, and these send any nonzero class to scalar multiples of every generating class. The simple quotient therefore has dimension $\operatorname{rank}_{\mathbb K}P^{\jcell}$. This also gives the simple right dimension by using the transpose pairing. See \cite[Section 5.4]{St-rep-monoid} or \cite{Tu-sandwich}.\end{proof}

\subsection{The RepGap}

We now define the main quantities that are important for cryptography. We define them slightly differently when compared to \cite{khovanov-monoidal-2024} (which is the correct definition), but that does not matter for our purposes; see \autoref{L:NotImportant} below.

\begin{Definition}\label{D:RepGap}
The \textit{RepGap} $\textnormal{Gap}_{\mathbb{K}}\,(\monoid)$ is the dimension of the smallest nontrivial simple representation of $\monoid$. 
The \textit{gap ratio} is the quantity $ \textnormal{Ratio}_{\mathbb{K}}\,(\monoid) := \textnormal{Gap}_{\mathbb{K}}\,(\monoid)/\sqrt{|\monoid|}$.
The \textit{semisimple RepGap} $\textnormal{ssGap}_{\mathbb{K}}(\monoid)$ is the minimal \textit{semisimple} dimension of the nontrivial simple $\monoid$-representations. 
The \textit{semisimple gap ratio} is $\textnormal{ssRatio}_{\mathbb{K}}(\monoid):= \textnormal{ssGap}_{\mathbb{K}}(\monoid)/\sqrt{|\monoid|}$
\end{Definition}

The semisimple RepGap is less important than the RepGap, but easier to compute.

\begin{Lemma}\label{L:NotImportant}
The RepGap from \autoref{D:RepGap} is always bigger or equal to the one in \cite{khovanov-monoidal-2024}. Moreover, for our reference diagram monoids, Temperley--Lieb, Motzkin and planar rook, the two definitions agree \cite{khovanov-monoidal-2024,Ar-RepGapMotzkin}.\null\hfill$\square$
\end{Lemma}

Our main results will imply that the rigid diagram monoids have a smaller RepGap than their pivotal counterparts, even when using our more generous RepGap definition compared to \cite{khovanov-monoidal-2024}. Hence, by \autoref{L:NotImportant}, we can still say that they are worse for our purposes, and therefore we do not have to (and will not) worry about the more subtle definition of the RepGap from \cite{khovanov-monoidal-2024}.

\begin{Remark}\label{R:GapRatio}
The gap ratio roughly measures the ratio of the security of the cryptographic protocol based on the monoid to the computational complexity of the monoid. More precisely, recall that the dimension of a simple representation over an algebraically closed field is at most the square root of the size of the monoid. 

If the RepGap is significantly smaller than this upper bound, then the impact from linear attacks is greater; if the monoid is relatively large, the computational complexity for utilizing any encryption protocols based on it is greater, however there is little security gain from this since it can be bypassed by a linear attack, due to the smaller simple representation. In other words, the closer the RepGap is to the square root of the size of the monoid, the less of an impact the linear attack will have. 
\end{Remark}

\begin{Remark}
\textit{Faithfulness} is the minimal dimension of nontrivial faithful representations. We will not compute or discuss it further in this paper, as the RepGap and gap ratio are sufficient for comparing the monoids we are interested in. 
\end{Remark}

In practice, for cryptographic purposes, we generally look for families of monoids indexed by $n \in \mathbb{Z}_{\geq 0}$ such that the \emph{RepGap tends to infinity exponentially with $n$, and the gap ratio does not exponentially tend to zero}. We will discuss examples of such monoids in the next section. 

\section{Diagram categories and monoids}

Martin \cite{Ma} and Jones \cite{Jo} discovered the partition algebra as a generalization of the Temperley–Lieb algebra, extending the construction that relates the Temperley–Lieb algebra to the Potts model in statistical mechanics; see \cite{HaRa-partition-algebras} for a self-contained summary. Specializing the parameter yields the partition monoid and its associated category.

\subsection{Pivotal diagram categories}

The \textit{diagram monoids} we refer to are submonoids of the partition monoid. For $m\in\N$ the picture to keep in mind, taken from \cite{Tu-sandwich}, is:
\begin{enumerate}[label=$\bullet$]

\item The \textit{partition monoid} $Pa_m$ of all diagrams of partitions of a $2m$-element set. The \textit{planar partition monoid} $pPa_m$ is the 
respective planar submonoid
of $Pa_m$.
\begin{gather*}
\begin{tikzpicture}[anchorbase]
\draw[usual] (0.5,0) to[out=90,in=180] (1.25,0.45) to[out=0,in=90] (2,0);
\draw[usual] (0.5,0) to[out=90,in=180] (1,0.35) to[out=0,in=90] (1.5,0);
\draw[usual] (0,1) to[out=270,in=180] (0.75,0.55) to[out=0,in=270] (1.5,1);
\draw[usual] (1.5,1) to[out=270,in=180] (2,0.55) to[out=0,in=270] (2.5,1);
\draw[usual] (0,0) to (0.5,1);
\draw[usual] (1,0) to (1,1);
\draw[usual] (2.5,0) to (2.5,1);
\draw[usual,dot] (2,1) to (2,0.8);
\end{tikzpicture}
\in Pa_m
,\quad
\begin{tikzpicture}[anchorbase]
\draw[usual] (0.5,0) to[out=90,in=180] (1.25,0.45) to[out=0,in=90] (2,0);
\draw[usual] (0.5,0) to[out=90,in=180] (1,0.35) to[out=0,in=90] (1.5,0);
\draw[usual] (0.5,1) to[out=270,in=180] (1,0.55) to[out=0,in=270] (1.5,1);
\draw[usual] (1.5,1) to[out=270,in=180] (2,0.55) to[out=0,in=270] (2.5,1);
\draw[usual] (0,0) to (0,1);
\draw[usual] (2.5,0) to (2.5,1);
\draw[usual,dot] (1,0) to (1,0.2);
\draw[usual,dot] (1,1) to (1,0.8);
\draw[usual,dot] (2,1) to (2,0.8);
\end{tikzpicture}
\in pPa_m
.
\end{gather*}

\item The \textit{rook-Brauer monoid} $RoBr_m$ consisting of all diagrams with components of size $1$ or $2$. The planar rook-Brauer monoid 
$pRoBr_m=Mo_m$ is also called the \textit{Motzkin monoid}.
\begin{gather*}
\begin{tikzpicture}[anchorbase]
\draw[usual] (1,0) to[out=90,in=180] (1.25,0.25) to[out=0,in=90] (1.5,0);
\draw[usual] (1,1) to[out=270,in=180] (1.75,0.55) to[out=0,in=270] (2.5,1);
\draw[usual] (0,0) to (0.5,1);
\draw[usual] (2.5,0) to (2,1);
\draw[usual,dot] (0.5,0) to (0.5,0.2);
\draw[usual,dot] (2,0) to (2,0.2);
\draw[usual,dot] (0,1) to (0,0.8);
\draw[usual,dot] (1.5,1) to (1.5,0.8);
\end{tikzpicture}
\in RoBr_m
,\quad
\begin{tikzpicture}[anchorbase]
\draw[usual] (0.5,0) to[out=90,in=180] (1.25,0.5) to[out=0,in=90] (2,0);
\draw[usual] (1,0) to[out=90,in=180] (1.25,0.25) to[out=0,in=90] (1.5,0);
\draw[usual] (2,1) to[out=270,in=180] (2.25,0.75) to[out=0,in=270] (2.5,1);
\draw[usual] (0,0) to (1,1);
\draw[usual,dot] (2.5,0) to (2.5,0.2);
\draw[usual,dot] (0,1) to (0,0.8);
\draw[usual,dot] (0.5,1) to (0.5,0.8);
\draw[usual,dot] (1.5,1) to (1.5,0.8);
\end{tikzpicture}
\in Mo_m
.
\end{gather*}

\item The \textit{Brauer monoid} $Br_m$ consisting of all diagrams with components of size $2$. The planar Brauer monoid $pBr_m=TL_m$ is known as the \textit{Temperley--Lieb monoid} (sometimes $TL_m$ is called \textit{Jones monoid} or \textit{Kauffman monoid}).
\begin{gather*}
\begin{tikzpicture}[anchorbase]
\draw[usual] (0.5,0) to[out=90,in=180] (1.25,0.45) to[out=0,in=90] (2,0);
\draw[usual] (1,0) to[out=90,in=180] (1.25,0.25) to[out=0,in=90] (1.5,0);
\draw[usual] (0,1) to[out=270,in=180] (0.75,0.55) to[out=0,in=270] (1.5,1);
\draw[usual] (1,1) to[out=270,in=180] (1.75,0.55) to[out=0,in=270] (2.5,1);
\draw[usual] (0,0) to (0.5,1);
\draw[usual] (2.5,0) to (2,1);
\end{tikzpicture}
\in Br_m
,\quad
\begin{tikzpicture}[anchorbase]
\draw[usual] (0.5,0) to[out=90,in=180] (1.25,0.5) to[out=0,in=90] (2,0);
\draw[usual] (1,0) to[out=90,in=180] (1.25,0.25) to[out=0,in=90] (1.5,0);
\draw[usual] (0,1) to[out=270,in=180] (0.25,0.75) to[out=0,in=270] (0.5,1);
\draw[usual] (2,1) to[out=270,in=180] (2.25,0.75) to[out=0,in=270] (2.5,1);
\draw[usual] (0,0) to (1,1);
\draw[usual] (2.5,0) to (1.5,1);
\end{tikzpicture}
\in TL_m
.
\end{gather*}

\item The \textit{rook monoid} or 
\textit{symmetric inverse semigroup} $Ro_m$ consisting of all diagrams with components of size $1$ or $2$, and all blocks have at most one component 
at the bottom and at most one at the top. The \textit{planar rook monoid}
$pRo_m$ is the corresponding planar submonoid.
\begin{gather*}
\begin{tikzpicture}[anchorbase]
\draw[usual] (0,0) to (1,1);
\draw[usual] (0.5,0) to (0,1);
\draw[usual] (2,0) to (2,1);
\draw[usual] (2.5,0) to (0.5,1);
\draw[usual,dot] (1,0) to (1,0.2);
\draw[usual,dot] (1.5,0) to (1.5,0.2);
\draw[usual,dot] (1.5,1) to (1.5,0.8);
\draw[usual,dot] (2.5,1) to (2.5,0.8);
\end{tikzpicture}
\in Ro_m
,\quad
\begin{tikzpicture}[anchorbase]
\draw[usual] (0,0) to (0.5,1);
\draw[usual] (0.5,0) to (1,1);
\draw[usual] (2,0) to (1.5,1);
\draw[usual] (2.5,0) to (2.5,1);
\draw[usual,dot] (1,0) to (1,0.2);
\draw[usual,dot] (1.5,0) to (1.5,0.2);
\draw[usual,dot] (0,1) to (0,0.8);
\draw[usual,dot] (2,1) to (2,0.8);
\end{tikzpicture}
\in pRo_m
.
\end{gather*}

\item The \textit{symmetric monoid/group} $S_m$ consisting of all matchings with components of size $1$. The \textit{planar symmetric group} is trivial $pS_m \cong\onemon$.
\begin{gather*}
\begin{tikzpicture}[anchorbase]
\draw[usual] (0,0) to (1,1);
\draw[usual] (0.5,0) to (0,1);
\draw[usual] (1,0) to (1.5,1);
\draw[usual] (1.5,0) to (2.5,1);
\draw[usual] (2,0) to (2,1);
\draw[usual] (2.5,0) to (0.5,1);
\end{tikzpicture}
\in S_m
,\quad
\begin{tikzpicture}[anchorbase]
\draw[usual] (0,0) to (0,1);
\draw[usual] (0.5,0) to (0.5,1);
\draw[usual] (1,0) to (1,1);
\draw[usual] (1.5,0) to (1.5,1);
\draw[usual] (2,0) to (2,1);
\draw[usual] (2.5,0) to (2.5,1);
\end{tikzpicture}
\in pS_m
.
\end{gather*}

\end{enumerate}
We will only be considering the planar monoids, however we included the non-planar (``symmetric'') monoids above for completeness. We also focus on $TL_m$, $Mo_m$ and $pRo_m$, since $pS_m$ is trivial and $pPa_m \cong TL_{2m}$ by \cite[(1.5)]{HaRa-partition-algebras}.

These monoids each have corresponding monoidal categories, with objects made up of words in a single letter $\bullet$ of length $m \in \mathbb{Z}_{\geq 0}$, and morphism sets consisting of diagrams like the above with each letter in the object acting as a node. The endomorphism sets on objects of length $m$ recover the original monoids. Composition is given by stacking diagrams on top of each other, as long as the objects in the middle are the same, then reducing based on the final pairings of each node (i.e. the partition of nodes). When the diagrams are endomorphisms, this is the same composition as in the original monoids.  

For example, for the Temperley--Lieb category, $TL$, we have:
\begin{gather*}
\begin{tikzpicture}[anchorbase]
\draw[usual] (0.5,0) to[out=90,in=180] (1.25,0.5) to[out=0,in=90] (2,0);
\draw[usual] (1,0) to[out=90,in=180] (1.25,0.25) to[out=0,in=90] (1.5,0);
\draw[usual] (0,0) to (1,1);
\draw[usual] (2.5,0) to (1.5,1);
\draw[usual] (0,0) node [below] {$\bullet$};
\draw[usual] (0.5,0) node [below] {$\bullet$};
\draw[usual] (1,0) node [below] {$\bullet$};
\draw[usual] (1.5,0) node [below] {$\bullet$};
\draw[usual] (2,0) node [below] {$\bullet$};
\draw[usual] (2.5,0) node [below] {$\bullet$};
\draw[usual] (1,1) node [above] {$\bullet$};
\draw[usual] (1.5,1) node [above] {$\bullet$};
\end{tikzpicture}
\in \textnormal{Hom}_{TL}(\bullet\hspace{0.05cm}\bullet\hspace{0.05cm}\bullet\hspace{0.05cm}\bullet\hspace{0.05cm}\bullet\hspace{0.05cm}\bullet,\bullet \hspace{0.05cm}\bullet)
.
\end{gather*}

\begin{Remark}
Further details about these monoids can be found in numerous sources; see, for example, \cite{HaRa-partition-algebras}. Since the partition category can be identified with cobordisms modulo certain relations, \cite{Ko} provides a convenient framework for translating these monoids into categorical language, see \cite{Hu-diagram-categories} for an explicit summary.
\end{Remark}

We also pick particular values $m$ for these monoids, for reasons that will become clear in the following section. From now
\begin{gather*}
\fcolorbox{black!50}{white!50}{we will simply write $\textnormal{TL}_{n},\textnormal{Mo}_{n}$, and $\textnormal{pRo}_{n}$ to mean endomorphisms on $(\bullet \hspace{0.05cm} \bullet)^{\otimes n}$ instead,}
\end{gather*}
which is equivalent to $\textnormal{TL}_{m},\textnormal{Mo}_{m}$, and $\textnormal{pRo}_{m}$ for $m=2n$ using the previous definitions. We also truncate the monoids to keep only the representations that are sufficiently large, following \cite[Section 4]{khovanov-monoidal-2024}. 
Write $D_n^{a\le k\le b}$ for the submonoid of diagrams of rank at most $b$ together with the identity, with the ideal of ranks below $a$ collapsed to zero when it is nonempty. Only admissible integral ranks are retained. For the classical monoids we use
\[
TL_n^T=TL_n^{0\le k\le 2\sqrt{2n}},\qquad
Mo_n^T=Mo_n^{0\le k\le\sqrt{2n}},\qquad
pRo_n^T=pRo_n^{n-\sqrt{2n}\le k\le n+\sqrt{2n}}.
\]
The middle cutoff is the one in \cite[Theorem 4E.1]{Ar-RepGapMotzkin}. All gap ratios use the order of the corresponding truncation.
\begin{Theorem}\label{T:PivotalGaps}
Assume $\textnormal{char}\,\mathbb K=0$ for $Mo_n$ and $TL_n$. The exact equality and all order bounds below hold for sufficiently large $n$.
For the classical Temperley--Lieb gap we have (here $C_n=\frac{1}{n+1}\binom{2n}{n}$)
\[
\textnormal{Gap}_{\mathbb K}(TL_n^T)=C_n-1\sim\frac{4^n}{\sqrt{\pi}\,n^{3/2}}.
\]
The other gaps satisfy
\[
\begin{aligned}
\textnormal{Gap}_{\mathbb K}(Mo_n^T)^{1/n}&\longrightarrow9,&
\textnormal{Gap}_{\mathbb K}(pRo_n^T)&=\Theta(4^n n^{-1/2}).
\end{aligned}
\]
The corresponding gap ratios satisfy
\[
\textnormal{Ratio}_{\mathbb K}(TL_n^T)=\Theta(n^{-3/4}),\quad
\textnormal{Ratio}_{\mathbb K}(Mo_n^T)^{1/n}\longrightarrow1,\quad
\textnormal{Ratio}_{\mathbb K}(pRo_n^T)=\Theta(n^{-1/4}).
\]
\end{Theorem}

\begin{proof}
Write $a_k,d_k$ for the classical standard and simple dimensions at rank $k$ on $2n$ points. Truncation preserves the retained sandwich matrices, whose integer entries let us compute their ranks over $\mathbb C$. The rank-zero simple and the identity-supported character are trivial, and $d_0=1$ since the rank-zero matrix is all ones.

At $q=e^{\pi i/3}$, comparing Weyl and tilting multiplicities in $V_q^{\otimes 2n}$ via \cite[Propositions~3.1 and~6.2]{An-simple-tl} gives
\[
a_{6r+2}=d_{6r+2},\qquad a_{6r}=d_{6r}+d_{6r+4},\qquad
a_{6r+4}=d_{6r+4}+d_{6r+6}.
\]
Here the higher even tilting weights $6r+4$ and $6r+6$ contribute the additional Weyl factors at $6r$ and $6r+4$, respectively. In particular $d_4=a_0-1=C_n-1$.

Put $b_k=a_k/C_n=(k+1)p_{k/2}$. The ballot formula gives $p_0=1$ and
\[
\frac{p_{t+1}}{p_t}=\frac{n-t}{n+t+2}.
\]
For $t\le\sqrt{2n}$ these ratios are $1+O(n^{-1/2})$ and change by $O(n^{-1})$ per step. Thus $p_t\asymp1$ and its first and second differences are $O(n^{-1/2})$ and $O(n^{-1})$. Consequently $b_k\asymp k+1$, and its first and second differences in steps of two are $O(1)$ and $O(n^{-1/2})$.

In particular, $b_{6s+4}-b_{6s}=\tfrac23(b_{6s+6}-b_{6s})+O(n^{-1/2})$. Summing this in the recurrence $d_{6s+6}-d_{6s}=a_{6s+4}-a_{6s}$ gives $d_{6r}/C_n=\tfrac23 b_{6r}+O(1)$. Hence $d_{6r+4}/C_n=\tfrac13 b_{6r+4}+O(1)$, uniformly on the window. Together with the critical ranks $6r+2$, this gives $d_k/C_n\ge c(k+1)-C$ for constants $c>0$ and $C$, so the minimizing rank stays bounded. For fixed $k$, $b_k\to k+1$, and the recurrence gives limits $4r$, $6r+3$, and $2r+1$ at positive ranks $6r$, $6r+2$, and $6r+4$, respectively. Rank $4$ is the unique minimizer of these limits. Stirling's formula gives the asymptotic.

On the retained ranks with $k\asymp\sqrt n$, the classical ballot number $\binom{2n}{n-k/2}-\binom{2n}{n-k/2-1}$ has order $4^n/n$ by Stirling's formula. Squaring and summing over $\Theta(\sqrt n)$ ranks gives a lower bound of order $16^n n^{-3/2}$ for $|TL_n^T|$, which is also the order of the full monoid.

For $Mo_n^T$, \cite[Theorem 4E.1]{Ar-RepGapMotzkin} gives a gap lower bound with $n$th-root limit $9$. The full monoid order has $n$th-root limit $81$ by the Motzkin-number asymptotic \cite{Be-catalan-motzkin,FS-analytic-combinatorics}. Since $\textnormal{Gap}_{\mathbb K}(Mo_n^T)\leq\sqrt{|Mo_n^T|}\leq\sqrt{|Mo_n|}$, this gives both limits.

For $pRo_n^T$, the simple dimension at rank $k$ is $\binom{2n}{k}$ \cite[Proposition 4F.7]{khovanov-monoidal-2024}. Uniformly over the window this has order $4^n/\sqrt n$. Squaring and summing over $\Theta(\sqrt n)$ ranks gives $|pRo_n^T|=\Theta(16^n/\sqrt n)$, and we are done.
\end{proof}

\subsection{Non-pivotal diagram categories}

Define left and right duals in the usual way, e.g. following \cite[Section 2]{egno-tensor-2015}.

\begin{Definition}
A monoidal category is \textit{rigid} if every object has a left and right dual. An object $X$ in a rigid monoidal category with a dual $X^*$ is \textit{pivotal} if $X^{**} \cong X$. A rigid monoidal category is called \textit{pivotal} if every object is pivotal.
\end{Definition}

\begin{Remark}
We have simplified the definition, however one should keep in mind that each dual is actually a monoidal functor, and the isomorphisms are isomorphisms of monoidal functors.  
\end{Remark}

A key feature of the categories described above is that they are all pivotal. This follows from the fact that the object $\bullet$ is self-dual, which is just the zigzag equations: 
\begin{gather*}
\begin{tikzpicture}[anchorbase]
\draw[usual] (0.5,1.25) to[out=90,in=180] (1,2) to[out=0,in=90] (1.5,1.25);
\draw[usual] (1.5,-0.75) to (1.5,0.75);
\draw[usual] (1.5,1) node {$\bullet$};
\draw[usual] (1.5,-0.75) node [below] {$\bullet$};
\draw[usual] (-0.5,0.75) to[out=270,in=180] (0,-0.005) to[out=0,in=270] (0.5,0.75);
\draw[usual] (-0.5,1.25) to (-0.5,2.75);
\draw[usual] (-0.5,1) node {$\bullet$};
\draw[usual] (-0.5,2.75) node [above] {$\bullet$};
\draw[usual] (0.5,1) node {$\bullet$};
\end{tikzpicture}
\;=\;
\begin{tikzpicture}[anchorbase]
\draw[usual] (0,0) to (0,1.5);
\draw[usual] (0,1.5) node [above] {$\bullet$};
\draw[usual] (0,0) node [below] {$\bullet$};
\end{tikzpicture}
\; ; \;
\begin{tikzpicture}[anchorbase]
\draw[usual] (-0.5,1.25) to[out=90,in=180] (0,2) to[out=0,in=90] (0.5,1.25);
\draw[usual] (-0.5,-0.75) to (-0.5,0.75);
\draw[usual] (1.5,1) node {$\bullet$};
\draw[usual] (-0.5,-0.75) node [below] {$\bullet$};
\draw[usual] (0.5,0.75) to[out=270,in=180] (1,0) to[out=0,in=270] (1.5,0.75);
\draw[usual] (1.5,1.25) to (1.5,2.75);
\draw[usual] (-0.5,1) node {$\bullet$};
\draw[usual] (1.5,2.75) node [above] {$\bullet$};
\draw[usual] (0.5,1) node {$\bullet$};
\end{tikzpicture}
\;=\;
\begin{tikzpicture}[anchorbase]
\draw[usual] (0,0) to (0,1.5);
\draw[usual] (0,1.5) node [above] {$\bullet$};
\draw[usual] (0,0) node [below] {$\bullet$};
\end{tikzpicture}
.
\end{gather*}
In the following sections, we will compute the RepGaps and gap ratios for \textit{rigid, non-pivotal} counterparts to the familiar categories above, and compare them to the results for the original pivotal versions. This will give some indication of which, between pivotal and non-pivotal, are more suitable for use in cryptography.

The exploration of non-pivotal analogs is justified with the following result:

\begin{Lemma}
Let $R, P$ be rigid monoidal categories and $F: R \rightarrow P$ a fully faithful monoidal functor. Then $P$ is a pivotal monoidal category $\implies$ $R$ is a pivotal monoidal category. 
\end{Lemma}

\begin{proof}
Since $Y \in P$ is pivotal, $Y$ is the left dual of $Y^*$, so there is an evaluation morphism $\textnormal{ev}_Y : Y \otimes Y^* \rightarrow \mathds{1}_P$. Since $F$ is a fully faithful monoidal functor, $\exists X \in Obj(R)$ and $\textnormal{ev}_X : X \otimes X^* \rightarrow \mathds{1}_R$ such that $F(\textnormal{ev}_X) : F(X) \otimes F(X^*) \rightarrow F(\mathds{1}_R) = \textnormal{ev}_Y : Y \otimes Y^* \rightarrow \mathds{1}_P$. Similarly, there is a corresponding co-evaluation morphism in $R$, which means that $X$ is a left dual of $X^*$, and since duals are unique up to (unique) isomorphism \cite[Proposition 2.10.5]{egno-tensor-2015}, this means $X \cong X^{**}$ so $X$ is pivotal. 

If $P$ is pivotal, the above argument gives a pivotal $X \in Obj(R)$ for every $Y \in Obj(P)$, and since $F$ is faithful then this means that this is all the objects in $R$, i.e., every $X \in Obj(R)$ is pivotal. 
\end{proof}

Let $G$ be a group and $Rep(G)$ be the monoidal category of all finite dimensional representations of $G$. Since $V^{**}\cong V$ for any finite dimensional vector space $V$, $Rep(G)$ is pivotal.

\begin{Corollary}
Let $R$ be a non-pivotal rigid monoidal category. Then there is no fully faithful monoidal functor from $R$ to $Rep(G)$.
\end{Corollary}

This shows that a non-pivotal counterpart to the Temperley--Lieb, Motzkin, and planar rook categories is meaningfully different, since the pivotal versions do have category equivalences to $Rep(G)$ for some group $G$; for example, the Temperley--Lieb category (when the circle is evaluated to $-2$ instead of $1$) is equivalent to $Rep(SL_2)$ via the functor that sends $\bullet$ to the standard two-dimensional $SL_2$-representation. 

\begin{Remark}
The categories we use where the parameter is evaluated to $1$ are actually equivalent to representations of quantum groups (or rather quantum enveloping algebras) and not groups, but the difference will not play a role for us.
\end{Remark}

This result also suggests that the monoids from non-pivotal categories will in some sense be more complicated than those from the pivotal categories. When designing encryption protocols, it would be more desirable to have more complicated, less easily described monoids, since that would hopefully make it harder to decrypt. This is incredibly imprecise; however, the RepGap and gap ratio are useful concrete tools to determine which monoids will likely be more effective for use in cryptography.

In the next three sections, we prove the following bounds on the RepGaps and gap ratios for the non-pivotal $TL$, $Mo$, and $pRo$ monoids:

\begin{Theorem}\label{T:main}
Over every field, the following estimates hold as $n\to\infty$. Write $b_n=\lfloor 2\sqrt{2n}\rfloor$ and $\rho=2^{3/2}3^{3/4}5^{-5/4}<0.87$.
\[
\begin{aligned}
\textnormal{Gap}_{\mathbb K}(rTL_n^T)&=\Theta(2^n n^{-1/2}),&
\textnormal{Ratio}_{\mathbb K}(rTL_n^T)&=\Theta(n^{-1/4}),\\
\textnormal{ssGap}_{\mathbb K}(rMo_n^T)&=\Theta(4^n n^{-3/2}),&
\textnormal{ssRatio}_{\mathbb K}(rMo_n^T)&=\Theta(n^{-1/2}3^{-b_n/2}),\\
\textnormal{Gap}_{\mathbb K}(rpRo_n^T)&=\Theta(2^n n^{-1/2}),&
\textnormal{Ratio}_{\mathbb K}(rpRo_n^T)&=O(\rho^n).
\end{aligned}
\]
\end{Theorem}

\autoref{T:main} will be proven case by case in the next sections. Note that \autoref{T:main} implies that the table in the introduction holds.

\section{The rigid Temperley--Lieb monoid}

We start with our first non-pivotal diagram monoid.

\subsection{Definition of \texorpdfstring{$rTL_n$}{rTLn}}

We now define a non-pivotal analog of the Temperley--Lieb category, a rigid monoidal category that we call the \emph{rigid Temperley--Lieb category} and denote by $rTL$.

\begin{Definition}
The objects, morphisms, and tensor product of $rTL$ are defined as follows:
\begin{enumerate}

\item The collection of \textit{objects} $Obj(rTL)$ consists of all finite words $X_{i_1}X_{i_2}...X_{i_n}$ in the alphabet $\{ X_i :i \in \mathbb{Z} \}$. We will typically only write the subscripts, for example the object $X_1X_2X_1X_2$ will often simply be written $1 \hspace{0.05cm} 2\hspace{0.05cm}1\hspace{0.05cm}2$. The \textit{empty word} is denoted $\mathds{1}$.

\item The \textit{length} of an object $X = X_{i_1}X_{i_2}...X_{i_n}$ is $\textnormal{length}(X)=n$, that is, the number of letters used to make up the word. The length of $\mathds{1}$ is zero. 

\item The collection of \textit{morphisms} between two objects $Hom(X,Y)$, for $X=X_{i_1}X_{i_2}...X_{i_n}$ and $Y=X_{j_1}X_{j_2}...X_{j_m}$, consists of all partitions of $\{1,2,...,m,-1,-2,...,-n \}$ into ordered pairs $(a,b)$, with the following properties:
\begin{enumerate}
\item $(a,-b)$, $a,b>0$, is a valid pair $\iff j_a = i_b$,
\item $(a,b)$, $a>0,b>0$ and $a<b$, is a valid pair $\iff j_b = j_{a}+1$,
\item $(-a,-b)$, $a>0,b>0$ and $a<b$, is a valid pair $\iff i_{b} = i_{a}-1$,
\item the partition is \textit{planar}, meaning there are no pairs $(a,b),(u,v)$ such that $a<u<b<v$, 
\item there are no other valid pairs.
\end{enumerate}

\item It will be useful to represent such morphisms as \textit{diagrams}, drawn vertically from $j_1,...,j_m$ to $i_1,...,i_n$ (i.e. just keeping the subscripts of the two objects, drawing from bottom to top), connecting each pair via an unbroken line, we call a \textit{string}. For example, the morphism $[(1,-1),(2,-6),(-2,-5),(-3,-4)]$ from $X_1X_2X_1X_0X_1X_2 \rightarrow X_1X_2$ would be represented by:
\begin{center}
\begin{tikzpicture}[anchorbase]
\draw[usual] (0.5,0) to[out=90,in=180] (1.25,0.5) to[out=0,in=90] (2,0);
\draw[usual] (1,0) to[out=90,in=180] (1.25,0.25) to[out=0,in=90] (1.5,0);
\draw[usual] (0,0) to (1,1);
\draw[usual] (2.5,0) to (1.5,1);
\draw[usual] (0,0) node [below] {$1$};
\draw[usual] (0.5,0) node [below] {$2$};
\draw[usual] (1,0) node [below] {$1$};
\draw[usual] (1.5,0) node [below] {$0$};
\draw[usual] (2,0) node [below] {$1$};
\draw[usual] (2.5,0) node [below] {$2$};
\draw[usual] (1,1) node [above] {$1$};
\draw[usual] (1.5,1) node [above] {$2$};
\end{tikzpicture}
$\in \textnormal{Hom}_{rTL}(1\hspace{0.05cm}2\hspace{0.05cm}1\hspace{0.05cm}0\hspace{0.05cm}1\hspace{0.05cm}2,1 \hspace{0.05cm}2)$.
\end{center}
The above conditions, when translated to the language of these diagrams, are:
\begin{enumerate}
\item $(a,b)$, with $ab<0$ and $j_{-b} = i_a$, is called a \textit{through strand},
\item $(a,b)$, with $a>0,b>0$, $a<b$ and $i_b = i_{a}+1$, is called a \textit{cup}, denoted cup$_{x,x+1}$, 
\item $(a,b)$, with $a<0,b<0$, $a>b$ and $j_{-b} = j_{-a}-1$, is called a \textit{cap}, denoted cap$_{x,x-1}$,
\item the partition being planar simply means that the diagram is planar, that is, no strings are allowed to intersect,
\item and the diagram is made of only these through strands, cups, and caps. 
\end{enumerate}
Henceforth we will consider diagrams and partitions interchangeably. 

\item The \textit{size} of a diagram $f:X \rightarrow Y$ is $\textnormal{size}(f)=\big(\textnormal{length}(X),\textnormal{length}(Y)\big)$.

\item \textit{Composition} of two morphisms $f:X \rightarrow Y$ and $g:Y \rightarrow Z$, $f\circ g:X\rightarrow Z$, is given by vertically placing the diagram $f$ on top of the diagram $g$, then obtaining a new partition by creating pairs that are connected by strings (noting that a string is not broken by a subscript). For example:
\begin{gather*}
\begin{tikzpicture}[anchorbase]
\draw[usual] (0.5,2) node [above] {$1$};
\draw[usual] (1,2) node [above] {$0$};
\draw[usual] (1.5,2) node [above] {$1$};
\draw[usual] (0,1) node [below] {$1$};
\draw[usual] (0.5,1) node [below] {$2$};
\draw[usual] (1,1) node [below] {$1$};
\draw[usual] (1.5,1) node [below] {$0$};
\draw[usual] (2,1) node [below] {$1$};
\draw[usual] (0.5,1) to[out=90,in=180] (1.25,1.5) to[out=0,in=90] (2,1);
\draw[usual] (1,1) to[out=90,in=180] (1.25,1.25) to[out=0,in=90] (1.5,1);
\draw[usual] (1,2) to[out=270,in=180] (1.25,1.75) to[out=0,in=270] (1.5,2);
\draw[usual] (0,1) to (0.5,2);
\end{tikzpicture}
\circ
\begin{tikzpicture}[anchorbase]
\draw[usual] (0.5,-0.5) to[out=90,in=180] (1.25,0) to[out=0,in=90] (2,-0.5);
\draw[usual] (1,-0.5) to[out=90,in=180] (1.25,-0.25) to[out=0,in=90] (1.5,-0.5);
\draw[usual] (0,0.5) to[out=270,in=180] (0.25,0.25) to[out=0,in=270] (0.5,0.5);
\draw[usual] (1.5,0.5) to[out=270,in=180] (1.75,0.25) to[out=0,in=270] (2,0.5);
\draw[usual] (0,-0.5) to (1,0.5);
\draw[usual] (0,1) node [below] {$1$};
\draw[usual] (0.5,1) node [below] {$2$};
\draw[usual] (1,1) node [below] {$1$};
\draw[usual] (1.5,1) node [below] {$0$};
\draw[usual] (2,1) node [below] {$1$};
\draw[usual] (0,-0.5) node [below] {$1$};
\draw[usual] (0.5,-0.5) node [below] {$2$};
\draw[usual] (1,-0.5) node [below] {$1$};
\draw[usual] (1.5,-0.5) node [below] {$0$};
\draw[usual] (2,-0.5) node [below] {$1$};
\end{tikzpicture}
=
\begin{tikzpicture}[anchorbase]
\draw[usual] (0.5,-0.5) to[out=90,in=180] (1.25,0) to[out=0,in=90] (2,-0.5);
\draw[usual] (1,-0.5) to[out=90,in=180] (1.25,-0.25) to[out=0,in=90] (1.5,-0.5);
\draw[usual] (0,0.5) to[out=270,in=180] (0.25,0.25) to[out=0,in=270] (0.5,0.5);
\draw[usual] (1.5,0.5) to[out=270,in=180] (1.75,0.25) to[out=0,in=270] (2,0.5);
\draw[usual] (0,-0.5) to (1,0.5);
\draw[usual] (0.5,2) node [above] {$1$};
\draw[usual] (1,2) node [above] {$0$};
\draw[usual] (1.5,2) node [above] {$1$};
\draw[usual] (0,1) node [below] {$1$};
\draw[usual] (0.5,1) node [below] {$2$};
\draw[usual] (1,1) node [below] {$1$};
\draw[usual] (1.5,1) node [below] {$0$};
\draw[usual] (2,1) node [below] {$1$};
\draw[usual] (0,-0.5) node [below] {$1$};
\draw[usual] (0.5,-0.5) node [below] {$2$};
\draw[usual] (1,-0.5) node [below] {$1$};
\draw[usual] (1.5,-0.5) node [below] {$0$};
\draw[usual] (2,-0.5) node [below] {$1$};
\draw[usual] (0.5,1) to[out=90,in=180] (1.25,1.5) to[out=0,in=90] (2,1);
\draw[usual] (1,1) to[out=90,in=180] (1.25,1.25) to[out=0,in=90] (1.5,1);
\draw[usual] (1,2) to[out=270,in=180] (1.25,1.75) to[out=0,in=270] (1.5,2);
\draw[usual] (0,1) to (0.5,2);
\end{tikzpicture}
=
\begin{tikzpicture}[anchorbase]
\draw[usual] (0.5,0) to[out=90,in=180] (1.25,0.5) to[out=0,in=90] (2,0);
\draw[usual] (1,0) to[out=90,in=180] (1.25,0.25) to[out=0,in=90] (1.5,0);
\draw[usual] (1,1) to[out=270,in=180] (1.25,0.75) to[out=0,in=270] (1.5,1);
\draw[usual] (0,0) to (0.5,1);
\draw[usual] (0,0) node [below] {$1$};
\draw[usual] (0.5,0) node [below] {$2$};
\draw[usual] (1,0) node [below] {$1$};
\draw[usual] (1.5,0) node [below] {$0$};
\draw[usual] (2,0) node [below] {$1$};
\draw[usual] (0.5,1) node [above] {$1$};
\draw[usual] (1,1) node [above] {$0$};
\draw[usual] (1.5,1) node [above] {$1$};
\end{tikzpicture}
.
\end{gather*}

\item The identity morphism on $X \in Obj(rTL)$, denoted $1_X$, is simply the diagram consisting only of through strands, that is, every subscript in $X$ at the bottom of the diagram is connected to its corresponding subscript at the top via a vertical line. $1_{\mathds{1}}$ is defined as there being no diagram, the empty diagram (or empty partition $[\hspace{0.05cm}]$).

\item Define a functor $\bigotimes : rTL \times rTL \rightarrow rTL$ by:
\begin{enumerate}
\item If $X=X_{i_1}X_{i_2}...X_{i_n}$, $Y=X_{j_1}X_{j_2}...X_{j_m}$ $\in Obj(rTL)$, then $$X \otimes Y = X_{i_1}X_{i_2}...X_{i_n}X_{j_1}X_{j_2}...X_{j_m}$$
that is, we simply concatenate the words $X$ and $Y$, and
\item if $f:X \rightarrow Y \in Hom_{rTL}(X,Y)$, $g: A \rightarrow B \in Hom_{rTL}(A,B)$, then we obtain $f \otimes g \in Hom_{rTL}(X\otimes A,Y\otimes B)$ by placing the diagram $g$ to the right of $f$, that is, concatenating the diagrams horizontally.
\end{enumerate} 
\end{enumerate}
We call the resulting construction the (set-theoretic) rigid Temperley--Lieb category.
\end{Definition}

\begin{Lemma}
With the structure above, $rTL$ is a rigid monoidal category, and its only pivotal object is the tensor unit.
\end{Lemma}

\begin{proof}
The composition of morphisms in $rTL$ is clearly associative, since we can simply stack all multiplied diagrams and obtain the result by following the strings from their top-most objects to the bottom-most objects they are connected to.
Moreover,
$\bigotimes$ defines a monoidal product, with the empty word and diagram ($\mathds{1}$, $1_{\mathds{1}}$) being the monoidal unit. So $rTL$ is a monoidal category.

For every $i \in \mathbb{Z}$, we have:
\begin{enumerate}
\item The \textit{left dual} of $X_i$ is ($X_i^* = X_{i-1}$, cup$_{i-1,i}$, cap$_{i,i-1}$), and
\item the \textit{right dual} of $X_i$ is $(^*X_i=X_{i+1}, \textnormal{cup}_{i,i+1}, \textnormal{cap}_{i+1,i})$, 
\end{enumerate}
which follows from the familiar zigzag equations:
\begin{equation}\label{zigzagrTL}
\begin{tikzpicture}[anchorbase]
\draw[usual] (0.5,1.25) to[out=90,in=180] (1,2) to[out=0,in=90] (1.5,1.25);
\draw[usual] (1.5,-0.75) to (1.5,0.75);
\draw[usual] (1.5,1) node {$i$};
\draw[usual] (1.5,-0.75) node [below] {$i$};
\draw[usual] (-0.5,0.75) to[out=270,in=180] (0,-0.005) to[out=0,in=270] (0.5,0.75);
\draw[usual] (-0.5,1.25) to (-0.5,2.75);
\draw[usual] (-0.5,1) node {$i$};
\draw[usual] (-0.5,2.75) node [above] {$i$};
\draw[usual] (0.5,1) node {$i+1$};
\end{tikzpicture}
\;=\;
\begin{tikzpicture}[anchorbase]
\draw[usual] (0,0) to (0,1.5);
\draw[usual] (0,1.5) node [above] {$i$};
\draw[usual] (0,0) node [below] {$i$};
\end{tikzpicture}
\; ; \;
\begin{tikzpicture}[anchorbase]
\draw[usual] (-0.5,1.25) to[out=90,in=180] (0,2) to[out=0,in=90] (0.5,1.25);
\draw[usual] (-0.5,-0.75) to (-0.5,0.75);
\draw[usual] (1.5,1) node {$i$};
\draw[usual] (-0.5,-0.75) node [below] {$i$};
\draw[usual] (0.5,0.75) to[out=270,in=180] (1,0) to[out=0,in=270] (1.5,0.75);
\draw[usual] (1.5,1.25) to (1.5,2.75);
\draw[usual] (-0.5,1) node {$i$};
\draw[usual] (1.5,2.75) node [above] {$i$};
\draw[usual] (0.5,1) node {$i-1$};
\end{tikzpicture}
\;=\;
\begin{tikzpicture}[anchorbase]
\draw[usual] (0,0) to (0,1.5);
\draw[usual] (0,1.5) node [above] {$i$};
\draw[usual] (0,0) node [below] {$i$};
\end{tikzpicture}
.
\end{equation}
Therefore, given a nonempty word $X=X_{i_1}X_{i_2}...X_{i_n}$, the left dual of $X$ is the object $X^*=X_{i_n-1}\otimes X_{i_{n-1}-1}\otimes \dots \otimes X_{i_2-1} \otimes X_{i_1-1}$, along with the evaluation and coevaluation morphisms, cup$_{X^*,X}$ and cap$_{X,X^*}$, defined by vertically concatenating cups and caps on the corresponding subscripts of $X$ and $X^*$. The right dual is found similarly, and hence $rTL$ is a rigid monoidal category. 

The double left and right duals shift each label by $-2$ and $+2$, respectively. An isomorphism must consist of through strands, since both composites with its inverse must be identity diagrams. Planarity and label preservation then force the source and target words to coincide. Hence no nonempty word is isomorphic to its double dual, and the empty word is its own dual.
\end{proof}

We have defined $rTL$, the \emph{rigid Temperley--Lieb category}, a rigid, non-pivotal analog of the Temperley--Lieb category, and from this we pick the monoid of endomorphisms 
\begin{gather*}
rTL_n := End_{rTL}(1\hspace{0.05cm}2)^{\otimes n}
\end{gather*}
to analyze for cryptographic purposes. We call $rTL_n$ the \emph{rigid Temperley--Lieb monoid}.

\begin{Remark}
The above definition does not include all endomorphisms in $rTL$ of length $2n$. This is because, unlike in $TL$, we have multiple different objects of the same length, meaning the endomorphisms cannot be composed. As a result, the best we can do is pick one object for our endomorphism monoid; the alternative would be to take a direct sum $\bigoplus_{X \in Obj(rTL),\textnormal(X)=m}\textnormal{End}_{rTL}(X)$, however this would require performing a completely different, generally nontrivial, calculation for each different endomorphism monoid, for no discernible extra benefit. 

The choice of $(1\hspace{0.05cm}2)^{\otimes n}$ as our object is because it captures enough of the structure of endomorphisms in $rTL$, while being relatively easier to compute cells. 
\end{Remark}

\subsection{Computing cells for \texorpdfstring{$rTL_n$}{rTLn}}

We have the following sandwich factorization:

\begin{Lemma}\label{factorisationrTL}
Every diagram $d \in rTL_n$ with $k$ through strands can be written as $\tau \circ 1_{\alpha(k)} \circ \beta$, where $\tau$ and $\beta$ are diagrams of size $(2n,k)$ and $(k,2n)$ respectively, and $\alpha(k) = 1\hspace{0.05cm}2\dots 1\hspace{0.05cm}2$ is of length $k$, where $k$ is an even positive integer.
\end{Lemma}

\begin{proof}
Firstly, clearly it is impossible to have a diagram with zero through strands, since there is no valid cap arrangement on $1\hspace{0.05cm}2\dots 1\hspace{0.05cm}2$. Moreover, there cannot be a negative number of through strands, so $k$ must be a positive integer.

Now, consider the possibilities for through strands of a diagram in $rTL_n$. Let $k$ be the number of through strands, and $\alpha(k)$ the sequence, of length $k$, of the letters paired by the through strands, read from left to right. For example, the following diagram has the sequence $1\hspace{0.05cm}2$:
\begin{center}
\begin{tikzpicture}[anchorbase]
\draw[usual] (0.5,0) to[out=90,in=180] (0.75,0.25) to[out=0,in=90] (1,0);
\draw[usual] (1.5,0) to[out=90,in=180] (1.75,0.25) to[out=0,in=90] (2,0);
\draw[usual] (0,1) to[out=270,in=180] (0.25,0.75) to[out=0,in=270] (0.5,1);
\draw[usual] (2,1) to[out=270,in=180] (2.25,0.75) to[out=0,in=270] (2.5,1);
\draw[usual] (0,0) to (1,1);
\draw[usual] (2.5,0) to (1.5,1);
\draw[usual] (0,1) node [above] {$1$};
\draw[usual] (0.5,1) node [above] {$2$};
\draw[usual] (1,1) node [above] {$1$};
\draw[usual] (1.5,1) node [above] {$2$};
\draw[usual] (2,1) node [above] {$1$};
\draw[usual] (2.5,1) node [above] {$2$};
\draw[usual] (0,0) node [below] {$1$};
\draw[usual] (0.5,0) node [below] {$2$};
\draw[usual] (1,0) node [below] {$1$};
\draw[usual] (1.5,0) node [below] {$2$};
\draw[usual] (2,0) node [below] {$1$};
\draw[usual] (2.5,0) node [below] {$2$};
\end{tikzpicture}
.
\end{center}

If the leftmost through strand pairs two 2s, then on the left we will have an odd number of letters on both the top and bottom. Since this is the leftmost through strand, the letters in the top must be paired up independently from the letters in the bottom, and vice versa. This is clearly impossible, since the number of letters in each is odd. A similar argument holds for the rightmost through strand. 

If there is a through strand to the left that pairs the same letters as the through strand to its right, then we will also have an odd number of letters on both the top and bottom, and the same reasoning holds again. This therefore shows that the sequence $\alpha(k)$ must be of the form $1\hspace{0.05cm}2\dots 1\hspace{0.05cm}2$, which is necessarily of even length.

The above reasoning is illustrated below, first with a ``bad" sequence of through strands, which cannot occur for any diagram in $rTL_n$:
\begin{center}
\begin{tikzpicture}[anchorbase]
\draw[usual] (0,0) node {$1$};
\draw[usual] (0.25,0) node {$2$};
\draw[usual] (0.5,0) node {$1$};
\draw[usual] (0.75,0) node {$2$};
\draw[usual] (1,0) node {$1$};
\draw[usual] (1.25,0) node {$2$};
\draw[usual] (1.75,0) node {$\dots$};
\draw[usual] (2.25,0) node {$1$};
\draw[usual] (2.5,0) node {$2$};
\draw[usual] (2.75,0) node {$1$};
\draw[usual] (3,0) node {$2$};
\draw[usual] (3.5,0) node {$\dots$};
\draw[usual] (4,0) node {$1$};
\draw[usual] (4.25,0) node {$2$};
\draw[usual] (4.5,0) node {$1$};
\draw[usual] (4.75,0) node {$2$};
\draw[usual] (5,0) node {$1$};
\draw[usual] (5.25,0) node {$2$};
\draw[usual] (0.75,0.25) to (0.25,1.75);
\draw[usual] (2.5,0.25) to (2.5,1.75);
\draw[usual] (5,0.25) to (4.5,1.75);
\draw[usual] (0,2) node {$1$};
\draw[usual] (0.25,2) node {$2$};
\draw[usual] (0.5,2) node {$1$};
\draw[usual] (0.75,2) node {$2$};
\draw[usual] (1,2) node {$1$};
\draw[usual] (1.25,2) node {$2$};
\draw[usual] (1.75,2) node {$\dots$};
\draw[usual] (2.25,2) node {$1$};
\draw[usual] (2.5,2) node {$2$};
\draw[usual] (2.75,2) node {$1$};
\draw[usual] (3,2) node {$2$};
\draw[usual] (3.5,2) node {$\dots$};
\draw[usual] (4,2) node {$1$};
\draw[usual] (4.25,2) node {$2$};
\draw[usual] (4.5,2) node {$1$};
\draw[usual] (4.75,2) node {$2$};
\draw[usual] (5,2) node {$1$};
\draw[usual] (5.25,2) node {$2$};
\draw[usual] (7,1) node {$\alpha(k)=2\hspace{0.05cm}2\hspace{0.05cm}1$,};
\end{tikzpicture}
\end{center}
then with a ``good" sequence that does appear in some diagrams in $rTL_n$. 
\begin{center}
\begin{tikzpicture}[anchorbase]
\draw[usual] (0,0) node {$1$};
\draw[usual] (0.25,0) node {$2$};
\draw[usual] (0.5,0) node {$1$};
\draw[usual] (0.75,0) node {$2$};
\draw[usual] (1,0) node {$1$};
\draw[usual] (1.25,0) node {$2$};
\draw[usual] (1.75,0) node {$\dots$};
\draw[usual] (2.25,0) node {$1$};
\draw[usual] (2.5,0) node {$2$};
\draw[usual] (2.75,0) node {$1$};
\draw[usual] (3,0) node {$2$};
\draw[usual] (3.5,0) node {$\dots$};
\draw[usual] (4,0) node {$1$};
\draw[usual] (4.25,0) node {$2$};
\draw[usual] (4.5,0) node {$1$};
\draw[usual] (4.75,0) node {$2$};
\draw[usual] (5,0) node {$1$};
\draw[usual] (5.25,0) node {$2$};
\draw[usual] (0,0.25) to (0.5,1.75);
\draw[usual] (2.5,0.25) to (2.5,1.75);
\draw[usual] (4.5,0.25) to (4.5,1.75);
\draw[usual] (5.25,0.25) to (4.75,1.75);
\draw[usual] (0,2) node {$1$};
\draw[usual] (0.25,2) node {$2$};
\draw[usual] (0.5,2) node {$1$};
\draw[usual] (0.75,2) node {$2$};
\draw[usual] (1,2) node {$1$};
\draw[usual] (1.25,2) node {$2$};
\draw[usual] (1.75,2) node {$\dots$};
\draw[usual] (2.25,2) node {$1$};
\draw[usual] (2.5,2) node {$2$};
\draw[usual] (2.75,2) node {$1$};
\draw[usual] (3,2) node {$2$};
\draw[usual] (3.5,2) node {$\dots$};
\draw[usual] (4,2) node {$1$};
\draw[usual] (4.25,2) node {$2$};
\draw[usual] (4.5,2) node {$1$};
\draw[usual] (4.75,2) node {$2$};
\draw[usual] (5,2) node {$1$};
\draw[usual] (5.25,2) node {$2$};
\draw[usual] (7,1) node {$\alpha(k)=1\hspace{0.05cm}2\hspace{0.05cm}1\hspace{0.05cm}2$.};
\end{tikzpicture}
\end{center}
Note that this is only a necessary condition for the sequence of through strands; for example, it is also required that the leftmost through strand connects to the leftmost $1$ on the bottom, and the rightmost through strand connects to the rightmost $2$ on the bottom. 

Finally, the factorisation into $\tau \circ 1_{\alpha(k)} \circ \beta$ follows from observing that the result of a composition depends only on the final pairings between the numbers in the very top object and very bottom object, so we can ``stretch out" the through strands like so:
\begin{center}
\begin{tikzpicture}[anchorbase]
\draw[usual] (0.5,0) to[out=90,in=180] (0.75,0.25) to[out=0,in=90] (1,0);
\draw[usual] (1.5,0) to[out=90,in=180] (1.75,0.25) to[out=0,in=90] (2,0);
\draw[usual] (0,1) to[out=270,in=180] (0.25,0.75) to[out=0,in=270] (0.5,1);
\draw[usual] (2,1) to[out=270,in=180] (2.25,0.75) to[out=0,in=270] (2.5,1);
\draw[usual] (0,0) to (1,1);
\draw[usual] (2.5,0) to (1.5,1);
\draw[usual] (0,1) node [above] {$1$};
\draw[usual] (0.5,1) node [above] {$2$};
\draw[usual] (1,1) node [above] {$1$};
\draw[usual] (1.5,1) node [above] {$2$};
\draw[usual] (2,1) node [above] {$1$};
\draw[usual] (2.5,1) node [above] {$2$};
\draw[usual] (0,0) node [below] {$1$};
\draw[usual] (0.5,0) node [below] {$2$};
\draw[usual] (1,0) node [below] {$1$};
\draw[usual] (1.5,0) node [below] {$2$};
\draw[usual] (2,0) node [below] {$1$};
\draw[usual] (2.5,0) node [below] {$2$};
\end{tikzpicture}
=
\begin{tikzpicture}[anchorbase]
\draw[usual] (0,2) to[out=270,in=180] (0.25,1.75) to[out=0,in=270] (0.5,2);
\draw[usual] (2,2) to[out=270,in=180] (2.25,1.75) to[out=0,in=270] (2.5,2);
\draw[usual] (1,1) to (1,2);
\draw[usual] (1.5,1) to (1.5,2);
\draw[usual] (0,2) node [above] {$1$};
\draw[usual] (0.5,2) node [above] {$2$};
\draw[usual] (1,2) node [above] {$1$};
\draw[usual] (1.5,2) node [above] {$2$};
\draw[usual] (2,2) node [above] {$1$};
\draw[usual] (2.5,2) node [above] {$2$};
\draw[usual] (1,1) node [below] {$1$};
\draw[usual] (1.5,1) node [below] {$2$};
\draw[usual] (1,0.5) to (1,-0.5);
\draw[usual] (1.5,0.5) to (1.5,-0.5);
\draw[usual] (1,-0.5) node [below] {$1$};
\draw[usual] (1.5,-0.5) node [below] {$2$};
\draw[usual] (0.5,-2) to[out=90,in=180] (0.75,-1.75) to[out=0,in=90] (1,-2);
\draw[usual] (1.5,-2) to[out=90,in=180] (1.75,-1.75) to[out=0,in=90] (2,-2);
\draw[usual] (0,-2) to (1,-1);
\draw[usual] (2.5,-2) to (1.5,-1);
\draw[usual] (0,-2) node [below] {$1$};
\draw[usual] (0.5,-2) node [below] {$2$};
\draw[usual] (1,-2) node [below] {$1$};
\draw[usual] (1.5,-2) node [below] {$2$};
\draw[usual] (2,-2) node [below] {$1$};
\draw[usual] (2.5,-2) node [below] {$2$};

\end{tikzpicture}
$=\tau \circ 1_{1\hspace{0.05cm}2} \circ \beta .$
\end{center}
\end{proof}
We will write $1_{\alpha(k)}=1_k$, with $k$ the number of through strands, since the meaning is clear. For $x \in rTL_n$, write $k_x$ for the number of through strands in $x$.

\begin{Lemma}\label{rtlseqorder}
If $x,y \in rTL_n$, then $k_x \geq k_{xy}$ and $k_x \geq k_{yx}$. 
(The number of through strands cannot increase.)
\end{Lemma}

\begin{proof}
In the composition $xy$, each strand arises from a strand in $x$ and a strand in $y$, such that there is a continuous line between the top object of $x$ and the bottom object of $y$. Therefore, it is not possible to increase the number of through strands in $xy$ compared to $x$. It is possible to decrease the number of through strands, since for example we have
\begin{gather*}
\begin{tikzpicture}[anchorbase]
\draw[usual] (0,2) node {$1$};
\draw[usual] (0.5,2) node {$2$};
\draw[usual] (0,0) node {$1$};
\draw[usual] (0.5,0) node {$2$};
\draw[usual] (0,-2) node {$1$};
\draw[usual] (0.5,-2) node {$2$};
\draw[usual] (0,1.75) to (0,0.25);
\draw[usual] (0.5,1.75) to (0.5,0.25);
\draw[usual] (0,-0.25) to[out=270,in=180] (0.25,-0.5) to[out=0,in=270] (0.5,-0.25);
\end{tikzpicture}
=
\begin{tikzpicture}[anchorbase]
\draw[usual] (0,2) node {$1$};
\draw[usual] (0.5,2) node {$2$};
\draw[usual] (0,0) node {$1$};
\draw[usual] (0.5,0) node {$2$};
\draw[usual] (0,1.75) to[out=270,in=180] (0.25,1.5) to[out=0,in=270] (0.5,1.75);
\end{tikzpicture}
.
\end{gather*}
Hence, $1_{k_{xy}} \leq 1_{k_x}$, and similarly $1_{k_{yx}} \leq 1_{k_x}$.
\end{proof}

With the above two lemmas, we can describe the cells of $rTL_n$:

\begin{Proposition}\label{P:tlcells}
The cells of $rTL_n$ have the following properties:
\begin{enumerate}
\item The $J$-cells of $rTL_n$ are indexed by the number of through strands, with order $2n<_{lr}\cdots<_{lr}2$. That is, number of through strands $k$ is the apex of the $J$-cell.

\item The size of the left and right cells of $rTL_n$ is given by the number of diagrams possible when fixing the bottom and top half, respectively.

\item Every $J$-cell in $rTL_n$ is idempotent.
\end{enumerate}
\end{Proposition}

\begin{proof}
The sandwich pictures below, which we will use from now on, are meant to show the factorization in \autoref{factorisationrTL}.
The reader unfamiliar with these illustrations is referred to \cite{Tu-sandwich}.

Fix a $J$-cell $\mathcal{J}$, and suppose $a,b \in J$, meaning $a\leq_{lr} b$ and $b \leq_{lr} a$. This means $\exists c,c',d,d' \in rTL_n$ such that $b = cad$ and $a = c'bd'$.  

Using \autoref{factorisationrTL}, we have:
\begin{center}    
\begin{tikzpicture}[anchorbase,scale=1.5]
\draw[mor] (0,0) to (0.25,0.5) to (0.75,0.5) to (1,0) to (0,0);
\node at (0.5,0.25){$B_b$};
\draw[mor] (0,1.5) to (0.25,1) to (0.75,1) to (1,1.5) to (0,1.5);
\node at (0.5,1.25){$T_b$};
\draw[mor] (0.25,0.5) to (0.25,1) to (0.75,1) to (0.75,0.5) to (0.25,0.5);
\node at (0.5,0.75){\scalebox{1}{$1_{k_b}$}};
\end{tikzpicture}
=
\begin{tikzpicture}[anchorbase,scale=1,decoration={brace,amplitude=7pt}]
\draw[mor] (0,-2) to (0.25,-1.5) to (0.75,-1.5) to (1,-2) to (0,-2);
\node at (0.5,-1.75){$B_d$};
\draw[mor] (0,-0.5) to (0.25,-1) to (0.75,-1) to (1,-0.5) to (0,-0.5);
\node at (0.5,-0.75){$T_d$};
\draw[mor] (0.25,-1.5) to (0.25,-1) to (0.75,-1) to (0.75,-1.5) to (0.25,-1.5);
\node at (0.5,-1.25){\scalebox{0.65}{$1_{k_d}$}};
\draw[mor] (0,-0.5) to (0.25,0) to (0.75,0) to (1,-0.5) to (0,-0.5);
\node at (0.5,-0.25){$B_a$};
\draw[mor] (0,1) to (0.25,0.5) to (0.75,0.5) to (1,1) to (0,1);
\node at (0.5,0.75){$T_a$};
\draw[mor] (0.25,0) to (0.25,0.5) to (0.75,0.5) to (0.75,0) to (0.25,0);
\node at (0.5,0.25){\scalebox{0.65}{$1_{k_a}$}};
\draw[mor] (0,1) to (0.25,1.5) to (0.75,1.5) to (1,1) to (0,1);
\node at (0.5,1.25){$B_c$};
\draw[mor] (0,2.5) to (0.25,2) to (0.75,2) to (1,2.5) to (0,2.5);
\node at (0.5,2.25){$T_c$};
\draw[mor] (0.25,1.5) to (0.25,2) to (0.75,2) to (0.75,1.5) to (0.25,1.5);
\node at (0.5,1.75){\scalebox{0.5}{$1_{k_c}$}};
\draw[decorate] (1,2) to (1,-1.5);
\end{tikzpicture}
$=1_{k_b}$;
\begin{tikzpicture}[anchorbase,scale=1.5]
\draw[mor] (0,0) to (0.25,0.5) to (0.75,0.5) to (1,0) to (0,0);
\node at (0.5,0.25){$B_a$};
\draw[mor] (0,1.5) to (0.25,1) to (0.75,1) to (1,1.5) to (0,1.5);
\node at (0.5,1.25){$T_a$};
\draw[mor] (0.25,0.5) to (0.25,1) to (0.75,1) to (0.75,0.5) to (0.25,0.5);
\node at (0.5,0.75){\scalebox{1}{$1_{k_a}$}};
\end{tikzpicture}
=
\begin{tikzpicture}[anchorbase,scale=1,decoration={brace,amplitude=7pt}]
\draw[mor] (0,-2) to (0.25,-1.5) to (0.75,-1.5) to (1,-2) to (0,-2);
\node at (0.5,-1.75){$B_{d'}$};
\draw[mor] (0,-0.5) to (0.25,-1) to (0.75,-1) to (1,-0.5) to (0,-0.5);
\node at (0.5,-0.75){$T_{d'}$};
\draw[mor] (0.25,-1.5) to (0.25,-1) to (0.75,-1) to (0.75,-1.5) to (0.25,-1.5);
\node at (0.5,-1.25){\scalebox{0.5}{$1_{k_{d'}}$}};
\draw[mor] (0,-0.5) to (0.25,0) to (0.75,0) to (1,-0.5) to (0,-0.5);
\node at (0.5,-0.25){$B_b$};
\draw[mor] (0,1) to (0.25,0.5) to (0.75,0.5) to (1,1) to (0,1);
\node at (0.5,0.75){$T_b$};
\draw[mor] (0.25,0) to (0.25,0.5) to (0.75,0.5) to (0.75,0) to (0.25,0);
\node at (0.5,0.25){\scalebox{0.5}{$1_{k_b}$}};
\draw[mor] (0,1) to (0.25,1.5) to (0.75,1.5) to (1,1) to (0,1);
\node at (0.5,1.25){$B_{c'}$};
\draw[mor] (0,2.5) to (0.25,2) to (0.75,2) to (1,2.5) to (0,2.5);
\node at (0.5,2.25){$T_{c'}$};
\draw[mor] (0.25,1.5) to (0.25,2) to (0.75,2) to (0.75,1.5) to (0.25,1.5);
\node at (0.5,1.75){\scalebox{0.5}{$1_{k_{c'}}$}};
\draw[decorate] (1,2) to (1,-1.5);
\end{tikzpicture}
$=1_{k_a}$.
\end{center}
By \autoref{rtlseqorder}, equivalent diagrams have the same through rank. At a fixed rank each half admits a compatible opposite half. Pair equal cup and cap slots, shifting the final consecutive run of cup slots left by one if it ends at slot $n$. Gluing these opposite halves preserves all through strands, so either half can be replaced by multiplication. For the order, take $e_m$ with cups and caps in the first $m$ slots. If $m\ge l$, stacking $e_m$ above $e_l$ gives $e_m$: the first $l$ cup--cap pairs preserve the two strands emerging from the first $l+1$ bottom sites and the next $m-l$ caps extend that interval to the first $m+1$ sites. Thus $e_m=e_me_l$, giving the order asserted in (a).

The argument for right and left cells is similar. If we have $a,b \in J$, with $k$ through strands, then the condition for $a \sim_r b$ is:
\begin{center}
\begin{tikzpicture}[anchorbase,scale=1.5]
\draw[mor] (0,0) to (0.25,0.5) to (0.75,0.5) to (1,0) to (0,0);
\node at (0.5,0.25){$B_b$};
\draw[mor] (0,1.5) to (0.25,1) to (0.75,1) to (1,1.5) to (0,1.5);
\node at (0.5,1.25){$T_b$};
\draw[mor] (0.25,0.5) to (0.25,1) to (0.75,1) to (0.75,0.5) to (0.25,0.5);
\node at (0.5,0.75){\scalebox{1}{$1_{k}$}};
\end{tikzpicture}
=
\begin{tikzpicture}[anchorbase,scale=1]
\draw[mor] (0,-0.5) to (0.25,0) to (0.75,0) to (1,-0.5) to (0,-0.5);
\node at (0.5,-0.25){$B_c$};
\draw[mor] (0,1) to (0.25,0.5) to (0.75,0.5) to (1,1) to (0,1);
\node at (0.5,0.75){$T_c$};
\draw[mor] (0.25,0) to (0.25,0.5) to (0.75,0.5) to (0.75,0) to (0.25,0);
\node at (0.5,0.25){\scalebox{0.65}{$1_{k_c}$}};
\draw[mor] (0,1) to (0.25,1.5) to (0.75,1.5) to (1,1) to (0,1);
\node at (0.5,1.25){$B_a$};
\draw[mor] (0,2.5) to (0.25,2) to (0.75,2) to (1,2.5) to (0,2.5);
\node at (0.5,2.25){$T_a$};
\draw[mor] (0.25,1.5) to (0.25,2) to (0.75,2) to (0.75,1.5) to (0.25,1.5);
\node at (0.5,1.75){\scalebox{0.65}{$1_{k}$}};
\end{tikzpicture}
;
\begin{tikzpicture}[anchorbase,scale=1.5]
\draw[mor] (0,0) to (0.25,0.5) to (0.75,0.5) to (1,0) to (0,0);
\node at (0.5,0.25){$B_a$};
\draw[mor] (0,1.5) to (0.25,1) to (0.75,1) to (1,1.5) to (0,1.5);
\node at (0.5,1.25){$T_a$};
\draw[mor] (0.25,0.5) to (0.25,1) to (0.75,1) to (0.75,0.5) to (0.25,0.5);
\node at (0.5,0.75){\scalebox{1}{$1_{k}$}};
\end{tikzpicture}
=
\begin{tikzpicture}[anchorbase,scale=1]
\draw[mor] (0,-0.5) to (0.25,0) to (0.75,0) to (1,-0.5) to (0,-0.5);
\node at (0.5,-0.25){$B_{c'}$};
\draw[mor] (0,1) to (0.25,0.5) to (0.75,0.5) to (1,1) to (0,1);
\node at (0.5,0.75){$T_{c'}$};
\draw[mor] (0.25,0) to (0.25,0.5) to (0.75,0.5) to (0.75,0) to (0.25,0);
\node at (0.5,0.25){\scalebox{0.65}{$1_{k_{c'}}$}};
\draw[mor] (0,1) to (0.25,1.5) to (0.75,1.5) to (1,1) to (0,1);
\node at (0.5,1.25){$B_b$};
\draw[mor] (0,2.5) to (0.25,2) to (0.75,2) to (1,2.5) to (0,2.5);
\node at (0.5,2.25){$T_b$};
\draw[mor] (0.25,1.5) to (0.25,2) to (0.75,2) to (0.75,1.5) to (0.25,1.5);
\node at (0.5,1.75){\scalebox{1}{$1_{k}$}};
\end{tikzpicture}
.
\end{center}
which shows that $a \sim _r b \iff T_a=T_b$, i.e. the top halves of the diagrams are the same. Similar reasoning shows that $a\sim_l b \iff B_a = B_b$, proving (b).

Finally, suppose we have a fixed $J$-cell $\mathcal{J}$ of apex $k$ through strands. If $k=2n$ then the identity is an idempotent in $\mathcal{J}$. 

Suppose now that $2 \leq k < 2n$. Then there exists a diagram of the form:
\begin{center}
\begin{tikzpicture}[anchorbase]
\draw[usual] (0.5,0) to[out=90,in=180] (0.75,0.25) to[out=0,in=90] (1,0);
\draw[usual] (1.5,0) to[out=90,in=180] (1.75,0.25) to[out=0,in=90] (2,0);
\draw[usual] (3.5,0) to[out=90,in=180] (3.75,0.25) to[out=0,in=90] (4,0);
\draw[usual] (0,2) to[out=270,in=180] (0.25,1.75) to[out=0,in=270] (0.5,2);
\draw[usual] (1,2) to[out=270,in=180] (1.25,1.75) to[out=0,in=270] (1.5,2);
\draw[usual] (3,2) to[out=270,in=180] (3.25,1.75) to[out=0,in=270] (3.5,2);
\draw[usual] (0,0) to (4,2);
\draw[usual] (4.5,0) to (4.5,2);
\draw[usual] (5,0) to (5,2);
\draw[usual] (5.5,0) to (5.5,2);
\draw[usual] (6.5,0) to (6.5,2);
\draw[usual] (7,0) to (7,2);
\draw[usual] (0,2) node [above] {$1$};
\draw[usual] (0.5,2) node [above] {$2$};
\draw[usual] (1,2) node [above] {$1$};
\draw[usual] (1.5,2) node [above] {$2$};
\draw[usual] (3,2) node [above] {$1$};
\draw[usual] (3.5,2) node [above] {$2$};
\draw[usual] (4,2) node [above] {$1$};
\draw[usual] (4.5,2) node [above] {$2$};
\draw[usual] (5,2) node [above] {$1$};
\draw[usual] (5.5,2) node [above] {$2$};
\draw[usual] (6.5,2) node [above] {$1$};
\draw[usual] (7,2) node [above] {$2$};
\draw[usual] (0,0) node [below] {$1$};
\draw[usual] (0.5,0) node [below] {$2$};
\draw[usual] (1,0) node [below] {$1$};
\draw[usual] (1.5,0) node [below] {$2$};
\draw[usual] (2,0) node [below]{$1$};
\draw[usual] (2.75,0) node {$\dots$};
\draw[usual] (2.25,2) node {$\dots$};
\draw[usual] (3.5,0) node [below] {$2$};
\draw[usual] (4,0) node [below] {$1$};
\draw[usual] (4.5,0) node [below] {$2$};
\draw[usual] (5,0) node [below] {$1$};
\draw[usual] (5.5,0) node [below] {$2$};
\draw[usual] (6,1) node {$\dots$};
\draw[usual] (6.5,0) node [below] {$1$};
\draw[usual] (7,0) node [below] {$2$};
\end{tikzpicture}
\end{center}
where we have simply moved all caps and cups as far to the left as possible. Stacking this diagram on top of itself shows that this diagram is idempotent, and hence every $J$-cell contains an idempotent. 
\end{proof}

\begin{Remark}
$k=2n$ corresponds to the bottom $J$-cell and bottom trivial representation, and $k=2$ corresponds to the top $J$-cell and top trivial representation. Those familiar with the representations of $TL_n$ might be surprised, since the top trivial representation typically corresponds to $k=0$, which is missing for $rTL_n$. This can be reconciled by noting that $k=0$ corresponds to the maximal number of cups and caps in a Temperley--Lieb diagram, while $k=2$ does the same in $rTL_n$. 
\end{Remark}

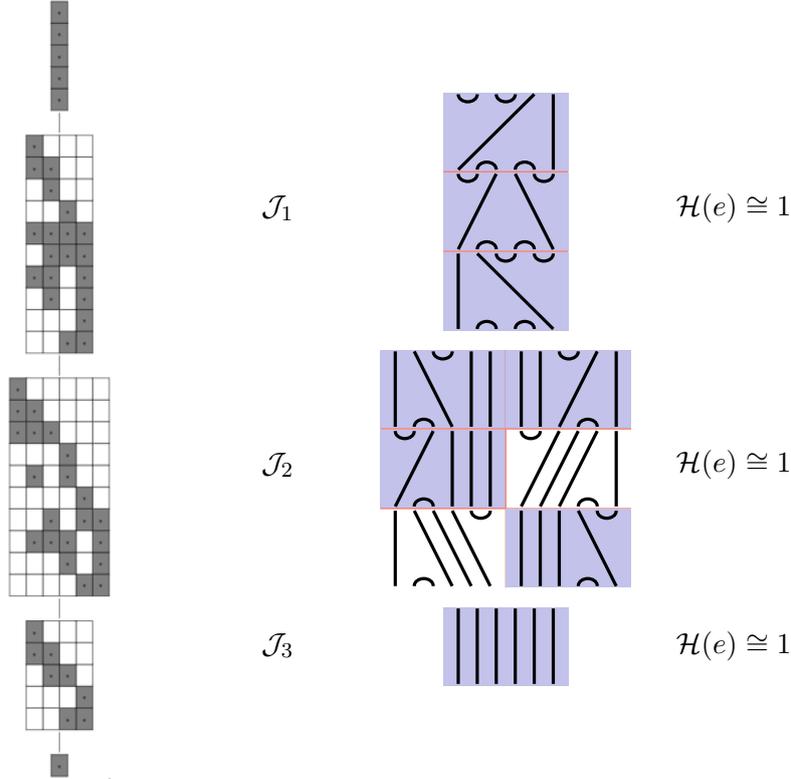
\begin{figure}[ht]
\begin{gather*}
\includegraphics[align=c,scale=0.85]{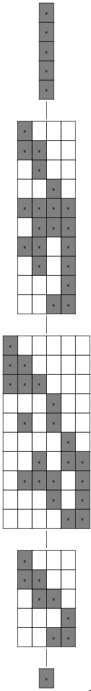}
\xy
(0,0)*{\begin{gathered}
\begin{tabular}{C}
\arrayrulecolor{tomato}
\cellcolor{mydarkblue!25}
\begin{tikzpicture}[anchorbase]
\draw[usual] (0.25,-0.5) to[out=90,in=180] (0.375,-0.4) to[out=0,in=90] (0.5,-0.5);
\draw[usual] (0.75,-0.5) to[out=90,in=180] (0.875,-0.4) to[out=0,in=90] (1,-0.5);
\draw[usual] (0,0.5) to[out=270,in=180] (0.125,0.4) to[out=0,in=270] (0.25,0.5);
\draw[usual] (0.5,0.5) to[out=270,in=180] (0.625,0.4) to[out=0,in=270] (0.75,0.5);
\draw[usual] (0,-0.5) to (1,0.5);
\draw[usual] (1.25,0.5) to (1.25,-0.5);
\end{tikzpicture} 
\\
\hline 
\cellcolor{mydarkblue!25}
\begin{tikzpicture}[anchorbase]
\draw[usual] (0.25,-0.5) to[out=90,in=180] (0.375,-0.4) to[out=0,in=90] (0.5,-0.5);
\draw[usual] (0.75,-0.5) to[out=90,in=180] (0.875,-0.4) to[out=0,in=90] (1,-0.5);
\draw[usual] (0,0.5) to[out=270,in=180] (0.125,0.4) to[out=0,in=270] (0.25,0.5);
\draw[usual] (1,0.5) to[out=270,in=180] (1.125,0.4) to[out=0,in=270] (1.25,0.5);
\draw[usual] (0.5,0.5) to (0,-0.5);
\draw[usual] (0.75,0.5) to (1.25,-0.5);
\end{tikzpicture} 
\\
\hline 
\cellcolor{mydarkblue!25}
\begin{tikzpicture}[anchorbase]
\draw[usual] (0.25,-0.5) to[out=90,in=180] (0.375,-0.4) to[out=0,in=90] (0.5,-0.5);
\draw[usual] (0.75,-0.5) to[out=90,in=180] (0.875,-0.4) to[out=0,in=90] (1,-0.5);
\draw[usual] (0.5,0.5) to[out=270,in=180] (0.625,0.4) to[out=0,in=270] (0.75,0.5);
\draw[usual] (1,0.5) to[out=270,in=180] (1.125,0.4) to[out=0,in=270] (1.25,0.5);
\draw[usual] (0,0.5) to (0,-0.5);
\draw[usual] (0.25,0.5) to (1.25,-0.5);
\end{tikzpicture}
\end{tabular}
\\[3pt]
\begin{tabular}{C|C}
\arrayrulecolor{tomato}
\cellcolor{mydarkblue!25}
\begin{tikzpicture}[anchorbase]
\draw[usual] (0.25,-0.5) to[out=90,in=180] (0.375,-0.4) to[out=0,in=90] (0.5,-0.5);
\draw[usual] (0.5,0.5) to[out=270,in=180] (0.625,0.4) to[out=0,in=270] (0.75,0.5);
\draw[usual] (0,-0.5) to (0,0.5);
\draw[usual] (0.75,-0.5) to (0.25,0.5);
\draw[usual] (1,-0.5) to (1,0.5);
\draw[usual] (1.25,-0.5) to (1.25,0.5);
\end{tikzpicture} &
\cellcolor{mydarkblue!25}
\begin{tikzpicture}[anchorbase]
\draw[usual] (0.75,-0.5) to[out=90,in=180] (0.875,-0.4) to[out=0,in=90] (1,-0.5);
\draw[usual] (0.5,0.5) to[out=270,in=180] (0.625,0.4) to[out=0,in=270] (0.75,0.5);
\draw[usual] (0,-0.5) to (0,0.5);
\draw[usual] (0.25,-0.5) to (0.25,0.5);
\draw[usual] (0.5,-0.5) to (1,0.5);
\draw[usual] (1.25,-0.5) to (1.25,0.5);
\end{tikzpicture} 
\\
\hline
\cellcolor{mydarkblue!25}
\begin{tikzpicture}[anchorbase]
\draw[usual] (0.25,-0.5) to[out=90,in=180] (0.375,-0.4) to[out=0,in=90] (0.5,-0.5);
\draw[usual] (0,0.5) to[out=270,in=180] (0.125,0.4) to[out=0,in=270] (0.25,0.5);
\draw[usual] (0,-0.5) to (0.5,0.5);
\draw[usual] (0.75,-0.5) to (0.75,0.5);
\draw[usual] (1,-0.5) to (1,0.5);
\draw[usual] (1.25,-0.5) to (1.25,0.5);
\end{tikzpicture} &
\begin{tikzpicture}[anchorbase]
\draw[usual] (0.75,-0.5) to[out=90,in=180] (0.875,-0.4) to[out=0,in=90] (1,-0.5);
\draw[usual] (0,0.5) to[out=270,in=180] (0.125,0.4) to[out=0,in=270] (0.25,0.5);
\draw[usual] (0,-0.5) to (0.5,0.5);
\draw[usual] (0.25,-0.5) to (0.75,0.5);
\draw[usual] (0.5,-0.5) to (1,0.5);
\draw[usual] (1.25,-0.5) to (1.25,0.5);
\end{tikzpicture}
\\
\hline
\begin{tikzpicture}[anchorbase]
\draw[usual] (0.25,-0.5) to[out=90,in=180] (0.375,-0.4) to[out=0,in=90] (0.5,-0.5);
\draw[usual] (1,0.5) to[out=270,in=180] (1.125,0.4) to[out=0,in=270] (1.25,0.5);
\draw[usual] (0,-0.5) to (0,0.5);
\draw[usual] (0.75,-0.5) to (0.25,0.5);
\draw[usual] (1,-0.5) to (0.5,0.5);
\draw[usual] (1.25,-0.5) to (0.75,0.5);
\end{tikzpicture} &
\cellcolor{mydarkblue!25}
\begin{tikzpicture}[anchorbase]
\draw[usual] (0.75,-0.5) to[out=90,in=180] (0.875,-0.4) to[out=0,in=90] (1,-0.5);
\draw[usual] (1,0.5) to[out=270,in=180] (1.125,0.4) to[out=0,in=270] (1.25,0.5);
\draw[usual] (0,-0.5) to (0,0.5);
\draw[usual] (0.25,-0.5) to (0.25,0.5);
\draw[usual] (0.5,-0.5) to (0.5,0.5);
\draw[usual] (1.25,-0.5) to (0.75,0.5);
\end{tikzpicture}
\end{tabular}
\\[3pt]
\begin{tabular}{C}
\arrayrulecolor{tomato}
\cellcolor{mydarkblue!25}
\begin{tikzpicture}[anchorbase]
\draw[usual] (0,-0.5) to (0,0.5);
\draw[usual] (0.25,-0.5) to (0.25,0.5);
\draw[usual] (0.5,-0.5) to (0.5,0.5);
\draw[usual] (0.75,-0.5) to (0.75,0.5);
\draw[usual] (1,-0.5) to (1,0.5);
\draw[usual] (1.25,-0.5) to (1.25,0.5);
\end{tikzpicture}
\end{tabular}
\end{gathered}};
(-30,24)*{\jcell_{1}};
(-30,-10)*{\jcell_{2}};
(-30,-34)*{\jcell_{3}};
(30,24)*{\hcell(e)\cong\onemon};
(30,-10)*{\hcell(e)\cong\onemon};
(30,-34)*{\hcell(e)\cong\onemon};
(-51,0)*{\phantom{a}};
\endxy
\end{gather*}
\caption{Left: the cell structure of $rTL_{5}$, produced with the Semigroups package in GAP (adjusted to match our conventions): every block represents a $J$-cell, and every shaded box an $\mathcal{H}$-cell containing an idempotent. Right: the cell structure of $rTL_3$, with each $\mathcal{H}$-cell written explicitly.}
\label{F:tl5}
\end{figure}

\begin{Example}
For $rTL_{3}$ and $rTL_{5}$ the cells are illustrated in \autoref{F:tl5}. It is explained in \cite{St-github-2025} how to produce these pictures in GAP.
\end{Example}

Now, we calculate the sizes of the left and right cells.

\begin{Proposition}\label{rtlcells}
For a fixed number of through strands $k$, we have:
\begin{gather*}
|\mathcal{R}_n^k| = \binom{n-1}{\frac{2n-k}{2}},\quad
|\mathcal{L}_n^k| = \binom{n}{\frac{2n-k}{2}},\quad
\end{gather*}
moreover, the total number of elements in $rTL_n$ is $\binom{2n-1}{n}$.
\end{Proposition}

\begin{proof}
Fix $k$ through strands, and fix the letters in the diagram that each through strand pairs up. Then, all that's left is to arrange the cups and caps on the top and bottom respectively, and the arrangements on the top and bottom are independent of each other. 

Consider the bottom first. There are $2n$ letters on the bottom, and only $\frac{2n-2}{2} = n-1$ spots a cap can be placed (with no nesting allowed), since on the far left and far right we have $1\hspace{0.05cm}2$. Of the $2n$ total letters in the diagram, $k$ are already paired up, so there are only $\frac{2n-k}{2}$ remaining on each of the bottom and top words. 

Once we place the caps, there is only one possible arrangement of the through strands, since the diagrams in $rTL_n$ are planar. Hence, the number of different bottom half diagrams is equal to the number of ways to choose $(2n-k)/2$ pairs of letters (i.e. $1$s and $2$s) to be placed in the $n-1$ valid cap positions. This is equal to $\binom{n-1}{\frac{2n-k}{2}}$, proving the formula for $|\mathcal{R}_n^k|$.

The size of the left cells can be computed similarly, noting that the far left $1$ and the far right $2$ can be used for cups on the top. 

Finally, since the top and bottom are independent once we've fixed $k$ through strands (note that $k$ must be even, by \autoref{factorisationrTL}), the total number of possible diagrams is 
\[
\begin{split}
\sum_{\substack{k=0 \\ k \textnormal{ even}}}^{2n} |\mathcal{R}_n^k| |\mathcal{L}_n^k| = \sum_{k=0}^n \binom{n-1}{\frac{2n-2k}{2}}\binom{n}{\frac{2n-2k}{2}} = \sum_{k=0}^n \binom{n-1}{n-k}\binom{n}{k}
= \binom{2n-1}{n},
\end{split}
\]
where the final step follows from the Chu--Vandermonde identity.
\end{proof}

\begin{Remark}
\autoref{P:tlreps} below shows that $|\mathcal{R}_n^k|$ gives the simple dimension. 
\end{Remark}



\subsection{RepGap and gap ratio of \texorpdfstring{$rTL_n$}{rTLn}}

We are on the final round in our proof of \autoref{T:main} for $rTL_n$. 
As a first step:

\begin{Proposition}\label{P:tlreps}
The simple right $rTL_n$-representations are indexed and ordered by $2n<_{lr}\cdots<_{lr}2$.
Moreover, the right representations $\Delta_{\mathcal{R}_n^k}$ are all simple, and provide a complete list, up to isomorphism, of the simple representations of $rTL_n$. 
\end{Proposition}

\begin{proof}
The first statement is immediate from \autoref{P:CellsSimples} and \autoref{P:tlcells}. Due to the classification \autoref{P:CellsSimples} and \autoref{P:CellsGram}, it is enough to show that the corresponding Gram matrices for each $k$ have maximum rank (the figures use the transpose of the matrix in \autoref{P:CellsGram}).   

Fix the number of through strands $k$. Since $rTL_n$ is planar, an $H$-cell is idempotent if and only if the diagram formed as a product of the bottom half from the corresponding $L$-cell stacked on top of the top half from the corresponding $R$-cell has $k$ through strands. An example is in \autoref{F:rTLgramEx}.

\begin{figure}[ht]
\begin{gather*}
\xy
(0,0)*{\begin{gathered}
\begin{tabular}{C||C|C}
\arrayrulecolor{tomato}
 & 
\begin{tikzpicture}[anchorbase]
\draw[usual] (0.25,-0.5) to[out=90,in=180] (0.375,-0.4) to[out=0,in=90] (0.5,-0.5);
\draw[usual] (0,-0.5) to (0,0.5);
\draw[usual] (0.75,-0.5) to (0.75,0.5);
\draw[usual] (1,-0.5) to (1,0.5);
\draw[usual] (1.25,-0.5) to (1.25,0.5);
\end{tikzpicture} &
\begin{tikzpicture}[anchorbase]
\draw[usual] (0.75,-0.5) to[out=90,in=180] (0.875,-0.4) to[out=0,in=90] (1,-0.5);
\draw[usual] (0,-0.5) to (0,0.5);
\draw[usual] (0.25,-0.5) to (0.25,0.5);
\draw[usual] (0.5,-0.5) to (0.5,0.5);
\draw[usual] (1.25,-0.5) to (1.25,0.5);
\end{tikzpicture}
\\[3pt]
\hline
\hline
\begin{tikzpicture}[anchorbase]
\draw[usual] (0,0.5) to[out=270,in=180] (0.125,0.4) to[out=0,in=270] (0.25,0.5);
\draw[usual] (0.5,-0.5) to (0.5,0.5);
\draw[usual] (0.75,-0.5) to (0.75,0.5);
\draw[usual] (1,-0.5) to (1,0.5);
\draw[usual] (1.25,-0.5) to (1.25,0.5);
\end{tikzpicture} &
\begin{tikzpicture}[anchorbase]
\draw[usual] (0,0) to (0,0.75);
\draw[usual] (0.25,0) to (0.25,0.75);
\draw[usual] (0.5,0) to (0.5,0.75);
\draw[usual] (0.75,0) to (0.75,0.75);
\end{tikzpicture} &
\begin{tikzpicture}[anchorbase]
\draw[usual] (0,0.75) to[out=270,in=180] (0.125,0.65) to[out=0,in=270] (0.25,0.75);
\draw[usual] (0.25,0) to[out=90,in=180] (0.375,0.1) to[out=0,in=90] (0.5,0);
\draw[usual] (0,0) to (0.5,0.75);
\draw[usual] (0.75,0) to (0.75,0.75);
\end{tikzpicture}
\\[3pt]
\hline
\begin{tikzpicture}[anchorbase]
\draw[usual] (0.5,0.5) to[out=270,in=180] (0.625,0.4) to[out=0,in=270] (0.75,0.5);
\draw[usual] (0,-0.5) to (0,0.5);
\draw[usual] (0.25,-0.5) to (0.25,0.5);
\draw[usual] (1,-0.5) to (1,0.5);
\draw[usual] (1.25,-0.5) to (1.25,0.5);
\end{tikzpicture} &
\begin{tikzpicture}[anchorbase]
\draw[usual] (0,0) to (0,0.75);
\draw[usual] (0.25,0) to (0.25,0.75);
\draw[usual] (0.5,0) to (0.5,0.75);
\draw[usual] (0.75,0) to (0.75,0.75);
\end{tikzpicture} &
\begin{tikzpicture}[anchorbase]
\draw[usual] (0,0) to (0,0.75);
\draw[usual] (0.25,0) to (0.25,0.75);
\draw[usual] (0.5,0) to (0.5,0.75);
\draw[usual] (0.75,0) to (0.75,0.75);
\end{tikzpicture}
\\[3pt]
\hline
\begin{tikzpicture}[anchorbase]
\draw[usual] (1,0.5) to[out=270,in=180] (1.125,0.4) to[out=0,in=270] (1.25,0.5);
\draw[usual] (0,-0.5) to (0,0.5);
\draw[usual] (0.25,-0.5) to (0.25,0.5);
\draw[usual] (0.5,-0.5) to (0.5,0.5);
\draw[usual] (0.75,-0.5) to (0.75,0.5);
\end{tikzpicture} &
\begin{tikzpicture}[anchorbase]
\draw[usual] (0,0) to (0,0.75);
\draw[usual] (0.25,0) to[out=90,in=180] (0.375,0.1) to[out=0,in=90] (0.5,0);
\draw[usual] (0.75,0) to (0.25,0.75);
\draw[usual] (0.5,0.75) to[out=270,in=180] (0.625,0.65) to[out=0,in=270] (0.75,0.75);
\end{tikzpicture} &
\begin{tikzpicture}[anchorbase]
\draw[usual] (0,0) to (0,0.75);
\draw[usual] (0.25,0) to (0.25,0.75);
\draw[usual] (0.5,0) to (0.5,0.75);
\draw[usual] (0.75,0) to (0.75,0.75);
\end{tikzpicture}
\end{tabular}
\end{gathered}};
\endxy
\quad,
\end{gather*}
\begin{gather*}
P^{\mathcal{J}_2} =
\begin{pmatrix}
1 & 0 \\
1 & 1 \\
0 & 1
\end{pmatrix}
\quad .
\end{gather*}
\caption{How to produce the $rTL_3$ Gram matrix for $k=4$.}
\label{F:rTLgramEx}
\end{figure}

In $rTL_n$, since a cup joins $1$ to $2$ and a cap joins $2$ to $1$ (reading left to right), cups and caps combine in ``runs" when taking products. A run is called ``even" when the number of caps equals the number of cups, producing a through strand, and ``odd" if there is an imbalance, producing a cap or a cup. An example of each is shown in \autoref{F:TLruns}. 

\begin{figure}[ht]
\begin{gather*}
\begin{tikzpicture}[anchorbase]
\draw[usual] (0,2) node {$1$};
\draw[usual] (0.5,2) node {$2$};
\draw[usual] (1,2) node {$1$};
\draw[usual] (1.5,2) node {$2$};
\draw[usual] (2,2) node {$1$};
\draw[usual] (2.5,2) node {$2$};
\draw[usual] (3,2) node {$1$};
\draw[usual] (0,0) node {$1$};
\draw[usual] (0.5,0) node {$2$};
\draw[usual] (1,0) node {$1$};
\draw[usual] (1.5,0) node {$2$};
\draw[usual] (2,0) node {$1$};
\draw[usual] (2.5,0) node {$2$};
\draw[usual] (3,0) node {$1$};
\draw[usual] (0,-2) node {$1$};
\draw[usual] (0.5,-2) node {$2$};
\draw[usual] (1,-2) node {$1$};
\draw[usual] (1.5,-2) node {$2$};
\draw[usual] (2,-2) node {$1$};
\draw[usual] (2.5,-2) node {$2$};
\draw[usual] (3,-2) node {$1$};
\draw[usual] (0,1.75) to (0,0.25);
\draw[usual] (0,-0.25) to[out=270,in=180] (0.25,-0.5) to[out=0,in=270] (0.5,-0.25);
\draw[usual] (0.5,0.25) to[out=90,in=180] (0.75,0.5) to[out=0,in=90] (1,0.25);
\draw[usual] (1,-0.25) to[out=270,in=180] (1.25,-0.5) to[out=0,in=270] (1.5,-0.25);
\draw[usual] (1.5,0.25) to[out=90,in=180] (1.75,0.5) to[out=0,in=90] (2,0.25);
\draw[usual] (2,-0.25) to[out=270,in=180] (2.25,-0.5) to[out=0,in=270] (2.5,-0.25);
\draw[usual] (2.5,0.25) to[out=90,in=180] (2.75,0.5) to[out=0,in=90] (3,0.25);
\draw[usual] (3,-0.25) to (3,-1.75);
\end{tikzpicture}
=
\begin{tikzpicture}[anchorbase]
\draw[usual] (0,2) node {$1$};
\draw[usual] (0.5,2) node {$2$};
\draw[usual] (1,2) node {$1$};
\draw[usual] (1.5,2) node {$2$};
\draw[usual] (2,2) node {$1$};
\draw[usual] (2.5,2) node {$2$};
\draw[usual] (3,2) node {$1$};
\draw[usual] (0,0) node {$1$};
\draw[usual] (0.5,0) node {$2$};
\draw[usual] (1,0) node {$1$};
\draw[usual] (1.5,0) node {$2$};
\draw[usual] (2,0) node {$1$};
\draw[usual] (2.5,0) node {$2$};
\draw[usual] (3,0) node {$1$};
\draw[usual] (0,1.75) to (3,0.25);
\end{tikzpicture}
\quad,\quad
\begin{tikzpicture}[anchorbase]
\draw[usual] (0,2) node {$1$};
\draw[usual] (0.5,2) node {$2$};
\draw[usual] (1,2) node {$1$};
\draw[usual] (1.5,2) node {$2$};
\draw[usual] (2,2) node {$1$};
\draw[usual] (2.5,2) node {$2$};
\draw[usual] (3,2) node {$1$};
\draw[usual] (0,0) node {$1$};
\draw[usual] (0.5,0) node {$2$};
\draw[usual] (1,0) node {$1$};
\draw[usual] (1.5,0) node {$2$};
\draw[usual] (2,0) node {$1$};
\draw[usual] (2.5,0) node {$2$};
\draw[usual] (3,0) node {$1$};
\draw[usual] (0,-2) node {$1$};
\draw[usual] (0.5,-2) node {$2$};
\draw[usual] (1,-2) node {$1$};
\draw[usual] (1.5,-2) node {$2$};
\draw[usual] (2,-2) node {$1$};
\draw[usual] (2.5,-2) node {$2$};
\draw[usual] (3,-2) node {$1$};
\draw[usual] (0,1.75) to (0,0.25);
\draw[usual] (0,-0.25) to[out=270,in=180] (0.25,-0.5) to[out=0,in=270] (0.5,-0.25);
\draw[usual] (0.5,0.25) to[out=90,in=180] (0.75,0.5) to[out=0,in=90] (1,0.25);
\draw[usual] (1,-0.25) to[out=270,in=180] (1.25,-0.5) to[out=0,in=270] (1.5,-0.25);
\draw[usual] (1.5,0.25) to[out=90,in=180] (1.75,0.5) to[out=0,in=90] (2,0.25);
\draw[usual] (2,-0.25) to[out=270,in=180] (2.25,-0.5) to[out=0,in=270] (2.5,-0.25);
\draw[usual] (2.5,0.25) to (2.5,1.75);
\end{tikzpicture}
=
\begin{tikzpicture}[anchorbase]
\draw[usual] (0,2) node {$1$};
\draw[usual] (0.5,2) node {$2$};
\draw[usual] (1,2) node {$1$};
\draw[usual] (1.5,2) node {$2$};
\draw[usual] (2,2) node {$1$};
\draw[usual] (2.5,2) node {$2$};
\draw[usual] (3,2) node {$1$};
\draw[usual] (0,0) node {$1$};
\draw[usual] (0.5,0) node {$2$};
\draw[usual] (1,0) node {$1$};
\draw[usual] (1.5,0) node {$2$};
\draw[usual] (2,0) node {$1$};
\draw[usual] (2.5,0) node {$2$};
\draw[usual] (3,0) node {$1$};
\draw[usual] (0,1.75) to[out=270,in=180] (1.25,1) to[out=0,in=270] (2.5,1.75);
\end{tikzpicture}
\quad .
\end{gather*}
\caption{Left: an even run; right: an odd run.}
\label{F:TLruns}
\end{figure}

It's clear that the $H$-cells are idempotent precisely when the corresponding product of top and bottom half diagrams has only even runs. 

Label each space that a cup can be placed $1,\dots,n$ and each space a cap can be placed $1,\dots,n-1$. Write $A_r=\{a_1,\dots,a_m\}$ and $B_c=\{b_1,\dots,b_l\}$ for each row $r$ and column $c$, where $a_i \in \{1,\dots,n\}$ are the positions of the cups in the bottom half-diagram in row $r$ written in increasing order, and $b_j \in \{1,\dots,n-1\}$ are defined similarly for the caps. For example, in the table in \autoref{F:rTLgramEx}, the bottom halves in each row have one cup in positions $1$ on nodes $\{1,2\}$, $2$ on nodes $\{3,4\}$, and $3$ on nodes $\{5,6\}$, and the top halves in each column have one cap in positions $1$ on nodes $\{2,3\}$ and $2$ on nodes $\{4,5\}$. Note that a cap in position $b_j$ will connect with a cup in position $a_i\in\{b_j,b_j+1\}$.

Then, the order we use for constructing the Gram matrix is increasing sum of cup/cap index, breaking ties with lexicographical ordering on the nodes for each cup/cap. This ordering is used in the example \autoref{F:rTLgramEx}.

We show that restricting the Gram matrix to the rows $r$ with bottom half cup ordering $A_r \subseteq \{1,\dots,n-1\}$ gives a lower-triangular $|\mathcal{R}_n^k| \times |\mathcal{R}_n^k|$ minor with $1$s along the diagonal, which will complete the proof. 

Firstly, along the diagonal, every cup in space $a_i$ will have a corresponding cap in space $b_i$, forming only even runs and thus giving a $1$ in the Gram matrix.  

Above the diagonal, observe if an $(r,c)$ entry is equal to $1$, then for each $b_j \in B_c$ there must be a corresponding cup in position $b_j$ or $b_j+1$ to maintain an even run. Therefore, $\sum A_r \geq \sum B_c$ with equality when $A_r=B_c$, which makes it impossible for any $(r,c)$ entry to be $1$ when $c>r$, i.e. above the diagonal. This is because in the ordering defined above we would require $\sum B_c \geq \sum A_r$ which forces $B_c=A_r$ when the entry is $1$, i.e. $r=c$ since they are identical in lexicographical ordering.

Finally, this minor is clearly a square matrix of size $|\mathcal{R}_n^k| = \dbinom{n-1}{\#\text{caps}} = \dbinom{n-1}{\frac{2n-k}{2}}$.
\end{proof}

\begin{Remark}\label{P:tlrepsss}\label{R:semisimple}
All right cell representations are simple, but
the monoid $rTL_n$ is not semisimple: being semisimple would imply that the sum of the squares of the simple representations is the order of the monoid, but that is not true as our $J$-cells are not squares; {\cf} \autoref{rtlcells}. However, we can think of $rTL_n$ as being almost semisimple: only the slight asymmetry in its left and right cells sizes prevents it from being semisimple.
\end{Remark}

\begin{figure*}[ht]
\begin{gather*}
\includegraphics[scale=0.65]{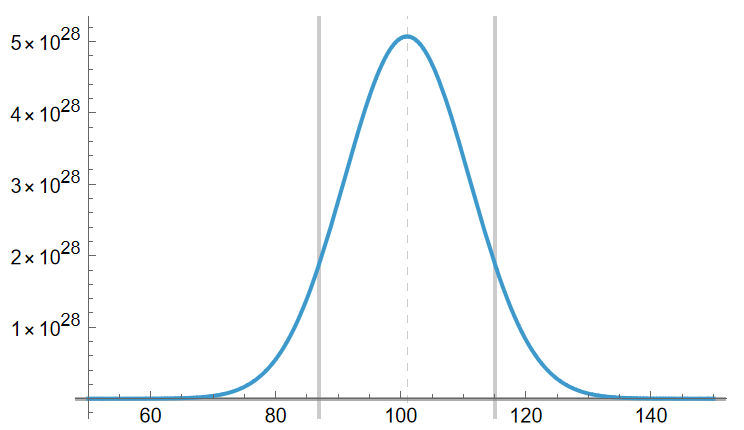}
\includegraphics[scale=0.65]{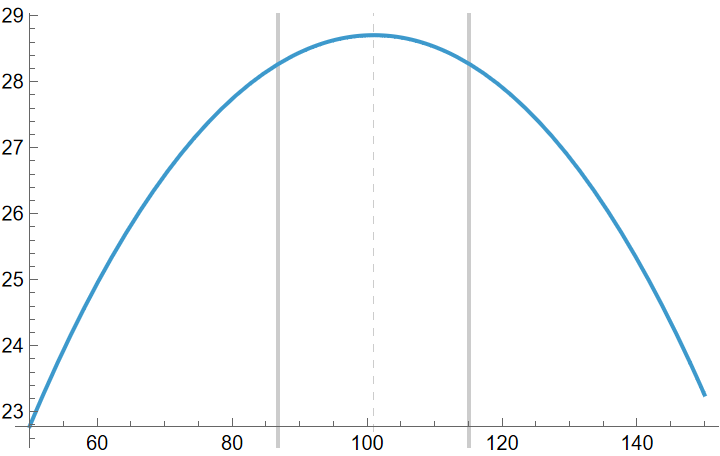}
\end{gather*}
\caption{Left: The simple dimensions over $k$, for $n=100$, with vertical lines marking the truncation end points. Right: A Log10 version of the same graph.}
\label{rTL_truncation}
\end{figure*}

Following \cite{khovanov-monoidal-2024},
we will truncate $rTL_n$, similarly to $TL_n$. 
From \autoref{rtlcells} we find
that $|\mathcal{R}_n^k|$ is smaller than $|\mathcal{L}_n^k|$ and its maximum occurs at $k=n+1$ (simply by the properties of the binomial coefficient).  
We consider the interval $n+1-\sqrt{2n} \leq k \leq n+1+ \sqrt{2n}$, which has a width on the order of $\sqrt{2n}$, up to a scalar, and is centered around the peak at $k=n+1$.
This is illustrated in \autoref{rTL_truncation}. To avoid unwieldy superscripts, as with the pivotal monoids, we will denote this truncated monoid by $rTL_n^T$.

Now we can find the RepGap and gap ratio, and compare the results with those for $TL_n$. 

\begin{Theorem} 
Over every field $\mathbb K$,
\[
\textnormal{Gap}_{\mathbb K}(rTL_n^T)=\Theta(2^n n^{-1/2}),\qquad
\textnormal{Ratio}_{\mathbb K}(rTL_n^T)=\Theta(n^{-1/4}).
\]
\end{Theorem}

\begin{proof}
The right cell representations are simple by \autoref{P:tlreps} and have dimension $\binom{n-1}{m}$ by \autoref{rtlcells}, where $m=(2n-k)/2$. Uniformly over the retained window, Stirling's formula gives order $2^n/\sqrt n$. The order of the truncation is $2+\sum_k\binom n m\binom{n-1}m$, summed over the retained even ranks, with the extra two terms accounting for the identity and zero for sufficiently large $n$. Each product has order $4^n/n$ and there are $\Theta(\sqrt n)$ ranks, giving the claimed ratio.
\end{proof}

Combining this with \autoref{T:PivotalGaps}, we have the following conclusion for $\textnormal{char}\,\mathbb{K}=0$ (which is only needed for $TL_n$ itself, see \autoref{R:Tl}).

\begin{Corollary}
Assume $\textnormal{char}\,\mathbb{K}=0$. Then
\[
\frac{\textnormal{Gap}_{\mathbb K}(rTL_n^T)}{\textnormal{Gap}_{\mathbb K}(TL_n^T)}=\Theta(n 2^{-n}),\qquad
\frac{\textnormal{Ratio}_{\mathbb K}(rTL_n^T)}{\textnormal{Ratio}_{\mathbb K}(TL_n^T)}=\Theta(n^{1/2}).
\]
Thus, we get the table entries as in the introduction.
\end{Corollary}

Although the gap ratio for $rTL_n^T$ is larger than the gap ratio for $TL_n^T$, the fact that it is significantly worse in terms of RepGap makes it less suitable for cryptographic purposes. 

\begin{Remark}\label{R:Tl}
The behavior of the monoid $TL_n$/category $TL$ is very different in prime characteristic compared to $\textnormal{char}\,\mathbb{K}=0$, see e.g. \cite{An-simple-tl,Sp-modular-tl,SuTuWeZh-mixed-tilting} (where we specialize so that the circle parameter is one in these papers). The simple dimensions of $rTL_n$ are field-independent by \autoref{P:tlreps}. 
\end{Remark}

\section{The rigid Motzkin monoid}

In this section we define a non-pivotal version of the Motzkin algebra/monoid as, for example, in \cite{BH}. This will be done in a similar way to the previous section, so we will be brief whenever appropriate.

\subsection{Definition of \texorpdfstring{$rMo_n$}{rMon}}

Similarly to $rTL_n$ we define:

\begin{Definition}
The \emph{rigid Motzkin category} $rMo$ has objects the same as those in $rTL$, and morphisms defined as partitions of those objects, similarly to $rTL$, however we also allow parts of size $1$, i.e. we allow unpaired partitions. In diagram form, this would be denoted by a \textit{top dot} or \textit{bottom dot}, generally referred to as \textit{dots}, as seen in the following example:

\begin{center}
\begin{tikzpicture}[anchorbase]
\draw[usual] (0.5,0) to[out=90,in=180] (1.25,0.5) to[out=0,in=90] (2,0);
\draw[usual] (1,0) to[out=90,in=180] (1.25,0.25) to[out=0,in=90] (1.5,0);
\draw[usual] (2,1) to[out=270,in=180] (2.25,0.75) to[out=0,in=270] (2.5,1);
\draw[usual] (0,0) to (1,1);
\draw[usual,dot] (2.5,0) to (2.5,0.2);
\draw[usual,dot] (0,1) to (0,0.8);
\draw[usual,dot] (0.5,1) to (0.5,0.8);
\draw[usual,dot] (1.5,1) to (1.5,0.8);
\draw[usual] (0,1) node [above] {$1$};
\draw[usual] (0.5,1) node [above] {$2$};
\draw[usual] (1,1) node [above] {$1$};
\draw[usual] (1.5,1) node [above] {$0$};
\draw[usual] (2,1) node [above] {$1$};
\draw[usual] (2.5,1) node [above] {$2$};
\draw[usual] (0,0) node [below] {$1$};
\draw[usual] (0.5,0) node [below] {$2$};
\draw[usual] (1,0) node [below] {$1$};
\draw[usual] (1.5,0) node [below] {$0$};
\draw[usual] (2,0) node [below] {$1$};
\draw[usual] (2.5,0) node [below] {$2$};
\end{tikzpicture}
.
\end{center}
Composition of diagrams is defined similarly to $rTL$, for example: 
\begin{gather*}
\begin{tikzpicture}[anchorbase]
\draw[usual] (0.5,-0.5) to[out=90,in=180] (1.25,0) to[out=0,in=90] (2,-0.5);
\draw[usual] (1,-0.5) to[out=90,in=180] (1.25,-0.25) to[out=0,in=90] (1.5,-0.5);
\draw[usual] (2,0.5) to[out=270,in=180] (2.25,0.25) to[out=0,in=270] (2.5,0.5);
\draw[usual] (0,-0.5) to (1,0.5);
\draw[usual,dot] (2.5,-0.5) to (2.5,-0.3);
\draw[usual,dot] (0,0.5) to (0,0.3);
\draw[usual,dot] (0.5,0.5) to (0.5,0.3);
\draw[usual,dot] (1.5,0.5) to (1.5,0.3);
\draw[usual] (0,2) node [above] {$1$};
\draw[usual] (0.5,2) node [above] {$2$};
\draw[usual] (1,2) node [above] {$1$};
\draw[usual] (1.5,2) node [above] {$0$};
\draw[usual] (2,2) node [above] {$1$};
\draw[usual] (2.5,2) node [above] {$2$};
\draw[usual] (0,1) node [below] {$1$};
\draw[usual] (0.5,1) node [below] {$2$};
\draw[usual] (1,1) node [below] {$1$};
\draw[usual] (1.5,1) node [below] {$0$};
\draw[usual] (2,1) node [below] {$1$};
\draw[usual] (2.5,1) node [below] {$2$};
\draw[usual] (0,-0.5) node [below] {$1$};
\draw[usual] (0.5,-0.5) node [below] {$2$};
\draw[usual] (1,-0.5) node [below] {$1$};
\draw[usual] (1.5,-0.5) node [below] {$0$};
\draw[usual] (2,-0.5) node [below] {$1$};
\draw[usual] (2.5,-0.5) node [below] {$2$};
\draw[usual] (0.5,1) to[out=90,in=180] (1.25,1.5) to[out=0,in=90] (2,1);
\draw[usual] (1,1) to[out=90,in=180] (1.25,1.25) to[out=0,in=90] (1.5,1);
\draw[usual] (2,2) to[out=270,in=180] (2.25,1.75) to[out=0,in=270] (2.5,2);
\draw[usual] (0,1) to (1,2);
\draw[usual,dot] (2.5,1) to (2.5,1.2);
\draw[usual,dot] (0,2) to (0,1.8);
\draw[usual,dot] (0.5,2) to (0.5,1.8);
\draw[usual,dot] (1.5,2) to (1.5,1.8);
\end{tikzpicture}
=
\begin{tikzpicture}[anchorbase]
\draw[usual] (0.5,0) to[out=90,in=180] (1.25,0.5) to[out=0,in=90] (2,0);
\draw[usual] (1,0) to[out=90,in=180] (1.25,0.25) to[out=0,in=90] (1.5,0);
\draw[usual] (2,1) to[out=270,in=180] (2.25,0.75) to[out=0,in=270] (2.5,1);
\draw[usual,dot] (0,0) to (0,0.2);
\draw[usual,dot] (1,1) to (1,0.8);
\draw[usual,dot] (2.5,0) to (2.5,0.2);
\draw[usual,dot] (0,1) to (0,0.8);
\draw[usual,dot] (0.5,1) to (0.5,0.8);
\draw[usual,dot] (1.5,1) to (1.5,0.8);
\draw[usual] (0,0) node [below] {$1$};
\draw[usual] (0.5,0) node [below] {$2$};
\draw[usual] (1,0) node [below] {$1$};
\draw[usual] (1.5,0) node [below] {$0$};
\draw[usual] (2,0) node [below] {$1$};
\draw[usual] (2.5,0) node [below] {$2$};
\draw[usual] (0,1) node [above] {$1$};
\draw[usual] (0.5,1) node [above] {$2$};
\draw[usual] (1,1) node [above] {$1$};
\draw[usual] (1.5,1) node [above] {$0$};
\draw[usual] (2,1) node [above] {$1$};
\draw[usual] (2.5,1) node [above] {$2$};
\end{tikzpicture}
.
\end{gather*}
Everything else, such as the tensor product, tensor unit, duals, etc., are defined in the same way as in $rTL$. Finally, the \emph{rigid Motzkin monoid} is
\begin{gather*}
rMo_n := \textnormal{End}_{rMo}(1\hspace{0.05cm}2)^{\otimes n},
\end{gather*}
where, as before, our object of choice is $(1\hspace{0.05cm}2)^{\otimes n}$.
\end{Definition}

We leave it to the reader to formulate and verify that $rMo$ shares the same categorical properties as $rTL$.

\subsection{Computing cells for \texorpdfstring{$rMo_n$}{rMon}}

We proceed analogously to $rTL$:

\begin{Lemma}\label{factorisationRMO}
Any $d \in rMo_n$ can be written 
\begin{equation}
d = \tau \circ 1_{\alpha(k)} \circ \beta
\end{equation}
for diagrams $\tau, \beta \in rMo$ and some sequence $\alpha(k)$ of length $k$ in $1$ and $2$, $0 \leq k \leq 2n$.
\end{Lemma}
\begin{proof}
If $d \in rMo_n$ has $k$ through strands, define the sequence $\alpha(k)$ as the sequence of through strands read from left to right; explicitly, $\alpha(k)_i = 1$ if the $i$th through strand pairs two $1$s, and $\alpha(k)_i = 2$ if the $i$th through strand pairs two $2$s. For example, the following diagram has the sequence $1\hspace{0.05cm}2$:
\begin{center}
\begin{tikzpicture}[anchorbase]

\draw[usual] (1.5,0) to[out=90,in=180] (1.75,0.25) to[out=0,in=90] (2,0);
\draw[usual] (0,1) to[out=270,in=180] (0.25,0.75) to[out=0,in=270] (0.5,1);
\draw[usual] (2,1) to[out=270,in=180] (2.25,0.75) to[out=0,in=270] (2.5,1);
\draw[usual] (0,0) to (1,1);
\draw[usual] (0.5,0) to (1.5,1);
\draw[usual,dot] (2.5,0) to (2.5,0.2);
\draw[usual,dot] (1,0) to (1,0.2);
\draw[usual] (0,1) node [above] {$1$};
\draw[usual] (0.5,1) node [above] {$2$};
\draw[usual] (1,1) node [above] {$1$};
\draw[usual] (1.5,1) node [above] {$2$};
\draw[usual] (2,1) node [above] {$1$};
\draw[usual] (2.5,1) node [above] {$2$};
\draw[usual] (0,0) node [below] {$1$};
\draw[usual] (0.5,0) node [below] {$2$};
\draw[usual] (1,0) node [below] {$1$};
\draw[usual] (1.5,0) node [below] {$2$};
\draw[usual] (2,0) node [below] {$1$};
\draw[usual] (2.5,0) node [below] {$2$};
\end{tikzpicture}
.
\end{center}
Note that, unlike $rTL_n$, the number of through strands can be any integer $0 \leq k \leq 2n$. 

Given the sequence $\alpha(k)$, the result follows via the same argument as that used for \autoref{factorisationrTL}, as once again only the final partition matters when determining the result of a composition. For example, the above diagram with $\alpha(k)=1\hspace{0.05cm}2$ becomes: 
\begin{center}
\begin{tikzpicture}[anchorbase]
\draw[usual] (0,2) to[out=270,in=180] (0.25,1.75) to[out=0,in=270] (0.5,2);
\draw[usual] (2,2) to[out=270,in=180] (2.25,1.75) to[out=0,in=270] (2.5,2);
\draw[usual] (1,1) to (1,2);
\draw[usual] (1.5,1) to (1.5,2);
\draw[usual] (0,2) node [above] {$1$};
\draw[usual] (0.5,2) node [above] {$2$};
\draw[usual] (1,2) node [above] {$1$};
\draw[usual] (1.5,2) node [above] {$2$};
\draw[usual] (2,2) node [above] {$1$};
\draw[usual] (2.5,2) node [above] {$2$};
\draw[usual] (1,1) node [below] {$1$};
\draw[usual] (1.5,1) node [below] {$2$};
\draw[usual] (1,0.5) to (1,-0.5);
\draw[usual] (1.5,0.5) to (1.5,-0.5);
\draw[usual] (1,-0.5) node [below] {$1$};
\draw[usual] (1.5,-0.5) node [below] {$2$};
\draw[usual] (1.5,-2) to[out=90,in=180] (1.75,-1.75) to[out=0,in=90] (2,-2);
\draw[usual,dot] (2.5,-2) to (2.5,-1.8);
\draw[usual,dot] (1,-2) to (1,-1.8);
\draw[usual] (0,-2) to (1,-1);
\draw[usual] (0.5,-2) to (1.5,-1);
\draw[usual] (0,-2) node [below] {$1$};
\draw[usual] (0.5,-2) node [below] {$2$};
\draw[usual] (1,-2) node [below] {$1$};
\draw[usual] (1.5,-2) node [below] {$2$};
\draw[usual] (2,-2) node [below] {$1$};
\draw[usual] (2.5,-2) node [below] {$2$};

\end{tikzpicture}
$=\tau \circ 1_{1\hspace{0.05cm}2} \circ \beta$.
\end{center}
The proof is complete.
\end{proof}

Now, for $n \geq 0$ fix a sequence of $k$ through strands $\alpha(k)$. Since the diagrams in $rMo_n$ are planar, the through strands create an ordered partition of both the top and bottom objects. For the above example, $\alpha(k)$ partitions the bottom and top object into the blocks $[[],[],[1\hspace{0.05cm}2\hspace{0.05cm}1\hspace{0.05cm}2]]$ and $[[1\hspace{0.05cm}2],[],[1\hspace{0.05cm}2]]$ respectively, including the empty blocks. 

Denote by $j_i$, $1\leq i \leq k+1$, the size of the block to the left of the through strand $\alpha(k)_i$. Then we have an ordered partition 
\begin{equation}\label{rmopartition}
j_1+j_2+\dots j_{k+1} = 2n - k.
\end{equation}
Define a partial order on the sequences $\alpha(k)$ by $\alpha(k) \leq \beta(k') \iff \alpha(k)$ is a subsequence of $\beta(k')$. Let $k_x$ denote the number of through strands in the diagram $x$, and $\alpha(k_x)$ the corresponding sequence. 
\begin{Lemma}\label{rmoseqorder}
For any diagrams $x,y \in rMo_n$, $\alpha(k_x) \geq \alpha(k_{xy})$ and $\alpha(k_x) \geq \alpha(k_{yx})$.
\end{Lemma}
\begin{proof}
This follows by observing that composing two diagrams can never add any through strands, while they can remove them. Through strands only occur in the composition $xy$ when a number in the top object of $x$ is connected to a number in the bottom object of $y$ via a continuous line. Therefore, each through strand in $xy$ is a result of a composition between a through strand in $x$ and a through strand in $y$, possibly with some caps and cups in the middle. 

This naturally means that there is no way to add any new through strands via composition, since they each must correspond to a pair of through strands from $x,y$. In all other cases, the through strand from $x$ or $y$ will disappear in $xy$, to become either a cap, cup, or dot, since there is no continuous line between the top object of $x$ and the bottom object of $y$.

As a result, the resulting sequence of through strands of $xy$ must be a subsequence of both $\alpha(k_x)$ and $\alpha(k_y)$.
The diagrams below summarize some of the possible cases for compositions involving through strands:
\begin{gather*}
\begin{tikzpicture}[anchorbase]
\draw[usual] (0,2) node {$1$};
\draw[usual] (0.5,2) node {$1$};
\draw[usual] (1,2) node {$1$};
\draw[usual] (1.5,2) node {$2$};
\draw[usual] (2,2) node {$1$};
\draw[usual] (2.5,2) node {$1$};
\draw[usual] (3,2) node {$2$};
\draw[usual] (0,0) node {$1$};
\draw[usual] (0.5,0) node {$1$};
\draw[usual] (1,0) node {$1$};
\draw[usual] (1.5,0) node {$2$};
\draw[usual] (2,0) node {$1$};
\draw[usual] (2.5,0) node {$1$};
\draw[usual] (3,0) node {$2$};
\draw[usual] (0,-2) node {$1$};
\draw[usual] (0.5,-2) node {$1$};
\draw[usual] (1,-2) node {$1$};
\draw[usual] (1.5,-2) node {$2$};
\draw[usual] (2,-2) node {$1$};
\draw[usual] (2.5,-2) node {$1$};
\draw[usual] (3,-2) node {$2$};
\draw[usual] (0,-1.75) to (0,-0.25);
\draw[usual] (0,0.25) to (0,1.75);
\draw[usual] (0.5,0.25) to (0.5,1.75);
\draw[usual] (1,0.25) to (1,1.75);
\draw[usual,dot] (0.5,-0.25) to (0.5,-0.45);
\draw[usual] (1.5,0.25) to[out=90,in=180] (1.75,0.5) to[out=0,in=90] (2,0.25);
\draw[usual] (1,-0.25) to[out=270,in=180] (1.25,-0.5) to[out=0,in=270] (1.5,-0.25);
\draw[usual] (2,-0.25) to (2,-1.75);
\draw[usual] (2.5,0.25) to (2.5,1.75);
\draw[usual] (3,0.25) to (3,1.75);
\draw[usual] (2.5,-0.25) to[out=270,in=180] (2.75,-0.5) to[out=0,in=270] (3,-0.25);
\end{tikzpicture}
=
\begin{tikzpicture}[anchorbase]
\draw[usual] (0,2) node {$1$};
\draw[usual] (0.5,2) node {$1$};
\draw[usual] (1,2) node {$1$};
\draw[usual] (1.5,2) node {$2$};
\draw[usual] (2,2) node {$1$};
\draw[usual] (2.5,2) node {$1$};
\draw[usual] (3,2) node {$2$};
\draw[usual] (0,0) node {$1$};
\draw[usual] (0.5,0) node {$1$};
\draw[usual] (1,0) node {$1$};
\draw[usual] (1.5,0) node {$2$};
\draw[usual] (2,0) node {$1$};
\draw[usual] (2.5,0) node {$1$};
\draw[usual] (3,0) node {$2$};
\draw[usual] (0,0.25) to (0,1.75);
\draw[usual,dot] (0.5,1.75) to (0.5,1.55);
\draw[usual] (1,1.75) to (2,0.25);
\draw[usual] (2.5,1.75) to[out=270,in=180] (2.75,1.5) to[out=0,in=270] (3,1.75);
\end{tikzpicture}
.
\end{gather*}
\begin{gather*}
\begin{tikzpicture}[anchorbase]
\draw[usual] (0,2) node {$1$};
\draw[usual] (0.5,2) node {$2$};
\draw[usual] (1,2) node {$1$};
\draw[usual] (1.5,2) node {$2$};
\draw[usual] (2,2) node {$1$};
\draw[usual] (2.5,2) node {$2$};
\draw[usual] (3,2) node {$1$};
\draw[usual] (0,0) node {$1$};
\draw[usual] (0.5,0) node {$2$};
\draw[usual] (1,0) node {$1$};
\draw[usual] (1.5,0) node {$2$};
\draw[usual] (2,0) node {$1$};
\draw[usual] (2.5,0) node {$2$};
\draw[usual] (3,0) node {$1$};
\draw[usual] (0,-2) node {$1$};
\draw[usual] (0.5,-2) node {$2$};
\draw[usual] (1,-2) node {$1$};
\draw[usual] (1.5,-2) node {$2$};
\draw[usual] (2,-2) node {$1$};
\draw[usual] (2.5,-2) node {$2$};
\draw[usual] (3,-2) node {$1$};
\draw[usual] (0,1.75) to (0,0.25);
\draw[usual] (0,-0.25) to[out=270,in=180] (0.25,-0.5) to[out=0,in=270] (0.5,-0.25);
\draw[usual] (0.5,0.25) to[out=90,in=180] (0.75,0.5) to[out=0,in=90] (1,0.25);
\draw[usual] (1,-0.25) to[out=270,in=180] (1.25,-0.5) to[out=0,in=270] (1.5,-0.25);
\draw[usual] (1.5,0.25) to[out=90,in=180] (1.75,0.5) to[out=0,in=90] (2,0.25);
\draw[usual] (2,-0.25) to[out=270,in=180] (2.25,-0.5) to[out=0,in=270] (2.5,-0.25);
\draw[usual] (2.5,0.25) to[out=90,in=180] (2.75,0.5) to[out=0,in=90] (3,0.25);
\draw[usual] (3,-0.25) to (3,-1.75);
\end{tikzpicture}
=
\begin{tikzpicture}[anchorbase]
\draw[usual] (0,2) node {$1$};
\draw[usual] (0.5,2) node {$2$};
\draw[usual] (1,2) node {$1$};
\draw[usual] (1.5,2) node {$2$};
\draw[usual] (2,2) node {$1$};
\draw[usual] (2.5,2) node {$2$};
\draw[usual] (3,2) node {$1$};
\draw[usual] (0,0) node {$1$};
\draw[usual] (0.5,0) node {$2$};
\draw[usual] (1,0) node {$1$};
\draw[usual] (1.5,0) node {$2$};
\draw[usual] (2,0) node {$1$};
\draw[usual] (2.5,0) node {$2$};
\draw[usual] (3,0) node {$1$};
\draw[usual] (0,1.75) to (3,0.25);
\end{tikzpicture}
.
\end{gather*}
The general pattern is similar.
\end{proof}

\begin{Proposition}\label{P:rmocells}
The cells of $rMo_n$ have the following properties:
\begin{enumerate}
\item The $J$-cells of $rMo_n$ are uniquely determined by the sequence of through strands $\alpha(k)$. That is, each $J$-cell is of apex $\alpha(k)$, with a partial order given by the partial order on sequences.
\item Within each $J$-cell, sizes of the right and left cells can be found by counting the number of diagrams when fixing the top respectively bottom half. Explicitly, we have 
\begin{equation*}
|\mathcal{R}_n^{\alpha(k)}| = \mathlarger{\sum}_{\substack{j_1+\dots+j_{k+1} \\=2n-k }}\prod_{i=1}^{k+1}B_{j_j}^0
\end{equation*}
and 
\begin{equation*}
|\mathcal{L}_n^{\alpha(k)}| = \mathlarger{\sum}_{\substack{j_1+\dots+j_{k+1} \\=2n-k }}\prod_{i=1}^{k+1}T_0^{j_j}
\end{equation*}
where the sum is over the partition of through strands \autoref{rmopartition}, and $B_{j_i}^0$ and $T_0^{j_i}$ denote the number of possible diagrams on the block of the bottom respectively top object of size $j_i$.
\item Every $J$-cell of $rMo_n$ is idempotent.
\end{enumerate}
\end{Proposition}
\begin{proof}
Fix a $J$-cell $\mathcal{J}$, and suppose $a,b \in J$, meaning $a\leq_{lr} b$ and $b \leq_{lr} a$. This means $\exists c,c',d,d' \in rMo_n$ such that $b = cad$ and $a = c'bd'$.  

Using \autoref{factorisationRMO}, we have:
\begin{center}    
\begin{tikzpicture}[anchorbase,scale=1.5]
\draw[mor] (0,0) to (0.25,0.5) to (0.75,0.5) to (1,0) to (0,0);
\node at (0.5,0.25){$B_b$};
\draw[mor] (0,1.5) to (0.25,1) to (0.75,1) to (1,1.5) to (0,1.5);
\node at (0.5,1.25){$T_b$};
\draw[mor] (0.25,0.5) to (0.25,1) to (0.75,1) to (0.75,0.5) to (0.25,0.5);
\node at (0.5,0.75){\scalebox{0.75}{$1_{\alpha(k_b)}$}};
\end{tikzpicture}
=
\begin{tikzpicture}[anchorbase,scale=1,decoration={brace,amplitude=7pt}]
\draw[mor] (0,-2) to (0.25,-1.5) to (0.75,-1.5) to (1,-2) to (0,-2);
\node at (0.5,-1.75){$B_d$};
\draw[mor] (0,-0.5) to (0.25,-1) to (0.75,-1) to (1,-0.5) to (0,-0.5);
\node at (0.5,-0.75){$T_d$};
\draw[mor] (0.25,-1.5) to (0.25,-1) to (0.75,-1) to (0.75,-1.5) to (0.25,-1.5);
\node at (0.5,-1.25){\scalebox{0.5}{$1_{\alpha(k_d)}$}};
\draw[mor] (0,-0.5) to (0.25,0) to (0.75,0) to (1,-0.5) to (0,-0.5);
\node at (0.5,-0.25){$B_a$};
\draw[mor] (0,1) to (0.25,0.5) to (0.75,0.5) to (1,1) to (0,1);
\node at (0.5,0.75){$T_a$};
\draw[mor] (0.25,0) to (0.25,0.5) to (0.75,0.5) to (0.75,0) to (0.25,0);
\node at (0.5,0.25){\scalebox{0.5}{$1_{\alpha(k_a)}$}};
\draw[mor] (0,1) to (0.25,1.5) to (0.75,1.5) to (1,1) to (0,1);
\node at (0.5,1.25){$B_c$};
\draw[mor] (0,2.5) to (0.25,2) to (0.75,2) to (1,2.5) to (0,2.5);
\node at (0.5,2.25){$T_c$};
\draw[mor] (0.25,1.5) to (0.25,2) to (0.75,2) to (0.75,1.5) to (0.25,1.5);
\node at (0.5,1.75){\scalebox{0.5}{$1_{\alpha(k_c)}$}};
\draw[decorate] (1,2) to (1,-1.5);
\end{tikzpicture}
$=1_{\alpha(k_b)}$;
\begin{tikzpicture}[anchorbase,scale=1.5]
\draw[mor] (0,0) to (0.25,0.5) to (0.75,0.5) to (1,0) to (0,0);
\node at (0.5,0.25){$B_a$};
\draw[mor] (0,1.5) to (0.25,1) to (0.75,1) to (1,1.5) to (0,1.5);
\node at (0.5,1.25){$T_a$};
\draw[mor] (0.25,0.5) to (0.25,1) to (0.75,1) to (0.75,0.5) to (0.25,0.5);
\node at (0.5,0.75){\scalebox{0.75}{$1_{\alpha(k_a)}$}};
\end{tikzpicture}
=
\begin{tikzpicture}[anchorbase,scale=1,decoration={brace,amplitude=7pt}]
\draw[mor] (0,-2) to (0.25,-1.5) to (0.75,-1.5) to (1,-2) to (0,-2);
\node at (0.5,-1.75){$B_{d'}$};
\draw[mor] (0,-0.5) to (0.25,-1) to (0.75,-1) to (1,-0.5) to (0,-0.5);
\node at (0.5,-0.75){$T_{d'}$};
\draw[mor] (0.25,-1.5) to (0.25,-1) to (0.75,-1) to (0.75,-1.5) to (0.25,-1.5);
\node at (0.5,-1.25){\scalebox{0.5}{$1_{\alpha(k_{d'})}$}};
\draw[mor] (0,-0.5) to (0.25,0) to (0.75,0) to (1,-0.5) to (0,-0.5);
\node at (0.5,-0.25){$B_b$};
\draw[mor] (0,1) to (0.25,0.5) to (0.75,0.5) to (1,1) to (0,1);
\node at (0.5,0.75){$T_b$};
\draw[mor] (0.25,0) to (0.25,0.5) to (0.75,0.5) to (0.75,0) to (0.25,0);
\node at (0.5,0.25){\scalebox{0.5}{$1_{\alpha(k_b)}$}};
\draw[mor] (0,1) to (0.25,1.5) to (0.75,1.5) to (1,1) to (0,1);
\node at (0.5,1.25){$B_{c'}$};
\draw[mor] (0,2.5) to (0.25,2) to (0.75,2) to (1,2.5) to (0,2.5);
\node at (0.5,2.25){$T_{c'}$};
\draw[mor] (0.25,1.5) to (0.25,2) to (0.75,2) to (0.75,1.5) to (0.25,1.5);
\node at (0.5,1.75){\scalebox{0.5}{$1_{\alpha(k_{c'})}$}};
\draw[decorate] (1,2) to (1,-1.5);
\end{tikzpicture}
$=1_{\alpha(k_a)}$.
\end{center}
In the left equation, the middle of the right hand side must resolve to $1_{\alpha(k_b)}$, and since $1_{\alpha(k_a)}$ is contained within that middle, $\alpha(k_b) \geq \alpha(k_a)$ (using \autoref{rmoseqorder}. Similarly, $\alpha(k_a) \geq \alpha(k_b)$, and thus $\alpha(k_a)=\alpha(k_b)$. 

Clearly, if $\alpha(k_a)=\alpha(k_b)$, then $a \leq_{lr} b$ and $b \leq_{lr} a$, and hence the $J$-cells of $rMo_n$ are completely determined by the sequences of through strands. 

The argument for right and left cells is similar. If we have $a,b \in J$, with sequence of through strands $\alpha(k)$, then the condition for $a \sim_r b$ is:
\begin{center}
\begin{tikzpicture}[anchorbase,scale=1.5]
\draw[mor] (0,0) to (0.25,0.5) to (0.75,0.5) to (1,0) to (0,0);
\node at (0.5,0.25){$B_b$};
\draw[mor] (0,1.5) to (0.25,1) to (0.75,1) to (1,1.5) to (0,1.5);
\node at (0.5,1.25){$T_b$};
\draw[mor] (0.25,0.5) to (0.25,1) to (0.75,1) to (0.75,0.5) to (0.25,0.5);
\node at (0.5,0.75){\scalebox{0.5}{$1_{\alpha(k)}$}};
\end{tikzpicture}
=
\begin{tikzpicture}[anchorbase,scale=1]
\draw[mor] (0,-0.5) to (0.25,0) to (0.75,0) to (1,-0.5) to (0,-0.5);
\node at (0.5,-0.25){$B_c$};
\draw[mor] (0,1) to (0.25,0.5) to (0.75,0.5) to (1,1) to (0,1);
\node at (0.5,0.75){$T_c$};
\draw[mor] (0.25,0) to (0.25,0.5) to (0.75,0.5) to (0.75,0) to (0.25,0);
\node at (0.5,0.25){\scalebox{0.5}{$1_{\alpha(k_c)}$}};
\draw[mor] (0,1) to (0.25,1.5) to (0.75,1.5) to (1,1) to (0,1);
\node at (0.5,1.25){$B_a$};
\draw[mor] (0,2.5) to (0.25,2) to (0.75,2) to (1,2.5) to (0,2.5);
\node at (0.5,2.25){$T_a$};
\draw[mor] (0.25,1.5) to (0.25,2) to (0.75,2) to (0.75,1.5) to (0.25,1.5);
\node at (0.5,1.75){\scalebox{0.5}{$1_{\alpha(k)}$}};
\end{tikzpicture}
;
\begin{tikzpicture}[anchorbase,scale=1.5]
\draw[mor] (0,0) to (0.25,0.5) to (0.75,0.5) to (1,0) to (0,0);
\node at (0.5,0.25){$B_a$};
\draw[mor] (0,1.5) to (0.25,1) to (0.75,1) to (1,1.5) to (0,1.5);
\node at (0.5,1.25){$T_a$};
\draw[mor] (0.25,0.5) to (0.25,1) to (0.75,1) to (0.75,0.5) to (0.25,0.5);
\node at (0.5,0.75){\scalebox{0.5}{$1_{\alpha(k)}$}};
\end{tikzpicture}
=
\begin{tikzpicture}[anchorbase,scale=1]
\draw[mor] (0,-0.5) to (0.25,0) to (0.75,0) to (1,-0.5) to (0,-0.5);
\node at (0.5,-0.25){$B_{c'}$};
\draw[mor] (0,1) to (0.25,0.5) to (0.75,0.5) to (1,1) to (0,1);
\node at (0.5,0.75){$T_{c'}$};
\draw[mor] (0.25,0) to (0.25,0.5) to (0.75,0.5) to (0.75,0) to (0.25,0);
\node at (0.5,0.25){\scalebox{0.5}{$1_{\alpha(k_{c'})}$}};
\draw[mor] (0,1) to (0.25,1.5) to (0.75,1.5) to (1,1) to (0,1);
\node at (0.5,1.25){$B_b$};
\draw[mor] (0,2.5) to (0.25,2) to (0.75,2) to (1,2.5) to (0,2.5);
\node at (0.5,2.25){$T_b$};
\draw[mor] (0.25,1.5) to (0.25,2) to (0.75,2) to (0.75,1.5) to (0.25,1.5);
\node at (0.5,1.75){\scalebox{0.5}{$1_{\alpha(k)}$}};
\end{tikzpicture}
.
\end{center}
which shows that $a \sim_r b \iff T_a=T_b$, i.e. the top half of the diagrams is the same. Similar reasoning shows that $a\sim_l b \iff B_a = B_b$. 

Given this, we know that the through strands of a diagram partition both the top and bottom into blocks, with sizes given by \autoref{rmopartition}. For a fixed partition, these blocks are independent of each other, thus the number of diagrams is the product of the number of diagrams within each block. Taking the sum over all partitions gives the desired result.

Finally, if there are $k$ through strands, the $J$-cell of apex $\alpha(k)$ contains an idempotent that has pairs of dots everywhere there are no through strands (in general there are multiple of these for each sequence of through strands). For example, if $n=3$ and $\alpha(k) = 1\hspace{0.05cm}2\hspace{0.05cm}2$, then the following diagram is an idempotent in the corresponding $J$-cell:
\begin{center}
\begin{tikzpicture}
\draw[usual] (0,0) to (0,1);
\draw[usual,dot] (0.5,0) to (0.5,0.2);
\draw[usual,dot] (0.5,1) to (0.5,0.8);
\draw[usual,dot] (1,0) to (1,0.2);
\draw[usual,dot] (1,1) to (1,0.8);
\draw[usual] (1.5,0) to (1.5,1);
\draw[usual,dot] (2,0) to (2,0.2);
\draw[usual,dot] (2,1) to (2,0.8);
\draw[usual] (2.5,0) to (2.5,1);
\draw[usual] (0,1) node [above] {$1$};
\draw[usual] (0.5,1) node [above] {$2$};
\draw[usual] (1,1) node [above] {$1$};
\draw[usual] (1.5,1) node [above] {$2$};
\draw[usual] (2,1) node [above] {$1$};
\draw[usual] (2.5,1) node [above] {$2$};
\draw[usual] (0,0) node [below] {$1$};
\draw[usual] (0.5,0) node [below] {$2$};
\draw[usual] (1,0) node [below] {$1$};
\draw[usual] (1.5,0) node [below] {$2$};
\draw[usual] (2,0) node [below] {$1$};
\draw[usual] (2.5,0) node [below] {$2$};
\end{tikzpicture}
.
\end{center}
The general case can be verified similarly.
\end{proof}

\begin{figure}[ht]
\includegraphics[scale=0.76]{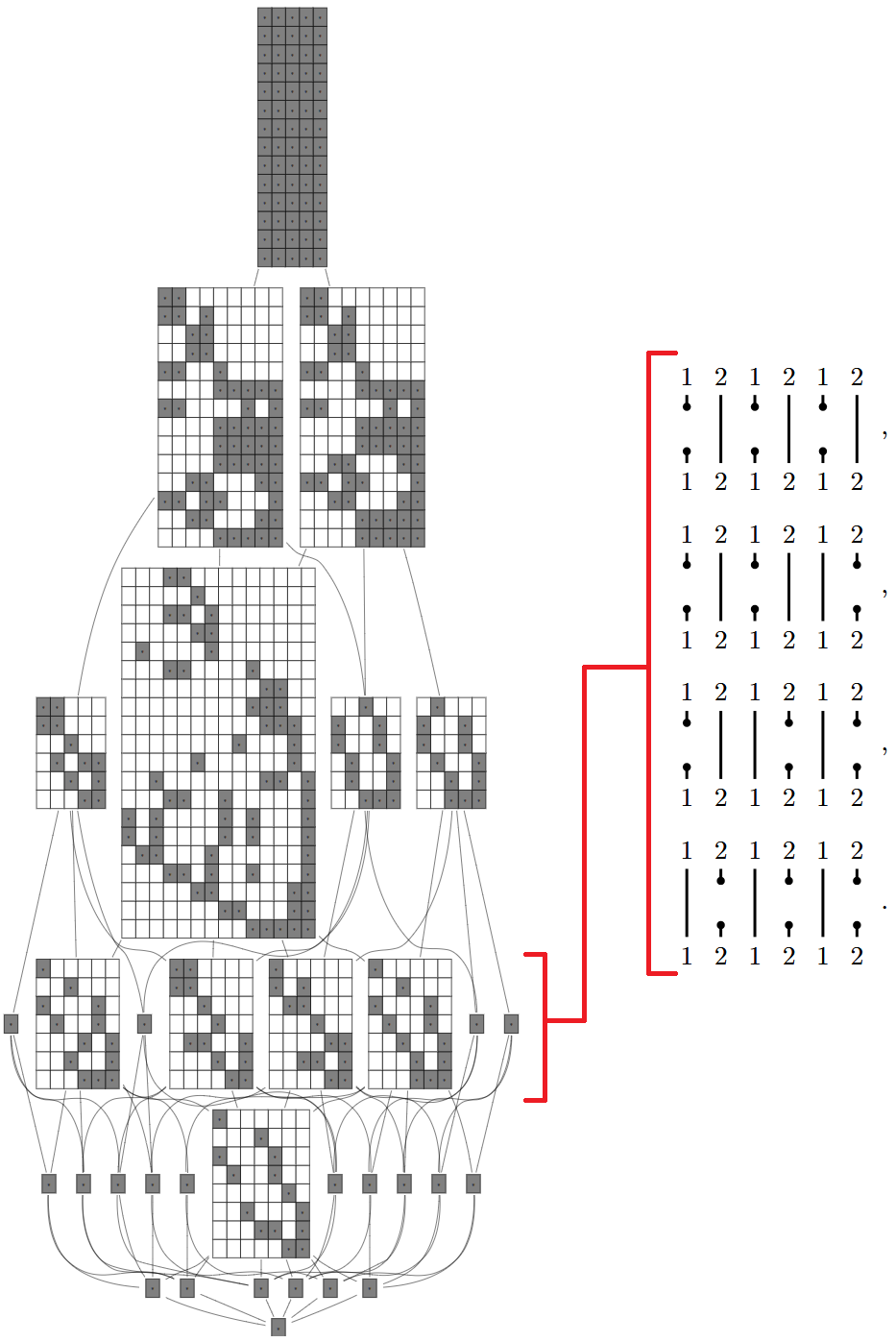}
\caption{Left: The cell structure for $rMo_3$. Right: explicit diagrams representing the single elements of the four $J$-cells of size $1$ for $k=3$.}
\label{rMo6_cells}
\end{figure}

\begin{Remark}
In general, we cannot compare all $J$-cells of $rMo_n$, so there is no total order like with $rTL_n$, only the partial order determined by the sequences of through strands: $J_{\alpha(k_1)} < J_{\alpha(k_2)} \iff \alpha(k_2)$ is a subsequence of $\alpha(k_1)$. As a result, we have a more complicated (or exciting) cell structure, seen in \autoref{rMo6_cells}.  
\end{Remark}

With this set up, we will continue to focus on the right cells, and the results for the left cells follow similarly. 

\begin{Proposition}\label{rmobottomdiag}
For all $m \in \mathbb{Z}_{\geq 0}$, the number of possible diagrams on the bottom object $1\hspace{0.05cm}2\hspace{0.05cm}1\hspace{0.05cm}2 \dots$ of length $m$ is given by:
\begin{equation}
A_m^0 = C\left(\left\lfloor \frac{m+1}{2} \right\rfloor \right)
\end{equation}
where $C(l)$ denotes the $l$th Catalan number. Moreover, if instead we have $2 \hspace{0.05cm} 1 \hspace{0.05cm} 2 \hspace{0.05cm} 1$ of length $m$, then the number of diagrams is $A_{m+2}^0$ if $m$ is even and $A_{m+1}$ if $m$ is odd.
\end{Proposition}

\begin{proof}
Let $m=2l$ be even, and note that $A_0^0 = 1$ (the only option is the empty diagram $1_{\mathds{1}}$). We claim that:
\begin{equation} \label{rmocatalanrec}
A_{2l}^0 = \sum_{i=1}^lA_{2i-2}^0A_{2l-2i}^0
\end{equation}
which is the familiar recurrence relation for the Catalan numbers, so assuming \autoref{rmocatalanrec}, we obtain $A_{2l}^0=C(l)=C(m/2)$. Observe also that if we consider $1\hspace{0.05cm}2\hspace{0.05cm}1\hspace{0.05cm}2\dots 1\hspace{0.05cm}2\hspace{0.05cm}1$ of length $2l-1$, adding a $2$ to the right changes nothing, since it cannot pair with anything to the left, and there is no $1$ to the right to pair with. Therefore, $A^0_{2l-1}=A^0_{2l}$ so we finally have $A_m^0=C\left(\left\lfloor \frac{m+1}{2} \right\rfloor\right)$. 

Now, to prove \autoref{rmocatalanrec}, consider the possibilities for the first $2$ in $1\hspace{0.05cm}2\hspace{0.05cm}1\hspace{0.05cm}2\dots 1\hspace{0.05cm}2$ (of length $2l$):
\begin{center}
\begin{tikzpicture}[anchorbase]
\draw[usual] (0,0) node {$1$};
\draw[usual] (0.25,0) node {$2$};
\draw[usual] (0.5,0) node {$1$};
\draw[usual] (0.75,0) node {$2$};
\draw[usual] (1.25,0) node [below] {$\dots$};
\draw[usual] (1.75,0) node {$1$};
\draw[usual] (2,0) node {$2$};
\draw[usual] (2.5,0) node [below] {$\dots$};
\draw[usual] (3,0) node {$1$};
\draw[usual] (3.25,0) node {$2$};
\draw[usual] (3.5,0) node {$1$};
\draw[usual] (3.75,0) node {$2$};
\draw[usual] (0.25,0.25) to[out=90,in=180] (0.75,0.6) to[out=0,in=90] (1.75,0.25);
\end{tikzpicture}
\end{center}
If we pair the first $2$ with the $i$th $1$ following that $2$, we have $1\hspace{0.05cm}2\dots 1\hspace{0.05cm}2$ of length $2i-2$ under the cap, and $1\hspace{0.05cm}2\dots 1\hspace{0.05cm}2$ of length $2l-2i$ remaining outside the cap, noting that the interruption where the cap is placed can be ignored, since the diagram is planar and thus we treat inside and outside the cap separately. These can be arranged independently in $A_{2i-2}^0A_{2l-2i}^0$ ways. Setting $i=l$ for the case where the first $2$ remains unpaired, we thus conclude the recurrence relation \autoref{rmocatalanrec}.

Finally, if we have $1\hspace{0.05cm}2\hspace{0.05cm}1\hspace{0.05cm}2\dots 1\hspace{0.05cm}2\hspace{0.05cm}$ of length $m+2$, the first $1$ and the final $2$ must remain unpaired, so the number of diagrams $A_{m+2}^0$ on this object is the same as the number of diagrams on $2\hspace{0.05cm}1\hspace{0.05cm}2\hspace{0.05cm}1\dots 2\hspace{0.05cm}1$ of length $m$. In the case $2\hspace{0.05cm}1\hspace{0.05cm}2\hspace{0.05cm}1\dots 2\hspace{0.05cm}1\hspace{0.05cm}2$ of length $m$ (i.e. odd), the final $2$ must remain unpaired, so the number of ways is the same as $2\hspace{0.05cm}1\hspace{0.05cm}2\hspace{0.05cm}1\dots 2\hspace{0.05cm}1$ of length $m-1$, which is $A_{m+1}^0$.   
\end{proof}

Consider now the following cases for $\alpha(k)$:
\begin{gather*}
\begin{tikzpicture}[anchorbase]
\draw[usual] (0,0) node {$1$};
\draw[usual] (0.25,0) node {$2$};
\draw[usual] (0.5,0) node {$1$};
\draw[usual] (0.75,0) node {$2$};
\draw[usual] (1.25,0) node [below] {$\dots$};
\draw[usual] (1.75,0) node {$1$};
\draw[usual] (2,0) node {$2$};
\draw[usual] (2.5,0) node [below] {$\dots$};
\draw[usual] (3,0) node {$1$};
\draw[usual] (3.25,0) node {$2$};
\draw[usual] (3.5,0) node {$1$};
\draw[usual] (3.75,0) node {$2$};
\draw[usual] (1.875,1.5) node {$\alpha(k)_1=1$};
\draw[usual] (0.5,0.25) to (1,1);
\draw[usual] (0,-2) node {$1$};
\draw[usual] (0.25,-2) node {$2$};
\draw[usual] (0.5,-2) node {$1$};
\draw[usual] (0.75,-2) node {$2$};
\draw[usual] (1.25,-2) node [below] {$\dots$};
\draw[usual] (1.75,-2) node {$1$};
\draw[usual] (2,-2) node {$2$};
\draw[usual] (2.5,-2) node [below] {$\dots$};
\draw[usual] (3,-2) node {$1$};
\draw[usual] (3.25,-2) node {$2$};
\draw[usual] (3.5,-2) node {$1$};
\draw[usual] (3.75,-2) node {$2$};
\draw[usual] (1.875,-1) node [above] {$\alpha(k)_{k}=2$};
\draw[usual] (3.25,-1.75) to (2.75,-1);
\draw[usual] (4,1) node {,}; 
\draw[usual] (4,-1) node {,}; 
\end{tikzpicture}
\begin{tikzpicture}[anchorbase]
\draw[usual] (0,0) node {$1$};
\draw[usual] (0.25,0) node {$2$};
\draw[usual] (0.5,0) node {$1$};
\draw[usual] (0.75,0) node {$2$};
\draw[usual] (1.25,0) node [below] {$\dots$};
\draw[usual] (1.75,0) node {$1$};
\draw[usual] (2,0) node {$2$};
\draw[usual] (2.5,0) node [below] {$\dots$};
\draw[usual] (3,0) node {$1$};
\draw[usual] (3.25,0) node {$2$};
\draw[usual] (3.5,0) node {$1$};
\draw[usual] (3.75,0) node {$2$};
\draw[usual] (1.875,1.5) node {$\alpha(k)_1=2$};
\draw[usual] (0.75,0.25) to (1.25,1);
\draw[usual] (0,-2) node {$1$};
\draw[usual] (0.25,-2) node {$2$};
\draw[usual] (0.5,-2) node {$1$};
\draw[usual] (0.75,-2) node {$2$};
\draw[usual] (1.25,-2) node [below] {$\dots$};
\draw[usual] (1.75,-2) node {$1$};
\draw[usual] (2,-2) node {$2$};
\draw[usual] (2.5,-2) node [below] {$\dots$};
\draw[usual] (3,-2) node {$1$};
\draw[usual] (3.25,-2) node {$2$};
\draw[usual] (3.5,-2) node {$1$};
\draw[usual] (3.75,-2) node {$2$};
\draw[usual] (1.875,-1) node [above] {$\alpha(k)_{k}=1$};
\draw[usual] (3,-1.75) to (2.5,-1);
\draw[usual] (4,1) node {;}; 
\draw[usual] (4,-1) node {;}; 
\end{tikzpicture}
\end{gather*}
\begin{gather*}
\begin{tikzpicture}[anchorbase]
\draw[usual] (0,0) node {$1$};
\draw[usual] (0.25,0) node {$2$};
\draw[usual] (0.5,0) node {$1$};
\draw[usual] (0.75,0) node {$2$};
\draw[usual] (1.25,0) node [below] {$\dots$};
\draw[usual] (1.75,0) node {$1$};
\draw[usual] (2,0) node {$2$};
\draw[usual] (2.5,0) node [below] {$\dots$};
\draw[usual] (3,0) node {$1$};
\draw[usual] (3.25,0) node {$2$};
\draw[usual] (3.5,0) node {$1$};
\draw[usual] (3.75,0) node {$2$};
\draw[usual] (1.875,1.5) node {$\alpha(k)_{i-1}\alpha(k)_i=1\hspace{0.05cm}2$};
\draw[usual] (0.5,0.25) to (1,1);
\draw[usual] (3.25,0.25) to (2.75,1);
\draw[usual] (0,-2) node {$1$};
\draw[usual] (0.25,-2) node {$2$};
\draw[usual] (0.5,-2) node {$1$};
\draw[usual] (0.75,-2) node {$2$};
\draw[usual] (1.25,-2) node [below] {$\dots$};
\draw[usual] (1.75,-2) node {$1$};
\draw[usual] (2,-2) node {$2$};
\draw[usual] (2.5,-2) node [below] {$\dots$};
\draw[usual] (3,-2) node {$1$};
\draw[usual] (3.25,-2) node {$2$};
\draw[usual] (3.5,-2) node {$1$};
\draw[usual] (3.75,-2) node {$2$};
\draw[usual] (1.875,-1) node [above] {$\alpha(k)_{i-1}\alpha(k)_i=1\hspace{0.05cm}1$};
\draw[usual] (0.5,-1.75) to (1,-1);
\draw[usual] (3.5,-1.75) to (3,-1);
\draw[usual] (4,1) node {,}; 
\draw[usual] (4,-1) node {,}; 
\end{tikzpicture}
\begin{tikzpicture}[anchorbase]
\draw[usual] (0,0) node {$1$};
\draw[usual] (0.25,0) node {$2$};
\draw[usual] (0.5,0) node {$1$};
\draw[usual] (0.75,0) node {$2$};
\draw[usual] (1.25,0) node [below] {$\dots$};
\draw[usual] (1.75,0) node {$1$};
\draw[usual] (2,0) node {$2$};
\draw[usual] (2.5,0) node [below] {$\dots$};
\draw[usual] (3,0) node {$1$};
\draw[usual] (3.25,0) node {$2$};
\draw[usual] (3.5,0) node {$1$};
\draw[usual] (3.75,0) node {$2$};
\draw[usual] (1.875,1.5) node {$\alpha(k)_{i-1}\alpha(k)_i=2\hspace{0.05cm}1$};
\draw[usual] (0.75,0.25) to (1.25,1);
\draw[usual] (3,0.25) to (2.5,1);
\draw[usual] (0,-2) node {$1$};
\draw[usual] (0.25,-2) node {$2$};
\draw[usual] (0.5,-2) node {$1$};
\draw[usual] (0.75,-2) node {$2$};
\draw[usual] (1.25,-2) node [below] {$\dots$};
\draw[usual] (1.75,-2) node {$1$};
\draw[usual] (2,-2) node {$2$};
\draw[usual] (2.5,-2) node [below] {$\dots$};
\draw[usual] (3,-2) node {$1$};
\draw[usual] (3.25,-2) node {$2$};
\draw[usual] (3.5,-2) node {$1$};
\draw[usual] (3.75,-2) node {$2$};
\draw[usual] (1.875,-1) node [above] {$\alpha(k)_{i-1}\alpha(k)_i=2\hspace{0.05cm}2$};
\draw[usual] (3.25,-1.75) to (2.75,-1);
\draw[usual] (0.75,-1.75) to (1.25,-1);
\draw[usual] (4,1) node {;}; 
\draw[usual] (4,-1) node {.}; 
\end{tikzpicture}
\end{gather*}
where $2 \leq i \leq k$. Recall from \autoref{rmopartition} that the through strands partition the bottom object into $k+1$ blocks of size $\geq 0$. Each block is of size $j_i$ and corresponds to the block to the left of the $i$th through strand, except for block $k+1$, which corresponds to the block to the right of the $k$th through strand. We summarise the above cases for each block, using \autoref{rmobottomdiag}, in the following table:
\begin{gather}\label{bottomtable}
\begin{tabular}{c|c|c|c}
\textbf{Case} & \textbf{Block} & \textbf{Block Shape} & \textbf{No. Diagrams} \\
\hline
\hline
$\alpha(k)_1=1$ & $j_1$ even, $j_1 \geq 0$ & $1\hspace{0.05cm}2\dots 1\hspace{0.05cm} 2\hspace{0.05cm}$ & $C(j_1/2)$ \\
\hline
$\alpha(k)_1=2$ & $j_1$ odd, $j_1\geq 1$ & $1\hspace{0.05cm}2\dots 1\hspace{0.05cm} 2\hspace{0.05cm}1$ & $C((j_1+1)/2)$ \\
\hline
$\alpha(k)_{k}=2$ & $j_{k+1}$ even, $j_{k+1} \geq 0$ & $1\hspace{0.05cm}2\dots 1\hspace{0.05cm} 2\hspace{0.05cm}$ & $C(j_{k+1}/2)$ \\
\hline
$\alpha(k)_{k}=1$ & $j_{k+1}$ odd, $j_{k+1} \geq 1$ & $2\hspace{0.05cm}1\dots 2\hspace{0.05cm} 1\hspace{0.05cm}2$ & $C((j_{k+1}+1)/2)$ \\
\hline
$\alpha(k)_{i-1}\alpha(k)_i=12$ & $j_i$ even, $j_i \geq 0$ & $2\hspace{0.05cm}1\dots 2\hspace{0.05cm} 1\hspace{0.05cm}$ & $C(j_i/2+1)$ \\
\hline
$\alpha(k)_{i-1}\alpha(k)_i=21$ & $j_i$ even, $j_i \geq 0$ & $1\hspace{0.05cm}2\dots 1\hspace{0.05cm} 2\hspace{0.05cm}$ & $C(j_i/2)$ \\
\hline
$\alpha(k)_{i-1}\alpha(k)_i=11$ & $j_i$ odd, $j_i \geq 1$ & $2\hspace{0.05cm}1\dots 2\hspace{0.05cm} 1\hspace{0.05cm}2$ & $C((j_i+1)/2)$ \\
\hline
$\alpha(k)_{i-1}\alpha(k)_i=22$ & $j_i$ odd, $j_i \geq 1$ & $1\hspace{0.05cm}2\dots 1\hspace{0.05cm} 2\hspace{0.05cm}1$ & $C((j_i+1)/2)$ 
\end{tabular}
\end{gather}
Observe that for all even $j_i$, except if $\alpha(k)_{i-1}\alpha(k)_i=12$, the number of diagrams is equal to $C(j_i/2)$, and for all odd $j_i$, the number of diagrams is equal to $C((j_i+1)/2)$. For this reason, we need the following lemma.

\begin{Lemma}\label{rmosamesequences}
Let $\alpha(k)$ be a sequence of $k$ through strands with $j$ occurrences of the block $1\hspace{0.05cm}2$. Then, we have:
\begin{enumerate}
\item all such sequences $\beta(k)$ with $j$ occurrences of $1\hspace{0.05cm}2$ have the same numbers of odd and even blocks, and
\item the number of such sequences is equal to $\binom{k+1}{2j+1}$.
\end{enumerate}
\end{Lemma}

\begin{proof}
For part (a), we begin with a sequence $\alpha(k)$ with $j$ $1\hspace{0.05cm}2$ blocks, then we show that if we create a new sequence $\beta(k)$ by changing a $1$ into a $2$ or a $2$ into a $1$, we either preserve the counts outlined in the proposition, or the number of $1\hspace{0.05cm}2$ blocks changes. We have the following cases:
\begin{enumerate}
\item The first number is a $1$ and $\alpha(k)=\textbf{1}\hspace{0.05cm}1\dots 1 \hspace{0.05cm} 2 \hspace{0.05cm} \alpha(k)_{i\geq m+1}$ where we have only considered up to the first appearance of a $2$ at $i=m$. Then if we set $\beta(k) = \textbf{2}\hspace{0.05cm}1\dots 1 \hspace{0.05cm} 2 \hspace{0.05cm}\alpha(k)_{i\geq m+1}$, the number of $1\hspace{0.05cm}2$ blocks hasn't changed, the first block has changed from even to odd, and the second block has changed from odd to even. Hence, the required counts remain the same. 
\item The first number is a $2$ and $\alpha(k)=\textbf{2}\hspace{0.05cm}2\dots 2 \hspace{0.05cm}1 \hspace{0.05cm}1\dots 1 \hspace{0.05cm} 2 \hspace{0.05cm} \alpha(k)_{i\geq m+1}$. 

Then, if $\beta(k)=\textbf{1}\hspace{0.05cm}2\dots 2 \hspace{0.05cm}1 \dots 1 \hspace{0.05cm} 2 \hspace{0.05cm} \alpha(k)_{i\geq m+1}$ we have increased the number of $1 \hspace{0.05cm}2$ blocks.
\item The last number is a $1$ and in $\beta(k)$ we change this $1$ into a $2$. Then there are two cases:
\begin{enumerate}
\item $\alpha(k) = \alpha(k)_{i\leq k-2}\hspace{0.05cm}2\hspace{0.05cm}\textbf{1}$, then if $\beta(k) = \alpha(k)_{i\leq k-2}\hspace{0.05cm}2\hspace{0.05cm}\textbf{2}$, we have changed the final block from odd to even and changed the second last block from even to odd, leaving the rest unchanged. 
\item $\alpha(k) = \alpha(k)_{i\leq k-2}\hspace{0.05cm}1\hspace{0.05cm}\textbf{1}$, then if $\beta(k) = \alpha(k)_{i\leq k-2}\hspace{0.05cm}1\hspace{0.05cm}\textbf{2}$, we have increased the number of $1\hspace{0.05cm}2$ blocks.
\end{enumerate}
\item The last number is a $2$ and in $\beta(k)$ we change this $2$ into a $1$. Then there are two cases:
\begin{enumerate}
\item $\alpha(k) = \alpha(k)_{i\leq k-2}\hspace{0.05cm}2\hspace{0.05cm}\textbf{2}$, then if $\beta(k) = \alpha(k)_{i\leq k-2}\hspace{0.05cm}2\hspace{0.05cm}\textbf{1}$, we have changed the final block from even to odd and changed the second last block from odd to even, so the overall counts remain unchanged. 
\item $\alpha(k) = \alpha(k)_{i\leq k-2}\hspace{0.05cm}1\hspace{0.05cm}\textbf{2}$, then if $\beta(k) = \alpha(k)_{i\leq k-2}\hspace{0.05cm}1\hspace{0.05cm}\textbf{1}$, we have decreased the number of $1\hspace{0.05cm}2$ blocks.
\end{enumerate}
\item $\alpha(k) = \alpha(k)_{i\leq m-1} \textbf{1} \alpha(k)_{i\geq m+1}$, and $\beta(k) = \alpha(k)_{i\leq m-1} \textbf{2} \alpha(k)_{i\geq m+1}$. Then, there are four cases:
\begin{enumerate}
\item $\alpha(k) = \alpha(k)_{i\leq m-2} 1\hspace{0.05cm}\textbf{1}\hspace{0.05cm}1 \alpha(k)_{i\geq m+2}$. If $\beta(k) = \alpha(k)_{i\leq m-2} 1\hspace{0.05cm}\textbf{2}\hspace{0.05cm}1 \alpha(k)_{i\geq m+2}$, then we have increased the number of $1\hspace{0.05cm}2$ blocks.
\item $\alpha(k) = \alpha(k)_{i\leq m-2} 1\hspace{0.05cm}\textbf{1}\hspace{0.05cm}2 \alpha(k)_{i\geq m+2}$. If $\beta(k) = \alpha(k)_{i\leq m-2} 1\hspace{0.05cm}\textbf{2}\hspace{0.05cm}2 \alpha(k)_{i\geq m+2}$, then we have lost a $1\hspace{0.05cm}2$ block in the the $m+1$st position and an odd block in the $m-1$st position, but gained a $1\hspace{0.05cm}2$ block in the $m-1$st position and an odd block in the $m+1$st position, so the overall counts remain unchanged.
\item $\alpha(k) = \alpha(k)_{i\leq m-2} 2\hspace{0.05cm}\textbf{1}\hspace{0.05cm}1 \alpha(k)_{i\geq m+2}$. If $\beta(k) = \alpha(k)_{i\leq m-2} 2\hspace{0.05cm}\textbf{2}\hspace{0.05cm}1 \alpha(k)_{i\geq m+2}$, then we have lost an even block in the $m-1$st position and an odd block in the $m+1$st position, but gained an odd block in the $m-1$st position and an even block in the $m+1$st position, so the overall counts remain unchanged.
\item $\alpha(k) = \alpha(k)_{i\leq m-2} 2\hspace{0.05cm}\textbf{1}\hspace{0.05cm}2 \alpha(k)_{i\geq m+2}$. If $\beta(k) = \alpha(k)_{i\leq m-2} 2\hspace{0.05cm}\textbf{2}\hspace{0.05cm}2 \alpha(k)_{i\geq m+2}$, then we have decreased the number of $1\hspace{0.05cm}2$ blocks.
\end{enumerate}
\item $\alpha(k) = \alpha(k)_{i\leq m-1} \textbf{2} \alpha(k)_{i\geq m+1}$, and $\beta(k) = \alpha(k)_{i\leq m-1} \textbf{1} \alpha(k)_{i\geq m+1}$. Then, there are four cases:
\begin{enumerate}
\item $\alpha(k) = \alpha(k)_{i\leq m-2} 1\hspace{0.05cm}\textbf{2}\hspace{0.05cm}1 \alpha(k)_{i\geq m+2}$. If $\beta(k) = \alpha(k)_{i\leq m-2} 1\hspace{0.05cm}\textbf{1}\hspace{0.05cm}1 \alpha(k)_{i\geq m+2}$, then we have decreased the number of $1\hspace{0.05cm}2$ blocks.
\item $\alpha(k) = \alpha(k)_{i\leq m-2} 1\hspace{0.05cm}\textbf{2}\hspace{0.05cm}2 \alpha(k)_{i\geq m+2}$. If $\beta(k) = \alpha(k)_{i\leq m-2} 1\hspace{0.05cm}\textbf{1}\hspace{0.05cm}2 \alpha(k)_{i\geq m+2}$, then we have lost a $1\hspace{0.05cm}2$ block in the $m-1$st position and an odd block in the $m+1$st position, but gained an odd block in the $m-1$st position and a $1\hspace{0.05cm}2$ block in the $m+1$st position, so the overall counts remain unchanged.
\item $\alpha(k) = \alpha(k)_{i\leq m-2} 2\hspace{0.05cm}\textbf{2}\hspace{0.05cm}1 \alpha(k)_{i\geq m+2}$. If $\beta(k) = \alpha(k)_{i\leq m-2} 2\hspace{0.05cm}\textbf{1}\hspace{0.05cm}1 \alpha(k)_{i\geq m+2}$, then we have lost an odd block in the $m-1$st position and an even block in the $m+1$st position, but gained an even block in the $m-1$st position and an odd block in the $m+1$st position, so the overall counts remain unchanged.
\item $\alpha(k) = \alpha(k)_{i\leq m-2} 2\hspace{0.05cm}\textbf{2}\hspace{0.05cm}2 \alpha(k)_{i\geq m+2}$. If $\beta(k) = \alpha(k)_{i\leq m-2} 2\hspace{0.05cm}\textbf{1}\hspace{0.05cm}2 \alpha(k)_{i\geq m+2}$, then we have increased the number of $1\hspace{0.05cm}2$ blocks.
\end{enumerate}
\end{enumerate}
Hence, since this covers every case for changing a $1$ into a $2$ or a $2$ into a $1$ for $\alpha(k)$, we have shown that the numbers of odd and even blocks remain unchanged if the number of $1\hspace{0.05cm}2$ blocks remains unchanged. 

Now, for part (b), consider a sequence $\alpha(k)$ of length $k$, made only of $1$s and $2$s, with $j$ occurrences of the substring $1\hspace{0.05cm}2$. Firstly, observe that if $j=0$, the number of possible sequences is clearly $k+1 = \binom{k+1}{1}$, since it's the number of sequences of the form $2\hspace{0.05cm}2\dots 2 \hspace{0.05cm}1\hspace{0.05cm}1\dots 1$ (including the sequences of all $2$s and all $1$s). 

Now suppose $1\leq j \leq \left\lfloor\frac{k}{2}\right\rfloor$. Suppose the first $1$ is in the $i$th position, $1\leq i \leq k-2j+1$, and suppose there is a sequence of $1$s of length $l$ (including the first $1$), $1\leq l \leq k-i-2j+1$:
\begin{equation*}
\overbrace{\underbrace{2\hspace{0.05cm}2\dots2}_{i-i}\hspace{0.05cm}\underbrace{1\hspace{0.05cm}1\dots1\hspace{0.05cm}1}_{l}\hspace{0.05cm}\underbrace{***}_{k-i-l+1}}^{k}
\end{equation*}
Then, we are placing $2j-1$ $1$s and $2$s (including the first $2$ after the $1$s) somewhere inside the remaining sequence of length $k-i-l+1$, which can be done in $\binom{k-i-l+1}{2j-1}$ ways. 

Thus, the total number of ways is $\sum_{i=1}^{k-2j+1}\sum_{l=1}^{k-i-2j+2}\binom{k-i-l+1}{2j-1}$. Using the Hockey-Stick identity, we obtain:
\[\sum_{i=1}^{k-2j+1}\sum_{l=1}^{k-i-2j+2}\binom{k-i-l+1}{2j-1}\] 
\[= \sum_{i=1}^{k-2j+1}\binom{k-i+1}{2j}\]
\[=\binom{k+1}{2j+1}\]
as required.
\end{proof}

\begin{Remark}\label{rmojincrease}
\autoref{rmosamesequences} shows that to find the sizes of right cells, we only need to consider the sequence of the form $1\hspace{0.05cm}2\hspace{0.05cm}1\hspace{0.05cm}2\dots 1\hspace{0.05cm}2 \hspace{0.05cm}2 \dots 2$, for each number of $1\hspace{0.05cm}2$ blocks $j$, $0\leq j \leq \left\lfloor\frac{k}{2}\right\rfloor$. Furthermore, we see that the $1\hspace{0.05cm}2$ blocks have $C(j_i/2+1)$ possible diagrams, which is larger than any of the other blocks, and the number of different partitions \autoref{rmopartition} also increases as $j$ increases. We prove the latter fact when discussing the planar rook monoid below (see \autoref{P:rprocellsizes}). Thus, the total number of possible diagrams increases as $j$ increases, with a minimum at $j=0$. This means that when finding the semisimple RepGap, we will need to focus on the $j=0$ case.
\end{Remark}

We denote by $\mathcal{R}_n^{k,j}$ the right cell arising from this particular sequence.

\begin{Proposition}\label{P:rmocellsizes}
Fix $n \in \mathbb{Z}_{\geq 0}$, $0\leq k \leq 2n$, and sequence of through strands $\alpha(k)$ with $j$ $1\hspace{0.05cm}2$ blocks, $0\leq j \leq \left\lfloor\frac{k}{2}\right\rfloor$. Then, the size of the corresponding right cell is:
\begin{equation}
|\mathcal{R}_n^{k,j}| = \sum_{r=0}^{k-j}(-1)^r\frac{k-r+1}{2n+k-r+1}\binom{k-j}{r}\binom{2n+k-r+1}{n}
\end{equation}
the size of the corresponding left cell is:
\begin{equation}
|\mathcal{L}_n^{k,j}| = \sum_{r=0}^{k-j+1}(-1)^r\frac{k-r+1}{2n+k-r+3}\binom{k-j+1}{r}\binom{2n+k-r+3}{n+1}
\end{equation}
and the size of the monoid is:
\begin{equation}
|rMo_n| = \sum_{k=0}^{2n}\sum_{j=0}^{\left\lfloor\frac{k}{2}\right\rfloor}\binom{k+1}{2j+1}|\mathcal{R}_n^{k,j}||\mathcal{L}_n^{k,j}|.
\end{equation}
\end{Proposition}

\begin{proof}
We have the sequence 
\begin{equation*}
\alpha(k) = \underbrace{1\hspace{0.05cm}2\hspace{0.05cm}1\hspace{0.05cm}2\dots 1\hspace{0.05cm}2}_{2j}\hspace{0.05cm}\underbrace{2\hspace{0.05cm}2\dots2}_{k-2j}
\end{equation*}
creating an ordered partition of the bottom object, with sizes $j_1,j_2,\dots,j_{k+1}$ such that
\begin{equation*}
j_1+j_2+\dots j_{k+1} = 2n-k.
\end{equation*}
With reference to the table \autoref{bottomtable}, the blocks of the partition are as follows:
\begin{enumerate}
\item $j_1$ is even, $j_1 \geq 0$, and the number of possible diagrams in the block is $C(j_1/2)$
\item $j_{k+1}$ is even, $j_{k+1} \geq 0$, and the number of possible diagrams in the block is $C(j_{k+1}/2)$
\item $j_{2l}$ is even, $j_{2l} \geq 0$, and the number of possible diagrams in the block is $C(j_{2l}/2+1)$ for $l=1,2,\dots,j$
\item  $j_{2m+1}$ is even, $j_{2m+1} \geq 0$, and the number of possible diagrams in the block is $C(j_{2m+1}/2)$ for $m=1,2,\dots,j-1$
\item $j_p$ is odd, $j_p \geq 1$, and the number of possible diagrams in the block is $C((j_p+1)/2)$ for $p=2j+1,2j+2,\dots,k$
\end{enumerate}
Hence, using \autoref{P:rmocells}, we obtain
\[
|\mathcal{R}_n^{k,j}|=\mathlarger{\sum}_{\substack{j_1+\dots+j_{k+1} \\=2n-k }} C(j_1/2)C(j_2/2+1)C(j_3/2)\dots C(j_{2j}/2+1)C((j_{2j+1}+1)/2)\dots C((j_k+1)/2)C(j_{k+1}/2).
\]
Now set $x_1=\frac{j_1}{2}$, $x_{2l}=\frac{j_{2l}}{2}+1$ for $l=1,\dots,j$, $x_{2m+1} = \frac{j_{2m+1}}{2}$ for $m=1,\dots,j-1$, $x_p=\frac{j_p+1}{2}$ for $p=2j+1,\dots,k$, and $x_{k+1}=\frac{x_{k+1}}{2}$. Then, the above becomes
\[
|\mathcal{R}_n^{k,j}|=\mathlarger{\sum}_{x_1+\dots+x_{k+1} =n}C(x_1)C(x_2)\dots C(x_{k+1}).
\]
This is almost the definition of the Catalan $k$-fold convolution, however we have $x_{2l}\geq 1$ and $x_p \geq 1$ for $l=1,\dots,j$ and $p=2j+1,\dots k$, so we instead need an alternating sum of Catalan $k$-fold convolutions (noting that $C(0)=1$):
\begin{equation*}
\begin{split}
|\mathcal{R}_n^{k,j}|=\mathlarger{\sum}_{\substack{x_1+\dots+x_{k+1} =n \\ x_i \geq 0}} C(x_1)\dots C(x_{k+1}) & - \binom{k-j}{1}\mathlarger{\sum}_{\substack{x_1+\dots+x_{k+1} =n \\ x_i \geq 0 \\ x_2=0}} C(x_1)\dots C(x_{k+1}) \\
+ \dots & + (-1)^{k-j} \binom{k-j}{k-j}\mathlarger{\sum}_{\substack{x_1+\dots+x_{k+1} =n \\ x_i \geq 0 \\ x_{2l},x_p=0}} C(x_1)\dots C(x_{k+1}) 
\end{split}
\end{equation*}
where in each term we have set $r$ of the $k-j$ $x_i$ that are $\geq 1$ to zero, noting that the formula will remain the same regardless of which we choose, so we can simply multiply by a binomial coefficient. From \cite[Proposition 1.2]{LS-catalan-convolution}, the formula for the Catalan $k$-fold convolution is
\[
\mathlarger{\sum}_{\substack{m_1+\dots m_k = n \\m_i \geq 0}}C(m_1)\dots C(m_k) = \frac{k}{2n+k}\binom{2n+k}{n}
\]
which finally gives
\[
|\mathcal{R}_n^{k,j}| = \sum_{r=0}^{k-j}(-1)^r\frac{k-r+1}{2n+k-r+1}\binom{k-j}{r}\binom{2n+k-r+1}{n}
\]
as required. The argument for left cells is similar. 

Finally, observe that if we fix the sequence $\alpha(k)$, the top half and bottom half of the diagrams become independent of each other. Thus, the number of possible diagrams in the monoid for each sequence $\alpha(k)$ is the number of top half diagrams multiplied by the number of bottom half diagrams, then to find the total $|rMo_n|$ by taking the sum over all sequences $\alpha(k)$, then another sum over all $k$, $0\leq k \leq 2n$. 

Using \autoref{rmosamesequences}, we can then simplify this to give
\[
|rMo_n| = \sum_{k=0}^{2n}\sum_{j=0}^{\left\lfloor\frac{k}{2}\right\rfloor}\binom{k+1}{2j+1}|\mathcal{R}_n^{k,j}||\mathcal{L}_n^{k,j}|.
\]
The proof is complete.
\end{proof}

\subsection{Semisimple RepGap and gap ratio of \texorpdfstring{$rMo_n$}{rMon}}

As we have seen above, although $rMo_n$ can largely be treated in the same way as $Mo_n$, certain differences make $rMo_n$ more intricate. This trend will continue in this section.

\begin{Remark}\label{R:SSRepGapWorse}
By \autoref{L:CellsDimensions}, the semisimple RepGap is an upper bound for the RepGap, and we will see that this upper bound alone is enough to determine that $rMo_n$ is less suitable for cryptography than $Mo_n$. The following example exhibits a strict inequality at individual apexes. 

We can compute the dimensions of the simple representations by computing the ranks of corresponding Gram matrices (using \autoref{P:CellsGram}). Setting $n=3$ and $k=2$, we obtain ranks $4,5,5$, and $12$, in increasing order of dimension (using GAP). The semisimple dimension can be easily computed via \autoref{P:rmocellsizes}, which gives $5$ for $j=0$ and $14$ for $j=1$.

This shows that the simple dimension can be strictly smaller than the semisimple dimension at an individual apex, and it also reveals that there are further distinguishing features of the sequence of through strands apart from simply the number of $1\hspace{0.05cm}2$ pairs when we are at the level of simple representations. In this case, the minimal simple dimension corresponds to the sequence $2\hspace{0.05cm}1$. 
\end{Remark}

We simplify the cell sizes found in \autoref{P:rmocellsizes}.

\begin{Lemma}\label{C:rmocellsizes}
The cell sizes are:
\[
|\mathcal{R}_n^{k,j}| = \frac{(1-j+2k)\Gamma(1+j+2n)}{\Gamma(1+j-k+n)\Gamma(2+k+n)},
\]
\[
|\mathcal{L}_n^{k,j}|=\frac{(2-j+2k)\Gamma(2+j+2n)}{\Gamma(1+j-k+n)\Gamma(3+k+n)}.
\]
\end{Lemma}

\begin{proof}
We prove this for the right cells as the left cells are similar. Let $p=k+1$, $q=k-j$, 
\[
F(n,p,q) := \sum_{r=0}^{q}(-1)^r \binom{q}{r}\frac{p-r}{2n+p-r}\binom{2n+p-r}{n},
\]
and 
\[
B(n,p) := \frac{p}{2n+p}\binom{2n+p}{n} = \binom{2n+p-1}{n}-\binom{2n+p-1}{n-1}.
\]
Straightforward binomial simplification gives $B(n,p)-B(n,p-1)=B(n-1,p+1)$. Define the backwards difference operator $\nabla B(n,p) := B(n,p)-B(n,p-1)$, then clearly the iterated operator $\nabla^q B(n,p) = B(n-q,p+q)$ (a simple induction shows this). Observe also by the binomial theorem 
\[
\nabla^q B(n,p) = \nabla^{q-1}(B(n,p)-B(n,p-1)) 
\] 
\[
= \sum_{r=0}^q (-1)^r \binom{q}{r} B(n,p-r)
\]
which is precisely $F(n,p,q)$. Finally, after some simplification,
\[
F(n,p,q) = {B(n-q,p+q) = \frac{p+q}{2n+p-q}\binom{2n+p-q}{n-q}} = \frac{(1-j+2k)\Gamma(1+j+2n)}{\Gamma(1+j-k+n)\Gamma(2+k+n)}.
\]
The proof is complete.
\end{proof}

\begin{figure}[ht]
\begin{gather*}
\includegraphics[scale=0.65]{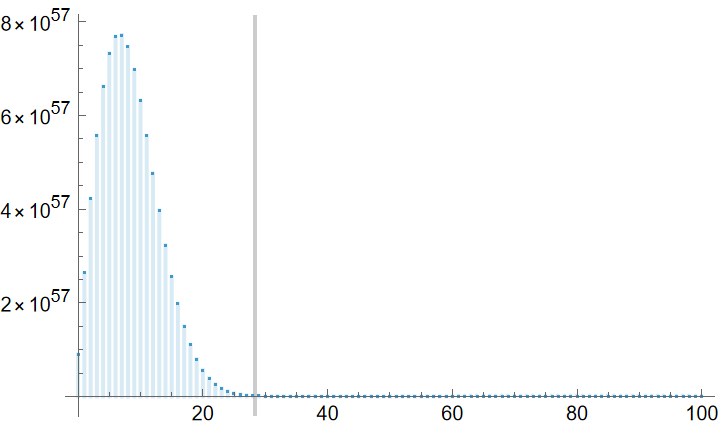}
\includegraphics[scale=0.65]{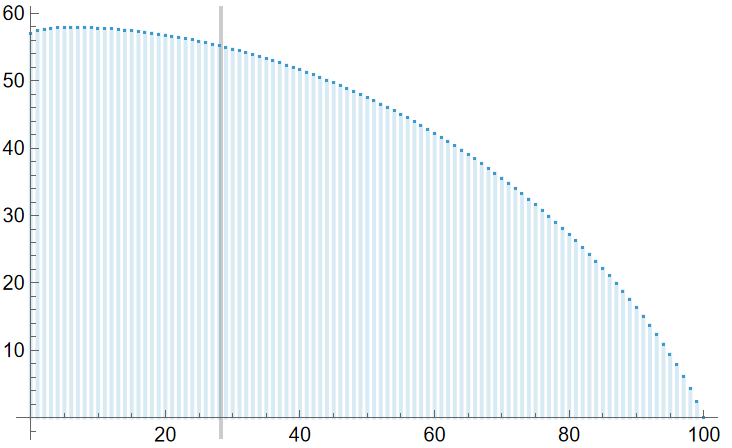}
\end{gather*}
\caption{Left: the semisimple dimension over $k$, for $n=100$, with a vertical line showing the rightmost truncation endpoint. Right: A Log10 version of the same graph.}
\label{rMo_truncation}
\end{figure}

We truncate $rMo_n$ at $0\le k\le b_n=\lfloor2\sqrt{2n}\rfloor$ to retain the range with large right cell dimensions. This is wider than the classical Motzkin window. See \autoref{rMo_truncation}.

We again denote this truncated monoid by $rMo_n^T$.

\begin{Theorem}
We have
\[
\textnormal{ssGap}_{\mathbb K}(rMo_n^T)=\Theta(4^n n^{-3/2}),\qquad
\textnormal{ssRatio}_{\mathbb K}(rMo_n^T)=\Theta(n^{-1/2}3^{-b_n/2}).
\]
\end{Theorem}

\begin{proof} 
By the formulas of \autoref{C:rmocellsizes}, both cell sizes have order $(k+1)2^j4^n n^{-3/2}$ uniformly for $k\le b_n$. Indeed, their factorial expressions are a factor of order $(k+1)/n$ times a binomial coefficient with upper argument $2n+j+O(1)$ and lower argument within $O(\sqrt n)$ of the centre. The rank-zero cell affords the trivial simple, so the lower bound from Stirling's formula and the cell at $k=1$, $j=0$ give the semisimple gap.

Now, 
\[
|rMo_n^T| = {1+}\sum_{k=0}^{{b_n}} \sum_{j=0}^{\lfloor \frac{k}{2}\rfloor} \dbinom{k+1}{2j+1}\frac{(1-j+2k)\Gamma(1+j+2n)}{\Gamma(1+j-k+n)\Gamma(2+k+n)} \frac{(2-j+2k)\Gamma(2+j+2n)}{\Gamma(1+j-k+n)\Gamma(3+k+n)}
\]
where the extra term accounts for the identity for large $n$. With the uniform cell estimates above, the inner sum is of order $(k+1)^2 16^n n^{-3}\sum_j\binom{k+1}{2j+1}4^j$. The binomial sum equals $(3^{k+1}-(-1)^{k+1})/4$. Summing over $k$ gives $\Theta(16^n n^{-2}3^{b_n})$, since the last term bounds the sum below and a geometric tail bounds it above. Dividing the semisimple gap by the square root of this order gives the semisimple gap ratio.
\end{proof}

\begin{Remark}
The asymptotic of the $n$th root of $|rMo_n^T|$ is $16$, however, unlike for $rTL$ and the pivotal planar monoids, this is significantly smaller than the size of the full monoid. Horizontal concatenation embeds $rMo_m\times rMo_n$ in $rMo_{m+n}$, so by Fekete's lemma the limit of $|rMo_n|^{1/n}$ is at least $|rMo_{500}|^{1/500}$, which the order formula gives numerically as larger than $20$. 

The reason for this discrepancy comes from the peculiar cell structure of $rMo_n$, which results in a large number of small cells. For $R_n^{k,0}$, $n$ fixed, the peak over $k$ is significantly different to the peak over $k$ of the monoid size $|rMo_n^k| = \sum_{j=0}^{\left\lfloor\frac{k}{2}\right\rfloor}\binom{k+1}{2j+1}|\mathcal{R}_n^{k,j}||\mathcal{L}_n^{k,j}|$. This is visualized in \autoref{rMo_bulk}.
\begin{figure}[ht]
\begin{gather*}
\includegraphics[scale=0.65]{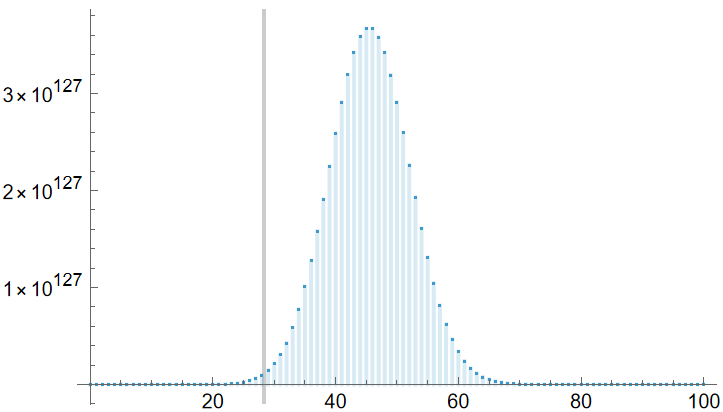}
\includegraphics[scale=0.65]{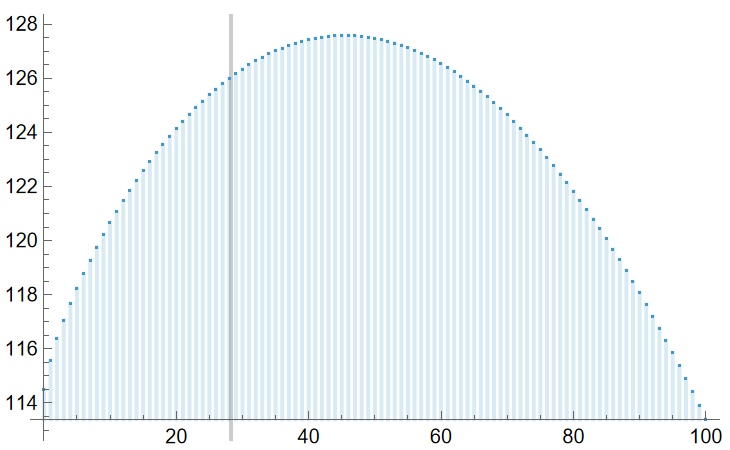}
\end{gather*}
\caption{Left: the bulk of the monoid $rMo_n$, with the rightmost truncation endpoint marked by a vertical line, showing the truncation is much smaller than the monoid size. Right: A Log10 version of the same graph.}
\label{rMo_bulk}
\end{figure}
Due to this, the size of the truncated monoid is relatively small, since we need to truncate where the peak of the right cells are. 

This also reveals that the limit superior of the $n$th root of the ratio
of the truncated RepGap to the square root of the size of the full monoid
is less than $1$. Therefore, the smallest nontrivial simple representations
of the truncation are significantly smaller than this upper bound, in the
context of \autoref{R:GapRatio}.
\end{Remark}

Combining this with \autoref{T:PivotalGaps}, we have the following conclusion:

\begin{Corollary}\label{C:MO}
Assume $\textnormal{char}\,\mathbb K=0$. Then $\textnormal{ssGap}_{\mathbb K}(rMo_n^T)^{1/n}\to4$ and $\textnormal{Gap}_{\mathbb K}(Mo_n^T)^{1/n}\to9$, so the rigid gap is exponentially smaller. The $n$th roots of $\textnormal{ssRatio}_{\mathbb K}(rMo_n^T)$ and $\textnormal{Ratio}_{\mathbb K}(Mo_n^T)$ both tend to $1$.
Thus, we get the table entries as in the introduction.
\end{Corollary}

Therefore, $rMo_n$ is significantly less suitable for cryptographic purposes than its pivotal counterpart, since the semisimple gap is much worse, and the semisimple gap ratio is at best the same, up to the exponential factor. Since the semisimple dimension is an upper bound for the simple dimension, and in light of \autoref{R:SSRepGapWorse}, we expect the RepGap and gap ratio to be even worse. 

\begin{figure}[ht]
\includegraphics[scale=0.75]{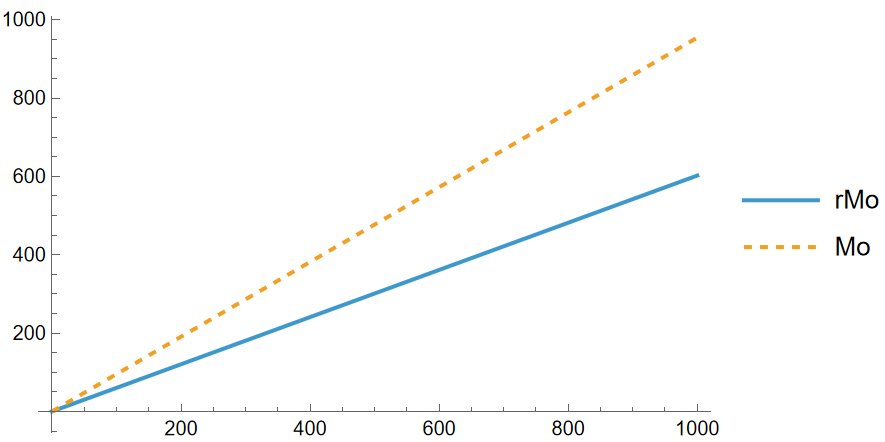}
\caption{Comparison of the exponential reference curves $4^n$ for the rigid semisimple gap and $9^n$ for the classical gap, on a base-ten logarithmic scale.}
\end{figure}

\begin{Remark}
As before, the characteristic assumption in \autoref{C:MO} is needed for $Mo_n$ itself. Indeed, the Motzkin category is Schur--Weyl dual to $\mathfrak{sl}_2$ \cite{BH,DG}, similarly to the Temperley--Lieb category but with a different generating representation (namely the vector plus the trivial representation), so the characteristic-zero proviso here is the same as for the classical Temperley--Lieb comparison. One could use \cite{Ar-RepGapMotzkin} to get the same result for arbitrary fields.
\end{Remark}

\section{The ``rigid'' planar rook monoid}

The planar rook monoid and its associated category have been discovered many times, see, e.g., \cite{So,KhSa-cat-polyring}, and we will now define a non-pivotal version of it.

\subsection{Definition of \texorpdfstring{$rpRo_n$}{rpRon}}

The planar rook monoid has no duals so it is naturally not possible to define a non-pivotal analog. We instead define it as follows, which is more of a two-color variant of the planar rook monoid:

\begin{Definition}
The \emph{``rigid'' planar rook category} (we drop the quotation marks from now on) $rpRo$ is the subcategory of $rMo$ that contains only through strands and dots. As before,
\begin{gather*}
rpRo_n := \textnormal{End}_{rpRo}(1\hspace{0.05cm}2)^{\otimes n},
\end{gather*}
is the \emph{rigid planar rook monoid}.
\end{Definition}

\begin{Remark}
Like the planar rook monoid, the planar symmetric monoid (the trivial group) also has no rigid analog, however the two-color submonoid of $rTL_n$ consisting of only through strands is isomorphic to the usual version, so we omit it from the discussion. 
\end{Remark}

\begin{Remark}
Inspired by the rigid planar rook monoid, it is potentially worth investigating two-color or $n$-color versions of the other diagram monoids, as somewhat of an intermediate between pivotal and rigid. 
\end{Remark}

\subsection{Computing cells for \texorpdfstring{$rpRo_n$}{rpRon}}

The arguments for the rigid planar rook monoid, $rpRo_n$, are very similar to those used for $rMo_n$. 

We find the sizes of the left and right cells of $rpRo_n$ via the same partition method used for $rMo_n$, however, it is much easier since each block has only one possible diagram (all dots). Moreover, there is no longer a difference between the left and right cells. 

\begin{Proposition}\label{P:rprocellsizes}
For each $n$, number of through strands $k$, and number of occurrences $j$ of $1 \hspace{0.05cm}2$ in the sequence of through strands $\alpha(k)$, the sizes of the left and right cells are:
\begin{equation*}
|L_n^{k,j}| = |R_n^{k,j}| = \binom{n+j}{k}.
\end{equation*}
Furthermore, we have 
\[
|rpRo_n| = \sum_{k=0}^{2n}\sum_{j=0}^{\left\lfloor \frac{k}{2} \right\rfloor} \binom{k+1}{2j+1}\binom{n+j}{k}^2.
\]
\end{Proposition}

\begin{proof}
Since each block only has one possible diagram, the count simply becomes the number of partitions $j_1 + \dots + j_{k+1} = 2n-k$. Using the results from $rMo_n$, we know that $j_1$, $j_{k+1}$, and $j_m$ for $m=2,3,\dots 2j$ are all even, and the remaining $k-2j$ are odd. Therefore, we can write $j_i = 2y_i$ for $i = 1,2,\dots, 2j,k+1$, and $j_i = 2y_i+1$ for $i = 2j+1, \dots k$, where each $y_i \geq 0$. So, the partition becomes:
\[
2y_1 + 2y_2 + \dots 2y_{2j} + 2y_{2j+1}+1 + \dots 2y_k + 1 + 2y_{k+1} = 2n-k
\implies y_1 + \dots y_{k+1} = n-k+j.
\]
Hence, $|\mathcal{L}_n^{k,j}|=|\mathcal{R}_n^{k,j}|$ is just the number of weak compositions of $n-k+j$ into $k+1$ parts, which is, by, for example, \cite[Theorem 5.2]{Bo-walk-combinatorics}:
\[
\binom{n-k+j+k+1-1}{k+1-1} = \binom{n+j}{k}.
\]
The final equation for $|rpRo_n|$ follows from the exact same argument used for $|rMo_n|$. 
\end{proof}

\begin{figure}[ht]
\includegraphics[scale=0.85]{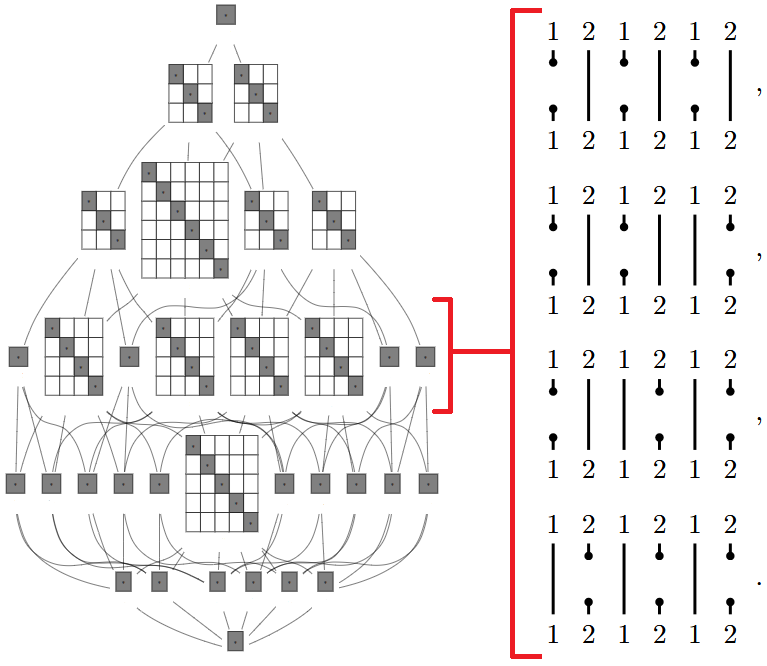}
\caption{Left: The cell structure for $rpRo_3$. Right: explicit diagrams representing the single elements of the four $J$-cells of size $1$ for $k=3$.}
\label{rpRo6_cells}
\end{figure}

\begin{Remark}
Once again, like in $rMo_n$, the sizes of the right cells (and left cells) increase with the number of $1\hspace{0.05cm}2$ blocks, meaning we need to focus on the $j=0$ case for the RepGap. Moreover, this proves the fact used previously in \autoref{rmojincrease} that the number of different partitions increases as $j$ increases. As seen in \autoref{rpRo6_cells}, we have a similar cell structure to $rMo$.
\end{Remark}

\subsection{RepGap and gap ratio of \texorpdfstring{$rpRo_n$}{rpRon}}

The key difference between $rpRo_n$ and $rMo_n$ is that $rpRo_n$ is semisimple, and this follows from the following.

\begin{Proposition}\label{P:rpRo_simples}
The right representations $\Delta _{R_n^{\alpha(k)}}$ are all simple, and form a complete list of the simple right representations of $rpRo_n$, up to isomorphism. 
\end{Proposition}

\begin{proof}
For a fixed through word, the same occupied subsets index both halves. Gluing preserves the word precisely when the two subsets coincide, so the sandwich matrix is the identity and has full rank over every field. This is the argument of \cite[Proposition 4F.7]{khovanov-monoidal-2024}.
\end{proof}

\begin{Remark}
The above proposition applies to left representations as well, as the arguments used apply similarly to multiplication on the left. 
\end{Remark}

\begin{Lemma}
$rpRo_n$ is semisimple over (arbitrary) $\mathbb{K}$.
\end{Lemma}

\begin{proof}
Every $J$-cell is idempotent by the same argument used in \autoref{P:rmocells} (in fact we construct the same idempotents) and every $J$-cell is square, by \autoref{P:rprocellsizes}. Finally, combining these facts with \autoref{P:rpRo_simples} and \autoref{P:CellsSemisimple} proves that $rpRo_n$ is semisimple over $\mathbb{K}$.
\end{proof}

\begin{figure}[ht]
\begin{gather*}
\includegraphics[scale=0.65]{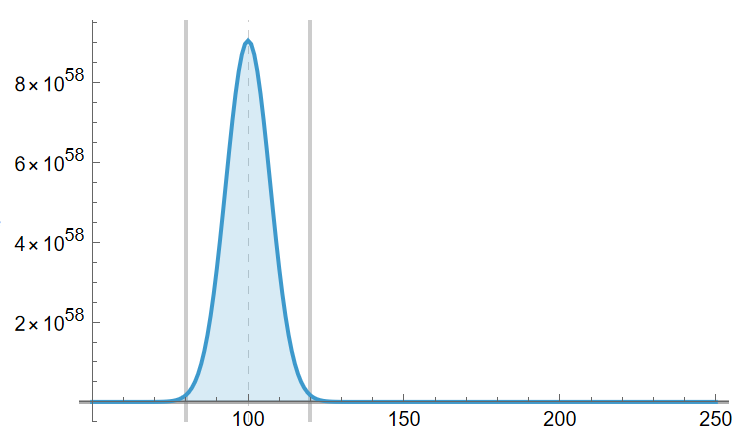}
\includegraphics[scale=0.6]{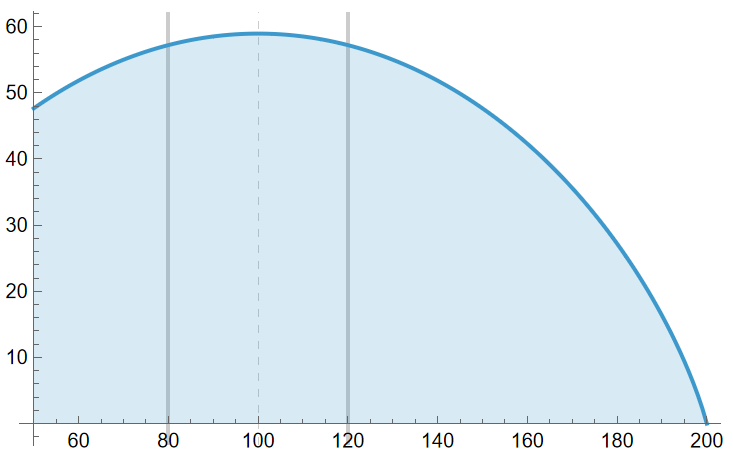}
\end{gather*}
\caption{Left: the dimension of the simple representations over $k$, for $n=200$, with vertical lines marking the truncation end points. Right: A Log10 version of the same graph.}
\label{rpRo_truncation}
\end{figure}

As with $pRo_n$, we truncate $rpRo_n^T := rpRo_n^{\frac{n}{2}-\sqrt{2n} \leq k \leq \frac{n}{2}+\sqrt{2n}}$. See \autoref{rpRo_truncation}.

\begin{Theorem} 
Over every field $\mathbb K$,
\[
\textnormal{Gap}_{\mathbb K}(rpRo_n^T)=\Theta(2^n n^{-1/2}),\qquad
\textnormal{Ratio}_{\mathbb K}(rpRo_n^T)=O(\rho^n),\qquad
\rho=2^{3/2}3^{3/4}5^{-5/4}<0.87 .
\]
\end{Theorem}

\begin{proof}
As with $rMo$, we take $j=0$ to get the minimal dimension of the representations.

Uniformly over the retained window, the smallest dimension at each rank is $\binom nk=\Theta(2^n n^{-1/2})$ by Stirling's formula.

For the ratio take $k=\lfloor n/2\rfloor$ and $j=\lfloor k/2\rfloor$. The corresponding summand in the order formula gives $\sqrt{|rpRo_n^T|}\ge\binom{n+j}{k}$. Stirling's formula gives $\binom nk/\binom{n+j}{k}=O(\rho^n)$, since rounding contributes bounded factors.
\end{proof}

\begin{Remark}
As with $rMo$, the truncated monoid is much smaller than the full monoid, since the location of the right cell peak is different to the location of the monoid peak, over $k$ through strands, seen in \autoref{rpRo_bulk}.
\begin{figure}[ht]
\begin{gather*}
\includegraphics[scale=0.65]{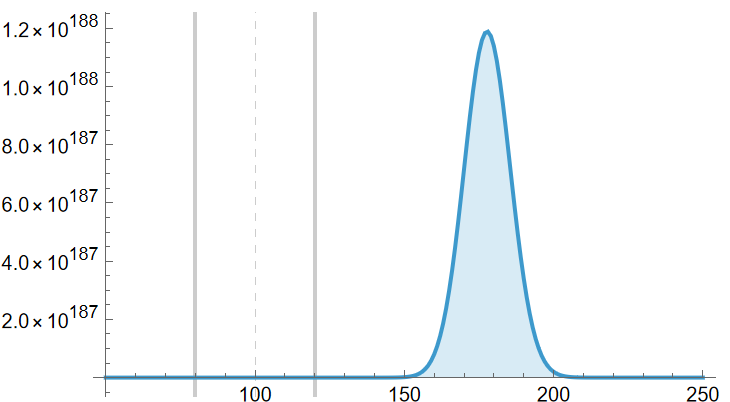}
\includegraphics[scale=0.6]{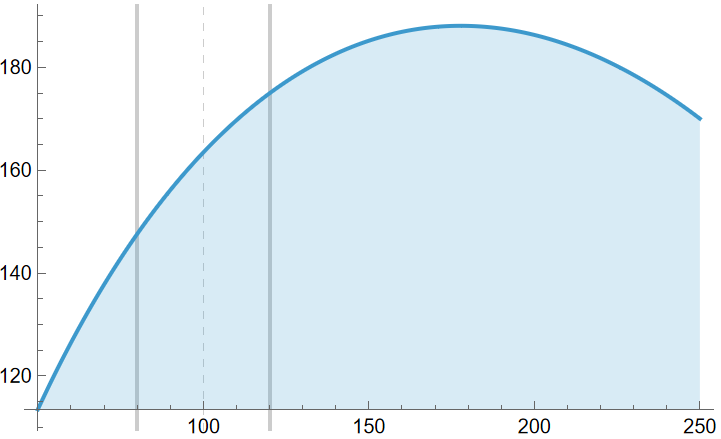}
\end{gather*}
\caption{Left: the bulk of the monoid $rpRo_n$, with truncation endpoints marked by vertical lines, showing the truncation is much smaller than the monoid size. Right: A Log10 version of the same graph.}
\label{rpRo_bulk}
\end{figure}
We can again find a lower bound for the asymptotic of $|rpRo_n|^{1/n}$ using the same horizontal concatenation argument as for $rMo_n$, and evaluating the order formula at $n=500$. We therefore have $|rpRo_n| \geq \Omega(8.9^n)$, meaning we once again have the $n$th root of the ratio of the RepGap to the square root of the size of the monoid being less than 1. 
\end{Remark}

\begin{Remark}
The upper bound for the gap ratio is quite coarse, since we are picking only one (non-maximal) term in the double sum that makes up $|rpRo_n^T|$. As a result, $rpRo_n$ is significantly worse for cryptographic purposes than $pRo_n$. 
\end{Remark}

\begin{Corollary}
The ratio of the two-colour gap to the classical gap is $\Theta(2^{-n})$. The ratio of their gap ratios is $O(n^{1/4}\rho^n)$.
Thus, we get the table entries as in the introduction.
\end{Corollary}

\begin{figure}[ht]
\includegraphics[scale=0.7]{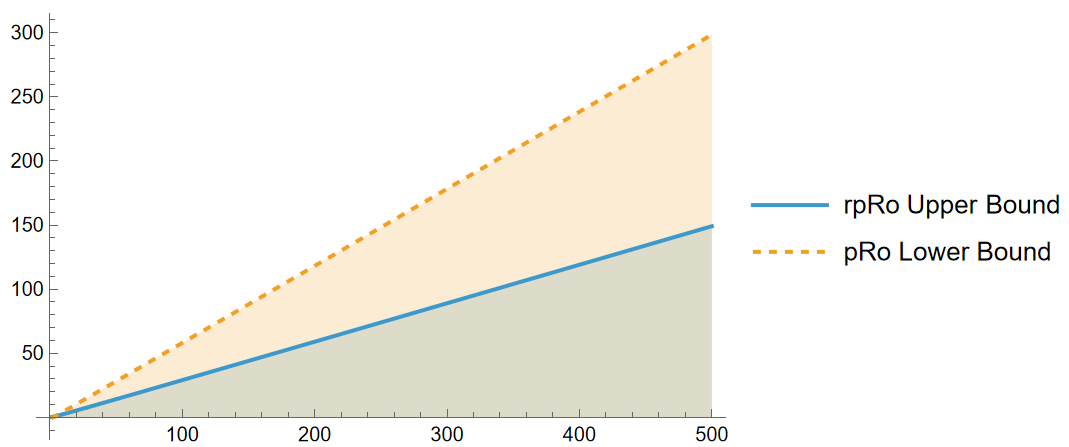}
\caption{Visual comparison of the RepGaps between the non-pivotal and pivotal planar rook monoid, on a Log10 plot.}
\end{figure}

\begin{figure}[ht]
\includegraphics[scale=0.7]{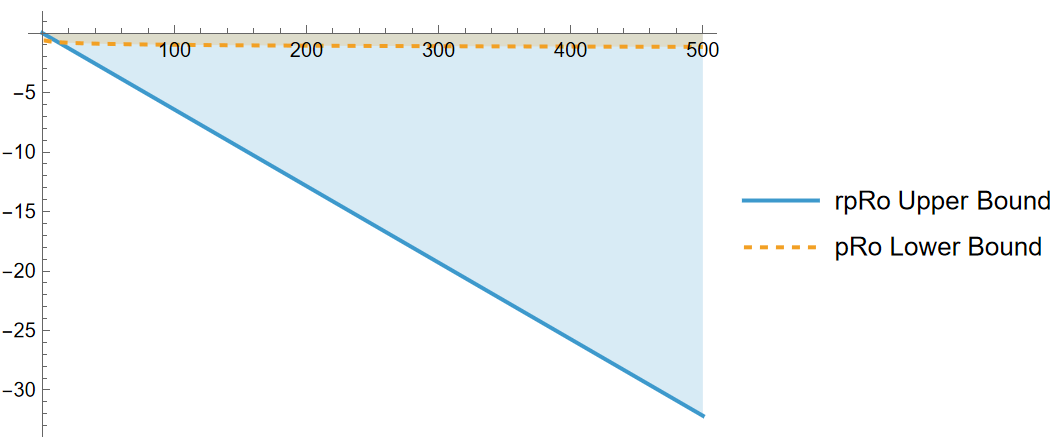}
\caption{Visual comparison of the gap ratios between the non-pivotal and pivotal planar rook monoid, on a Log10 plot.}
\end{figure}


\newpage
\newcommand{\etalchar}[1]{$^{#1}$}

\end{document}